# The Askey-scheme of hypergeometric orthogonal polynomials and its $q$-analogue

Roelof Koekoek     René F. Swarttouw

February 20, 1996

**Abstract**

We list the so-called Askey-scheme of hypergeometric orthogonal polynomials. In chapter 1 we give the definition, the orthogonality relation, the three term recurrence relation and generating functions of all classes of orthogonal polynomials in this scheme. In chapter 2 we give all limit relations between different classes of orthogonal polynomials listed in the Askey-scheme.

In chapter 3 we list the $q$-analogues of the polynomials in the Askey-scheme. We give their definition, orthogonality relation, three term recurrence relation and generating functions. In chapter 4 we give the limit relations between those basic hypergeometric orthogonal polynomials. Finally, in chapter 5 we point out how the 'classical' hypergeometric orthogonal polynomials of the Askey-scheme can be obtained from their $q$-analogues.

**Acknowledgement**

We would like to thank Professor Tom H. Koornwinder who suggested us to write a report like this. He also helped us solving many problems we encountered during the research and provided us with several references.

# Contents

# Preface

This report deals with orthogonal polynomials appearing in the so-called Askey-scheme of hypergeometric orthogonal polynomials and their $q$-analogues. Most formulas listed in this report can be found somewhere in the literature, but a handbook containing all these formulas did not exist. We collected known formulas for these hypergeometric orthogonal polynomials and we arranged them into the Askey-scheme and into a $q$-analogue of this scheme which we called the $q$-scheme. This $q$-scheme was not completely documented in the literature. So we filled in some gaps in order to get some sort of 'complete' scheme of $q$-hypergeometric orthogonal polynomials.

In chapter 0 we give some general definitions and formulas which can be used to transform several formulas into different forms of the same formula. In the other chapters we used the most common notations, but sometimes we had to change some notations in order to be consistent.

For each family of orthogonal polynomials listed in this report we give the conditions on the parameters for which the corresponding weight function is positive. These conditions are mentioned in the orthogonality relations. We remark that many of these orthogonal polynomials are still polynomials for other values of the parameters and that they can be defined for other values as well. That is why we gave no restrictions in the definitions. As pointed out in chapter 0 some definitions can be transformed into different forms so that they are valid for some values of the parameters for which the given form has no meaning. Other formulas, such as the generating functions, are only valid for some special values of parameters and arguments. These conditions are left out in this report.

We are aware of the fact that this report is by no means a full description of all that is known about (basic) hypergeometric orthogonal polynomials. More on each listed family of orthogonal polynomials can be found in the articles and books to which we refer.

In later versions of this report we want to add recurrence relations for the monic polynomials in each case, id est for the polynomials with leading coefficient equal to 1. We also hope to include more formulas containing quadratic transformations and we want to pay more attention to the transformation $q \leftrightarrow q^{-1}$.

Comments on this version of the report and suggestions for improvement are most welcome. If you find errors or gaps or if you have suggestions for inclusion of more formulas on (basic) hypergeometric orthogonal polynomials, please contact us and let us know.

<div align="right">Roelof Koekoek and René F. Swarttouw.</div>

| | |
|---|---|
| Roelof Koekoek | René F. Swarttouw |
| Delft University of Technology | Free University of Amsterdam |
| Faculty of Technical Mathematics and Informatics | Faculty of Mathematics and Informatics |
| Mekelweg 4 | De Boelelaan 1081 |
| 2628 CD Delft | 1081 HV Amsterdam |
| The Netherlands | The Netherlands |
| koekoek@twi.tudelft.nl | rene@cs.vu.nl |

# Definitions and miscellaneous formulas

## 0.1 Introduction

In this report we will list all known sets of orthogonal polynomials which can be defined in terms of a hypergeometric function or a basic hypergeometric function.

In the first part of the report we give a description of all classical hypergeometric orthogonal polynomials which appear in the so-called Askey-scheme. We give definitions, orthogonality relations, three term recurrence relations, differential or difference equations and generating functions for all families of orthogonal polynomials listed in this Askey-scheme of hypergeometric orthogonal polynomials.

In the second part we obtain a $q$-analogue of this scheme. We give definitions, orthogonality relations, three term recurrence relations, difference equations and generating functions for all known $q$-analogues of the hypergeometric orthogonal polynomials listed in the Askey-scheme.

Further we give limit relations between different families of orthogonal polynomials in both schemes and we point out how to obtain the classical hypergeometric orthogonal polynomials from their $q$-analogues.

The theory of $q$-analogues or $q$-extensions of classical formulas and functions is based on the observation that
$$\lim_{q \to 1} \frac{1-q^\alpha}{1-q} = \alpha.$$
Therefore the number $(1-q^\alpha)/(1-q)$ is sometimes called the basic number $[\alpha]$.

Now we can give a $q$-analogue of the Pochhammer-symbol $(a)_k$ which is defined by
$$(a)_0 := 1 \text{ and } (a)_k := a(a+1)(a+2)\cdots(a+k-1), \ k = 1, 2, 3, \ldots.$$
This $q$-extension is given by
$$(a;q)_0 := 1 \text{ and } (a;q)_k := (1-a)(1-aq)(1-aq^2)\cdots(1-aq^{k-1}), \ k = 1, 2, 3, \ldots.$$
It is clear that
$$\lim_{q \to 1} \frac{(q^\alpha;q)_k}{(1-q)^k} = (\alpha)_k.$$

In this report we will always assume that $0 < q < 1$.

For more details concerning the $q$-theory the reader is referred to the book [114] by G. Gasper and M. Rahman.

Since many formulas given in this report can be reformulated in many different ways we will give a selection of formulas, which can be used to obtain other forms of definitions, orthogonality relations and generating functions.

Most of these formulas given in this chapter can be found in [114].

We remark that in orthogonality relations we often have to add some condition(s) on the parameters of the orthogonal polynomials involved in order to have positive weight functions. By using the famous theorem of Favard these conditions can also be obtained from the three term recurrence relation.

In some cases, however, some conditions on the parameters may be needed in other formulas too. For instance, the definition (1.11.1) of the Laguerre polynomials has no meaning for negative integer values of the parameter $\alpha$. But in fact the Laguerre polynomials are also polynomials in the parameter $\alpha$. This can be seen by writing

$$L_n^{(\alpha)}(x) = \frac{1}{n!} \sum_{k=0}^{n} \frac{(-n)_k}{k!} (\alpha+k+1)_{n-k} x^k.$$

In this way the Laguerre polynomials are defined for all values of the parameter $\alpha$.

A similar remark holds for the Jacobi polynomials given by (1.8.1). We may also write (see section 0.4 for the definition of the hypergeometric function $_2F_1$)

$$P_n^{(\alpha,\beta)}(x) = (-1)^n \frac{(\beta+1)_n}{n!} {}_2F_1\left( \begin{array}{c} -n, n+\alpha+\beta+1 \\ \beta+1 \end{array} \middle| \frac{1+x}{2} \right),$$

which implies the well-known symmetry relation

$$P_n^{(\alpha,\beta)}(x) = (-1)^n P_n^{(\beta,\alpha)}(-x).$$

Even more general we have

$$P_n^{(\alpha,\beta)}(x) = \frac{1}{n!} \sum_{k=0}^{n} \frac{(-n)_k}{k!} (n+\alpha+\beta+1)_k (\alpha+k+1)_{n-k} \left(\frac{1-x}{2}\right)^k.$$

From this form it is clear that the Jacobi polynomials can be defined for all values of the parameters $\alpha$ and $\beta$ although the definition (1.8.1) is not valid for negative integer values of the parameter $\alpha$.

We will not indicate these difficulties in each formula.

Finally, we remark that in each recurrence relation listed in this report, except for (1.8.24) for the Chebyshev polynomials of the first kind, we may use $P_{-1}(x) = 0$ and $P_0(x) = 1$ as initial conditions.

## 0.2 The $q$-shifted factorials

The symbols $(a;q)_k$ defined in the preceding section are called $q$-shifted factorials. They can also be defined for negative values of $k$ as

$$(a;q)_k := \frac{1}{(1-aq^{-1})(1-aq^{-2})\cdots(1-aq^k)}, \ a \neq q, q^2, q^3, \ldots, q^{-k}, \ k = -1,-2,-3,\ldots. \quad (0.2.1)$$

Now we have

$$(a;q)_{-n} = \frac{1}{(aq^{-n};q)_n} = \frac{(-qa^{-1})^n}{(qa^{-1};q)_n} q^{\binom{n}{2}}, \ n = 0,1,2,\ldots, \quad (0.2.2)$$

where

$$\binom{n}{2} = \frac{1}{2}n(n-1).$$

We can also define

$$(a;q)_\infty = \prod_{k=0}^{\infty} (1 - aq^k).$$

This implies that

$$(a;q)_n = \frac{(a;q)_\infty}{(aq^n;q)_\infty}, \quad (0.2.3)$$

and, for any complex number $\lambda$,
$$(a;q)_\lambda = \frac{(a;q)_\infty}{(aq^\lambda;q)_\infty}, \tag{0.2.4}$$
where the principal value of $q^\lambda$ is taken.

If we change $q$ by $q^{-1}$ we obtain
$$(a;q^{-1})_n = (a^{-1};q)_n(-a)^n q^{-\binom{n}{2}},\ a \neq 0. \tag{0.2.5}$$

This formula can be used, for instance, to prove the following transformation formula between the little $q$-Laguerre (or Wall) polynomials given by (3.20.1) and the $q$-Laguerre polynomials defined by (3.21.1) :
$$p_n(x;q^{-\alpha}|q^{-1}) = \frac{(q;q)_n}{(q^{\alpha+1};q)_n} L_n^{(\alpha)}(-x;q)$$
or equivalently
$$L_n^{(\alpha)}(x;q^{-1}) = \frac{(q^{\alpha+1};q)_n}{(q;q)_n q^{n\alpha}} p_n(-x;q^\alpha|q).$$

By using (0.2.5) it is not very difficult to verify the following general transformation formula for $_4\phi_3$ polynomials (see section 0.4 for the definition of the basic hypergeometric function $_4\phi_3$) :
$$_4\phi_3\left(\begin{matrix}q^n,a,b,c\\d,e,f\end{matrix}\bigg|q^{-1};q^{-1}\right) = {}_4\phi_3\left(\begin{matrix}q^{-n},a^{-1},b^{-1},c^{-1}\\d^{-1},e^{-1},f^{-1}\end{matrix}\bigg|q;\frac{abcq^n}{def}\right),$$
where a limit is needed when one of the parameters is equal to zero. Other transformation formulas can be obtained from this one by applying limits as discussed in section 0.4.

Finally, we list a number of transformation formulas for the $q$-shifted factorials, where $k$ and $n$ are nonnegative integers :
$$(a;q)_{n+k} = (a;q)_n (aq^n;q)_k. \tag{0.2.6}$$

$$\frac{(aq^n;q)_k}{(aq^k;q)_n} = \frac{(a;q)_k}{(a;q)_n}. \tag{0.2.7}$$

$$(aq^k;q)_{n-k} = \frac{(a;q)_n}{(a;q)_k},\ k=0,1,2,\ldots,n. \tag{0.2.8}$$

$$(a;q)_n = (a^{-1}q^{1-n};q)_n(-a)^n q^{\binom{n}{2}},\ a \neq 0. \tag{0.2.9}$$

$$(aq^{-n};q)_n = (a^{-1}q;q)_n(-a)^n q^{-n-\binom{n}{2}},\ a \neq 0. \tag{0.2.10}$$

$$\frac{(aq^{-n};q)_n}{(bq^{-n};q)_n} = \frac{(a^{-1}q;q)_n}{(b^{-1}q;q)_n}\left(\frac{a}{b}\right)^n,\ a \neq 0,\ b \neq 0. \tag{0.2.11}$$

$$(a;q)_{n-k} = \frac{(a;q)_n}{(a^{-1}q^{1-n};q)_k}\left(-\frac{q}{a}\right)^k q^{\binom{k}{2}-nk},\ a \neq 0,\ k=0,1,2,\ldots,n. \tag{0.2.12}$$

$$\frac{(a;q)_{n-k}}{(b;q)_{n-k}} = \frac{(a;q)_n}{(b;q)_n}\frac{(b^{-1}q^{1-n};q)_k}{(a^{-1}q^{1-n};q)_k}\left(\frac{b}{a}\right)^k,\ a \neq 0,\ b \neq 0,\ k=0,1,2,\ldots,n. \tag{0.2.13}$$

$$(q^{-n};q)_k = \frac{(q;q)_n}{(q;q)_{n-k}}(-1)^k q^{\binom{k}{2}-nk},\ k=0,1,2,\ldots,n. \tag{0.2.14}$$

$$(aq^{-n};q)_k = \frac{(a^{-1}q;q)_n}{(a^{-1}q^{1-k};q)_n}(a;q)_k q^{-nk},\ a \neq 0. \tag{0.2.15}$$

$$(aq^{-n};q)_{n-k} = \frac{(a^{-1}q;q)_n}{(a^{-1}q;q)_k}\left(-\frac{a}{q}\right)^{n-k} q^{\binom{k}{2}-\binom{n}{2}},\ a \neq 0,\ k=0,1,2,\ldots,n. \tag{0.2.16}$$

$$(a;q)_{2n} = (a;q^2)_n(aq;q^2)_n. \tag{0.2.17}$$

$$(a^2;q^2)_n = (a;q)_n(-a;q)_n. \qquad (0.2.18)$$

$$(a;q)_\infty = (a;q^2)_\infty (aq;q^2)_\infty. \qquad (0.2.19)$$

$$(a^2;q^2)_\infty = (a;q)_\infty (-a;q)_\infty. \qquad (0.2.20)$$

Formula (0.2.18) can be used, for instance, to show that the generating function (3.10.29) for the continuous $q$-Legendre polynomials is a $q$-analogue of the generating function (1.8.44) for the Legendre polynomials. In fact we obtain (see section 0.4 for the definition of the basic hypergeometric functions $_r\phi_s$)

$$_2\phi_1\left(\begin{array}{c}q^{\frac{1}{2}}e^{i\theta},-q^{\frac{1}{2}}e^{i\theta}\\-q\end{array}\bigg|q;e^{-i\theta}t\right) {}_2\phi_1\left(\begin{array}{c}q^{\frac{1}{2}}e^{-i\theta},-q^{\frac{1}{2}}e^{-i\theta}\\-q\end{array}\bigg|q;e^{i\theta}t\right)$$
$$= {}_1\phi_0\left(\begin{array}{c}qe^{2i\theta}\\-\end{array}\bigg|q^2;e^{-i\theta}t\right) {}_1\phi_0\left(\begin{array}{c}qe^{-2i\theta}\\-\end{array}\bigg|q^2;e^{i\theta}t\right).$$

Now we use the $q$-binomial theorem (0.5.2) to show that this equals

$$\frac{(qe^{i\theta}t;q^2)_\infty(qe^{-i\theta}t;q^2)_\infty}{(e^{i\theta}t;q^2)_\infty(e^{-i\theta}t;q^2)_\infty} = {}_1\phi_0\left(\begin{array}{c}q\\-\end{array}\bigg|q^2;e^{i\theta}t\right) {}_1\phi_0\left(\begin{array}{c}q\\-\end{array}\bigg|q^2;e^{-i\theta}t\right).$$

If we let $q$ tend to one we now find by using the binomial theorem (0.5.1)

$$_1F_0\left(\begin{array}{c}\frac{1}{2}\\-\end{array}\bigg|e^{i\theta}t\right) {}_1F_0\left(\begin{array}{c}\frac{1}{2}\\-\end{array}\bigg|e^{-i\theta}t\right) = (1-e^{i\theta}t)^{-\frac{1}{2}}(1-e^{-i\theta}t)^{-\frac{1}{2}} = \frac{1}{\sqrt{1-2xt+t^2}},\ x=\cos\theta,$$

which equals (1.8.44).

## 0.3 The $q$-gamma function and the $q$-binomial coefficient

The $q$-gamma function is defined by

$$\Gamma_q(x) := \frac{(q;q)_\infty}{(q^x;q)_\infty}(1-q)^{1-x}. \qquad (0.3.1)$$

This is a $q$-analogue of the well-known gamma function since we have

$$\lim_{q\uparrow 1}\Gamma_q(x) = \Gamma(x).$$

Note that the $q$-gamma function satisfies the functional equation

$$\Gamma_q(z+1) = \frac{1-q^z}{1-q}\Gamma_q(z),\ \Gamma_q(1)=1,$$

which is a $q$-extension of the well-known functional equation

$$\Gamma(z+1) = z\Gamma(z),\ \Gamma(1)=1$$

for the ordinary gamma function. For nonintegral values of $z$ this ordinary gamma function also satisfies the relation

$$\Gamma(z)\Gamma(1-z) = \frac{\pi}{\sin\pi z},$$

which can be used to show that

$$\lim_{\alpha\to k}(1-q^{-\alpha+k})\Gamma(-\alpha)\Gamma(\alpha+1) = (-1)^{k+1}\ln q,\ k=0,1,2,\ldots.$$

This limit can be used to show that the orthogonality relation (3.27.2) for the Stieltjes-Wigert polynomials follows from the orthogonality relation (3.21.2) for the $q$-Laguerre polynomials.

The $q$-binomial coefficient is defined by

$$\begin{bmatrix} n \\ k \end{bmatrix}_q = \begin{bmatrix} n \\ n-k \end{bmatrix}_q := \frac{(q;q)_n}{(q;q)_k (q;q)_{n-k}}, \ k = 0, 1, 2, \ldots, n, \qquad (0.3.2)$$

where $n$ denotes a nonnegative integer.

This definition can be generalized in the following way. For arbitrary complex $\alpha$ we have

$$\begin{bmatrix} \alpha \\ k \end{bmatrix}_q := \frac{(q^{-\alpha};q)_k}{(q;q)_k}(-1)^k q^{k\alpha - \binom{k}{2}}. \qquad (0.3.3)$$

Or more general for all complex $\alpha$ and $\beta$ we have

$$\begin{bmatrix} \alpha \\ \beta \end{bmatrix}_q := \frac{\Gamma_q(\alpha+1)}{\Gamma_q(\beta+1)\Gamma_q(\alpha-\beta+1)} = \frac{(q^{\beta+1};q)_\infty (q^{\alpha-\beta+1};q)_\infty}{(q;q)_\infty (q^{\alpha+1};q)_\infty}. \qquad (0.3.4)$$

For instance this implies that

$$\frac{(q^{\alpha+1};q)_n}{(q;q)_n} = \begin{bmatrix} n+\alpha \\ n \end{bmatrix}_q.$$

Note that

$$\lim_{q \uparrow 1} \begin{bmatrix} \alpha \\ \beta \end{bmatrix}_q = \binom{\alpha}{\beta} = \frac{\Gamma(\alpha+1)}{\Gamma(\beta+1)\Gamma(\alpha-\beta+1)}.$$

For integer values of the parameter $\beta$ we have

$$\binom{\alpha}{k} = \frac{(-\alpha)_k}{k!}(-1)^k, \ k = 0, 1, 2, \ldots$$

and when the parameter $\alpha$ is an integer too we may write

$$\binom{n}{k} = \frac{n!}{k!(n-k)!}, \ k = 0, 1, 2, \ldots, n, \ n = 0, 1, 2, \ldots.$$

This latter formula can be used to show that

$$\binom{2n}{n} = \frac{\left(\frac{1}{2}\right)_n}{n!} 4^n, \ n = 0, 1, 2, \ldots.$$

This can be used to write the generating functions (1.8.29) and (1.8.35) for the Chebyshev polynomials of the first and the second kind in the following form :

$$R^{-1}\sqrt{\frac{1}{2}(1+R-xt)} = \sum_{n=0}^\infty \binom{2n}{n} T_n(x) \left(\frac{t}{4}\right)^n, \ R = \sqrt{1-2xt+t^2}$$

and

$$\frac{4}{R\sqrt{2(1+R-xt)}} = \sum_{n=0}^\infty \binom{2n+2}{n+1} U_n(x) \left(\frac{t}{4}\right)^n, \ R = \sqrt{1-2xt+t^2}$$

respectively.

Finally we remark that

$$(a;q)_n = \sum_{k=0}^n \begin{bmatrix} n \\ k \end{bmatrix}_q q^{\binom{k}{2}}(-a)^k. \qquad (0.3.5)$$

## 0.4 Hypergeometric and basic hypergeometric functions

The hypergeometric series $_rF_s$ is defined by

$$_rF_s\left(\begin{array}{c}a_1,\ldots,a_r\\b_1,\ldots,b_s\end{array}\bigg|\,z\right):=\sum_{k=0}^{\infty}\frac{(a_1,\ldots,a_r)_k}{(b_1,\ldots,b_s)_k}\frac{z^k}{k!},\qquad(0.4.1)$$

where

$$(a_1,\ldots,a_r)_k:=(a_1)_k\cdots(a_r)_k.$$

Of course, the parameters must be such that the denominator factors in the terms of the series are never zero. When one of the numerator parameters $a_i$ equals $-n$ where $n$ is a nonnegative integer this hypergeometric series is a polynomial in $z$. Otherwise the radius of convergence $\rho$ of the hypergeometric series is given by

$$\rho=\begin{cases}\infty & \text{if } r<s+1\\ 1 & \text{if } r=s+1\\ 0 & \text{if } r>s+1.\end{cases}$$

A hypergeometric series of the form (0.4.1) is called balanced (or Saalschützian) if $r=s+1$, $z=1$ and $a_1+a_2+\ldots+a_{s+1}+1=b_1+b_2+\ldots+b_s$.

The basic hypergeometric series (or $q$-hypergeometric series) $_r\phi_s$ is defined by

$$_r\phi_s\left(\begin{array}{c}a_1,\ldots,a_r\\b_1,\ldots,b_s\end{array}\bigg|\,q;z\right):=\sum_{k=0}^{\infty}\frac{(a_1,\ldots,a_r;q)_k}{(b_1,\ldots,b_s;q)_k}(-1)^{(1+s-r)k}q^{(1+s-r)\binom{k}{2}}\frac{z^k}{(q;q)_k},\qquad(0.4.2)$$

where

$$(a_1,\ldots,a_r;q)_k:=(a_1;q)_k\cdots(a_r;q)_k.$$

Again, we assume that the parameters are such that the denominator factors in the terms of the series are never zero. If one of the numerator parameters $a_i$ equals $q^{-n}$ where $n$ is a nonnegative integer this basic hypergeometric series is a polynomial in $z$. Otherwise the radius of convergence $\rho$ of the basic hypergeometric series is given by

$$\rho=\begin{cases}\infty & \text{if } r<s+1\\ 1 & \text{if } r=s+1\\ 0 & \text{if } r>s+1.\end{cases}$$

The special case $r=s+1$ reads

$$_{s+1}\phi_s\left(\begin{array}{c}a_1,\ldots,a_{s+1}\\b_1,\ldots,b_s\end{array}\bigg|\,q;z\right)=\sum_{k=0}^{\infty}\frac{(a_1,\ldots,a_{s+1};q)_k}{(b_1,\ldots,b_s;q)_k}\frac{z^k}{(q;q)_k}.$$

This basic hypergeometric series was first introduced by Heine in 1846. Therefore it is sometimes called Heine's series. A basic hypergeometric series of this form is called balanced (or Saalschützian) if $z=q$ and $a_1a_2\cdots a_{s+1}q=b_1b_2\cdots b_s$.

The $q$-hypergeometric series is a $q$-analogue of the hypergeometric series defined by (0.4.1) since

$$\lim_{q\uparrow 1}{}_r\phi_s\left(\begin{array}{c}q^{a_1},\ldots,q^{a_r}\\q^{b_1},\ldots,q^{b_s}\end{array}\bigg|\,q;(q-1)^{1+s-r}z\right)={}_rF_s\left(\begin{array}{c}a_1,\ldots,a_r\\b_1,\ldots,b_s\end{array}\bigg|\,z\right).$$

This limit will be used frequently in chapter 5.

We remark that

$$\lim_{a_r \to \infty} {}_r\phi_s \left( \begin{matrix} a_1, \ldots, a_r \\ b_1, \ldots, b_s \end{matrix} \bigg| q; \frac{z}{a_r} \right) = {}_{r-1}\phi_s \left( \begin{matrix} a_1, \ldots, a_{r-1} \\ b_1, \ldots, b_s \end{matrix} \bigg| q; z \right).$$

In fact, this is the reason for the factors $(-1)^{(1+s-r)k} q^{(1+s-r)\binom{k}{2}}$ in the definition (0.4.2) of the basic hypergeometric series.

The limit relations between hypergeometric orthogonal polynomials listed in chapter 2 of this report are based on the observations that

$$\,_rF_s \left( \begin{matrix} a_1, \ldots, a_{r-1}, \mu \\ b_1, \ldots, b_{s-1}, \mu \end{matrix} \bigg| z \right) = {}_{r-1}F_{s-1} \left( \begin{matrix} a_1, \ldots, a_{r-1} \\ b_1, \ldots, b_{s-1} \end{matrix} \bigg| z \right), \tag{0.4.3}$$

$$\lim_{\lambda \to \infty} {}_rF_s \left( \begin{matrix} a_1, \ldots, a_{r-1}, \lambda a_r \\ b_1, \ldots, b_s \end{matrix} \bigg| \frac{z}{\lambda} \right) = {}_{r-1}F_s \left( \begin{matrix} a_1, \ldots, a_{r-1} \\ b_1, \ldots, b_s \end{matrix} \bigg| a_r z \right), \tag{0.4.4}$$

$$\lim_{\lambda \to \infty} {}_rF_s \left( \begin{matrix} a_1, \ldots, a_r \\ b_1, \ldots, b_{s-1}, \lambda b_s \end{matrix} \bigg| \lambda z \right) = {}_rF_{s-1} \left( \begin{matrix} a_1, \ldots, a_r \\ b_1, \ldots, b_{s-1} \end{matrix} \bigg| \frac{z}{b_s} \right) \tag{0.4.5}$$

and

$$\lim_{\lambda \to \infty} {}_rF_s \left( \begin{matrix} a_1, \ldots, a_{r-1}, \lambda a_r \\ b_1, \ldots, b_{s-1}, \lambda b_s \end{matrix} \bigg| z \right) = {}_{r-1}F_{s-1} \left( \begin{matrix} a_1, \ldots, a_{r-1} \\ b_1, \ldots, b_{s-1} \end{matrix} \bigg| \frac{a_r z}{b_s} \right). \tag{0.4.6}$$

The limit relations between basic hypergeometric orthogonal polynomials described in chapter 4 of this report are based on the observations that

$$\,_r\phi_s \left( \begin{matrix} a_1, \ldots, a_{r-1}, \mu \\ b_1, \ldots, b_{s-1}, \mu \end{matrix} \bigg| q; z \right) = {}_{r-1}\phi_{s-1} \left( \begin{matrix} a_1, \ldots, a_{r-1} \\ b_1, \ldots, b_{s-1} \end{matrix} \bigg| q; z \right), \tag{0.4.7}$$

$$\lim_{\lambda \to \infty} {}_r\phi_s \left( \begin{matrix} a_1, \ldots, a_{r-1}, \lambda a_r \\ b_1, \ldots, b_s \end{matrix} \bigg| q; \frac{z}{\lambda} \right) = {}_{r-1}\phi_s \left( \begin{matrix} a_1, \ldots, a_{r-1} \\ b_1, \ldots, b_s \end{matrix} \bigg| q; a_r z \right), \tag{0.4.8}$$

$$\lim_{\lambda \to \infty} {}_r\phi_s \left( \begin{matrix} a_1, \ldots, a_r \\ b_1, \ldots, b_{s-1}, \lambda b_s \end{matrix} \bigg| q; \lambda z \right) = {}_r\phi_{s-1} \left( \begin{matrix} a_1, \ldots, a_r \\ b_1, \ldots, b_{s-1} \end{matrix} \bigg| q; \frac{z}{b_s} \right), \tag{0.4.9}$$

and

$$\lim_{\lambda \to \infty} {}_r\phi_s \left( \begin{matrix} a_1, \ldots, a_{r-1}, \lambda a_r \\ b_1, \ldots, b_{s-1}, \lambda b_s \end{matrix} \bigg| q; z \right) = {}_{r-1}\phi_{s-1} \left( \begin{matrix} a_1, \ldots, a_{r-1} \\ b_1, \ldots, b_{s-1} \end{matrix} \bigg| q; \frac{a_r z}{b_s} \right). \tag{0.4.10}$$

Mostly, the left-hand sides of the formulas (0.4.3) and (0.4.7) occur as limit cases when some numerator parameter and some denominator parameter tend to the same value.

Finally, we introduce a notation for the $N$th partial sum of a (basic) hypergeometric series. We will use this notation in the definitions of the discrete orthogonal polynomials. We also use it in order to write generating functions for discrete orthogonal polynomials in a compact way. We define

$$\,_r\tilde{F}_s \left( \begin{matrix} a_1, \ldots, a_r \\ b_1, \ldots, b_s \end{matrix} \bigg| z \right) := \sum_{k=0}^{N} \frac{(a_1, \ldots, a_r)_k}{(b_1, \ldots, b_s)_k} \frac{z^k}{k!}, \tag{0.4.11}$$

where $N$ denotes the nonnegative integer that appears in each definition of a family of discrete orthogonal polynomials. We also define

$$\,_r\tilde{\phi}_s \left( \begin{matrix} a_1, \ldots, a_r \\ b_1, \ldots, b_s \end{matrix} \bigg| q; z \right) := \sum_{k=0}^{N} \frac{(a_1, \ldots, a_r; q)_k}{(b_1, \ldots, b_s; q)_k} (-1)^{(1+s-r)k} q^{(1+s-r)\binom{k}{2}} \frac{z^k}{(q;q)_k}. \tag{0.4.12}$$

As an example of the use of these notations we remark that the definition (3.14.1) of the quantum $q$-Krawtchouk polynomials must be understood as

$$\,_2\tilde{\phi}_1 \left( \begin{matrix} q^{-n}, q^{-x} \\ q^{-N} \end{matrix} \bigg| q; pq^{n+1} \right) = \sum_{k=0}^{N} \frac{(q^{-n};q)_k (q^{-x};q)_k}{(q^{-N};q)_k (q;q)_k} \left( pq^{n+1} \right)^k.$$

In cases of discrete orthogonal polynomials, like the Racah, Hahn, dual Hahn and Krawtchouk polynomials, we need another special notation for the generating functions. In order to simplify the notation we write the generating functions as (products of) (truncated as above) power series in $t$ for which the $N$th partial sum equals the right-hand side. In these cases we use the notation $\simeq$ instead of the $=$ sign. As an example of this notation and the one mentioned above we note that the generating function (1.6.8) of the dual Hahn polynomials must be understood as follows. The $N$th partial sum of

$$e^t {}_2\tilde{F}_2\left(\begin{array}{c}-x,x+\gamma+\delta+1\\ \gamma+1,-N\end{array}\bigg|-t\right) = e^t \sum_{k=0}^{N} \frac{(-x)_k(x+\gamma+\delta+1)_k}{(\gamma+1)_k(-N)_k k!}(-t)^k$$

equals

$$\sum_{n=0}^{N} \frac{R_n(\lambda(x);\gamma,\delta,N)}{n!}t^n.$$

## 0.5 The $q$-binomial theorem and other summation formulas

One of the most important summation formulas for hypergeometric series is given by the binomial theorem :

$${}_1F_0\left(\begin{array}{c}a\\ -\end{array}\bigg|z\right) = \sum_{n=0}^{\infty}\frac{(a)_n}{n!}z^n = (1-z)^{-a}, \; |z|<1. \qquad (0.5.1)$$

A $q$-analogue of this formula is called the $q$-binomial theorem :

$${}_1\phi_0\left(\begin{array}{c}a\\ -\end{array}\bigg|q;z\right) = \sum_{n=0}^{\infty}\frac{(a;q)_n}{(q;q)_n}z^n = \frac{(az;q)_{\infty}}{(z;q)_{\infty}}, \; |z|<1. \qquad (0.5.2)$$

For $a=q^{-n}$ with $n$ a nonnegative integer we find

$${}_1\phi_0\left(\begin{array}{c}q^{-n}\\ -\end{array}\bigg|q;z\right) = (zq^{-n};q)_n, \; n=0,1,2,\ldots. \qquad (0.5.3)$$

In fact this is a $q$-analogue of Newton's binomium

$${}_1F_0\left(\begin{array}{c}-n\\ -\end{array}\bigg|z\right) = \sum_{k=0}^{n}\frac{(-n)_k}{k!}z^k = \sum_{k=0}^{n}\binom{n}{k}(-z)^k = (1-z)^n, \; n=0,1,2,\ldots. \qquad (0.5.4)$$

As an example of the use of these formulas we note that the generating function (1.6.5) of the dual Hahn polynomials can also be written as :

$${}_1F_0\left(\begin{array}{c}-x-\delta\\ -\end{array}\bigg|t\right) {}_2\tilde{F}_1\left(\begin{array}{c}x-N,x+\gamma+\delta+1\\ -N\end{array}\bigg|t\right).$$

In a similar way we find for the generating function (3.14.5) of the quantum $q$-Krawtchouk polynomials :

$$\frac{(q^{-x}t;q)_{\infty}}{(t;q)_{\infty}} {}_2\phi_1\left(\begin{array}{c}q^{x-N},0\\ pq\end{array}\bigg|q;q^{-x}t\right) = {}_1\phi_0\left(\begin{array}{c}q^{-x}\\ -\end{array}\bigg|q;t\right) {}_2\phi_1\left(\begin{array}{c}q^{x-N},0\\ pq\end{array}\bigg|q;q^{-x}t\right).$$

Another example of the use of the $q$-binomial theorem is the proof of the fact that the generating function (3.10.19) for the continuous $q$-ultraspherical (or Rogers) polynomials is a $q$-analogue of the generating function (1.8.14) for the Gegenbauer (or ultraspherical) polynomials. In fact we have, after the substitution $\beta = q^\lambda$ :

$$\left|\frac{(q^\lambda e^{i\theta}t;q)_\infty}{(e^{i\theta}t;q)_\infty}\right|^2 = \frac{(q^\lambda e^{i\theta}t, q^\lambda e^{-i\theta}t;q)_\infty}{(e^{i\theta}t, e^{-i\theta}t;q)_\infty} = {}_1\phi_0\left(\begin{array}{c}q^\lambda\\-\end{array}\bigg|\, q; e^{i\theta}t\right) {}_1\phi_0\left(\begin{array}{c}q^\lambda\\-\end{array}\bigg|\, q; e^{-i\theta}t\right),$$

which tends to (for $q \uparrow 1$)

$${}_1F_0\left(\begin{array}{c}\lambda\\-\end{array}\bigg|\, e^{i\theta}t\right) {}_1F_0\left(\begin{array}{c}\lambda\\-\end{array}\bigg|\, e^{-i\theta}t\right) = (1 - e^{i\theta}t)^{-\lambda}(1 - e^{-i\theta}t)^{-\lambda} = (1 - 2xt + t^2)^{-\lambda},\ x = \cos\theta,$$

which equals (1.8.14).

The well-known Gauss' summation formula

$${}_2F_1\left(\begin{array}{c}a,b\\c\end{array}\bigg|\, 1\right) = \frac{\Gamma(c)\Gamma(c-a-b)}{\Gamma(c-a)\Gamma(c-b)},\ \text{Re}(c-a-b) > 0 \tag{0.5.5}$$

and Vandermonde's summation formula

$${}_2F_1\left(\begin{array}{c}-n,b\\c\end{array}\bigg|\, 1\right) = \frac{(c-b)_n}{(c)_n},\ n = 0, 1, 2, \ldots \tag{0.5.6}$$

have the following $q$-analogues :

$${}_2\phi_1\left(\begin{array}{c}a,b\\c\end{array}\bigg|\, q; \frac{c}{ab}\right) = \frac{(a^{-1}c, b^{-1}c;q)_\infty}{(c, a^{-1}b^{-1}c;q)_\infty},\ \left|\frac{c}{ab}\right| < 1, \tag{0.5.7}$$

$${}_2\phi_1\left(\begin{array}{c}q^{-n},b\\c\end{array}\bigg|\, q; \frac{cq^n}{b}\right) = \frac{(b^{-1}c;q)_n}{(c;q)_n},\ n = 0, 1, 2, \ldots \tag{0.5.8}$$

and

$${}_2\phi_1\left(\begin{array}{c}q^{-n},b\\c\end{array}\bigg|\, q; q\right) = \frac{(b^{-1}c;q)_n}{(c;q)_n}b^n,\ n = 0, 1, 2, \ldots. \tag{0.5.9}$$

On the next level we have the summation formula

$${}_3F_2\left(\begin{array}{c}-n,a,b\\c, 1+a+b-c-n\end{array}\bigg|\, 1\right) = \frac{(c-a)_n(c-b)_n}{(c)_n(c-a-b)_n},\ n = 0, 1, 2, \ldots \tag{0.5.10}$$

which is called Saalschütz (or Pfaff-Saalschütz) summation formula. A $q$-analogue of this summation formula is

$${}_3\phi_2\left(\begin{array}{c}q^{-n},a,b\\c, abc^{-1}q^{1-n}\end{array}\bigg|\, q; q\right) = \frac{(a^{-1}c, b^{-1}c;q)_n}{(c, a^{-1}b^{-1}c;q)_n},\ n = 0, 1, 2, \ldots. \tag{0.5.11}$$

Finally, we have a summation formula for the ${}_1\phi_1$ series :

$${}_1\phi_1\left(\begin{array}{c}a\\c\end{array}\bigg|\, q; \frac{c}{a}\right) = \frac{(a^{-1}c;q)_\infty}{(c;q)_\infty}. \tag{0.5.12}$$

As an example of the use of this latter formula we remark that the $q$-Laguerre polynomials defined by (3.21.1) have the property that

$$L_n^{(\alpha)}(-1;q) = \frac{1}{(q;q)_n},\ n = 0, 1, 2, \ldots.$$

## 0.6 Transformation formulas

In this section we list a number of transformation formulas which can be used to transform definitions or other formulas into equivalent but different forms.

First of all we have Heine's transformation formulas for the $_2\phi_1$ series :

$$_2\phi_1\left(\begin{array}{c}a,b\\c\end{array}\bigg|\,q;z\right) = \frac{(az,b;q)_\infty}{(c,z;q)_\infty}{}_2\phi_1\left(\begin{array}{c}b^{-1}c,z\\az\end{array}\bigg|\,q;b\right) \qquad (0.6.1)$$

$$= \frac{(b^{-1}c,bz;q)_\infty}{(c,z;q)_\infty}{}_2\phi_1\left(\begin{array}{c}abc^{-1}z,b\\bz\end{array}\bigg|\,q;\frac{c}{b}\right) \qquad (0.6.2)$$

$$= \frac{(abc^{-1}z;q)_\infty}{(z;q)_\infty}{}_2\phi_1\left(\begin{array}{c}a^{-1}c,b^{-1}c\\c\end{array}\bigg|\,q;\frac{abz}{c}\right). \qquad (0.6.3)$$

The latter formula is a $q$-analogue of Euler's transformation formula :

$$_2F_1\left(\begin{array}{c}a,b\\c\end{array}\bigg|\,z\right) = (1-z)^{c-a-b}{}_2F_1\left(\begin{array}{c}c-a,c-b\\c\end{array}\bigg|\,z\right). \qquad (0.6.4)$$

Another transformation formula for the $_2F_1$ series which is also due to Euler is :

$$_2F_1\left(\begin{array}{c}a,b\\c\end{array}\bigg|\,z\right) = (1-z)^{-b}{}_2F_1\left(\begin{array}{c}c-a,b\\c\end{array}\bigg|\,\frac{z}{z-1}\right). \qquad (0.6.5)$$

This transformation formula is also known as the Pfaff-Kummer transformation formula.

As a limit case of this one we have Kummer's transformation formula for the confluent hypergeometric series :

$$_1F_1\left(\begin{array}{c}a\\c\end{array}\bigg|\,z\right) = e^z{}_1F_1\left(\begin{array}{c}c-a\\c\end{array}\bigg|\,-z\right). \qquad (0.6.6)$$

Limit cases of Heine's transformation formulas are

$$_2\phi_1\left(\begin{array}{c}0,0\\c\end{array}\bigg|\,q;z\right) = \frac{1}{(c,z;q)_\infty}{}_1\phi_1\left(\begin{array}{c}z\\0\end{array}\bigg|\,q;c\right) \qquad (0.6.7)$$

$$= \frac{1}{(z;q)_\infty}{}_0\phi_1\left(\begin{array}{c}-\\c\end{array}\bigg|\,q;cz\right), \qquad (0.6.8)$$

$$_2\phi_1\left(\begin{array}{c}a,0\\c\end{array}\bigg|\,q;z\right) = \frac{(az;q)_\infty}{(c,z;q)_\infty}{}_1\phi_1\left(\begin{array}{c}z\\az\end{array}\bigg|\,q;c\right) \qquad (0.6.9)$$

$$= \frac{1}{(z;q)_\infty}{}_1\phi_1\left(\begin{array}{c}a^{-1}c\\c\end{array}\bigg|\,q;az\right), \qquad (0.6.10)$$

$$_1\phi_1\left(\begin{array}{c}a\\c\end{array}\bigg|\,q;z\right) = \frac{(a,z;q)_\infty}{(c;q)_\infty}{}_2\phi_1\left(\begin{array}{c}a^{-1}c,0\\z\end{array}\bigg|\,q;a\right) \qquad (0.6.11)$$

$$= (ac^{-1}z;q)_\infty \cdot {}_2\phi_1\left(\begin{array}{c}a^{-1}c,0\\c\end{array}\bigg|\,q;\frac{az}{c}\right) \qquad (0.6.12)$$

and

$$_2\phi_1\left(\begin{array}{c}a,b\\0\end{array}\bigg|\,q;z\right) = \frac{(az,b;q)_\infty}{(z;q)_\infty}{}_2\phi_1\left(\begin{array}{c}z,0\\az\end{array}\bigg|\,q;b\right) \qquad (0.6.13)$$

$$= \frac{(bz;q)_\infty}{(z;q)_\infty}{}_1\phi_1\left(\begin{array}{c}b\\bz\end{array}\bigg|\,q;az\right). \qquad (0.6.14)$$

If we reverse the order of summation in a terminating $_1F_1$ series we obtain a $_2F_0$ series, in fact we have

$$_1F_1\left(\begin{array}{c}-n\\a\end{array}\bigg|\,x\right)=\frac{(-x)^n}{(a)_n}\,_2F_0\left(\begin{array}{c}-n,-a-n+1\\-\end{array}\bigg|-\frac{1}{x}\right),\ n=0,1,2,\ldots. \qquad (0.6.15)$$

If we apply this technique to a terminating $_2F_1$ series we find

$$_2F_1\left(\begin{array}{c}-n,b\\c\end{array}\bigg|\,x\right)=\frac{(b)_n}{(c)_n}(-x)^n\,_2F_1\left(\begin{array}{c}-n,-c-n+1\\-b-n+1\end{array}\bigg|\,\frac{1}{x}\right),\ n=0,1,2,\ldots. \qquad (0.6.16)$$

The $q$-analogues of these formulas are

$$_1\phi_1\left(\begin{array}{c}q^{-n}\\a\end{array}\bigg|\,q;z\right)=\frac{(q^{-1}z)^n}{(a;q)_n}\,_2\phi_1\left(\begin{array}{c}q^{-n},a^{-1}q^{1-n}\\0\end{array}\bigg|\,q;\frac{aq^{n+1}}{z}\right),\ n=0,1,2,\ldots \qquad (0.6.17)$$

and

$$_2\phi_1\left(\begin{array}{c}q^{-n},b\\c\end{array}\bigg|\,q;z\right)$$
$$=\frac{(b;q)_n}{(c;q)_n}q^{-n-\binom{n}{2}}(-z)^n\,_2\phi_1\left(\begin{array}{c}q^{-n},c^{-1}q^{1-n}\\b^{-1}q^{1-n}\end{array}\bigg|\,q;\frac{cq^{n+1}}{bz}\right),\ n=0,1,2,\ldots. \qquad (0.6.18)$$

A limit case of the latter formula is

$$_2\phi_0\left(\begin{array}{c}q^{-n},b\\-\end{array}\bigg|\,q;zq^n\right)=(b;q)_nz^n\,_2\phi_1\left(\begin{array}{c}q^{-n},0\\b^{-1}q^{1-n}\end{array}\bigg|\,q;\frac{q}{bz}\right),\ n=0,1,2,\ldots. \qquad (0.6.19)$$

The next transformation formula is due to Jackson :

$$_2\phi_1\left(\begin{array}{c}q^{-n},b\\c\end{array}\bigg|\,q;z\right)=\frac{(bc^{-1}q^{-n}z;q)_\infty}{(bc^{-1}z;q)_\infty}\,_3\phi_2\left(\begin{array}{c}q^{-n},b^{-1}c,0\\c,b^{-1}cqz^{-1}\end{array}\bigg|\,q;q\right),\ n=0,1,2,\ldots. \qquad (0.6.20)$$

Equivalently we have

$$_3\phi_2\left(\begin{array}{c}q^{-n},a,0\\b,c\end{array}\bigg|\,q;q\right)=\frac{(b^{-1}q;q)_\infty}{(b^{-1}q^{1-n};q)_\infty}\,_2\phi_1\left(\begin{array}{c}q^{-n},a^{-1}c\\c\end{array}\bigg|\,q;\frac{aq}{b}\right),\ n=0,1,2,\ldots. \qquad (0.6.21)$$

Other transformation formulas of this kind are given by :

$$_2\phi_1\left(\begin{array}{c}q^{-n},b\\c\end{array}\bigg|\,q;z\right)$$
$$=\frac{(b^{-1}c;q)_n}{(c;q)_n}\left(\frac{bz}{q}\right)^n\,_3\phi_2\left(\begin{array}{c}q^{-n},qz^{-1},c^{-1}q^{1-n}\\bc^{-1}q^{1-n},0\end{array}\bigg|\,q;q\right) \qquad (0.6.22)$$
$$=\frac{(b^{-1}c;q)_n}{(c;q)_n}\,_3\phi_2\left(\begin{array}{c}q^{-n},b,bc^{-1}q^{-n}z\\bc^{-1}q^{1-n},0\end{array}\bigg|\,q;q\right),\ n=0,1,2,\ldots, \qquad (0.6.23)$$

or equivalently

$$_3\phi_2\left(\begin{array}{c}q^{-n},a,b\\c,0\end{array}\bigg|\,q;q\right)=\frac{(b;q)_n}{(c;q)_n}a^n\,_2\phi_1\left(\begin{array}{c}q^{-n},b^{-1}c\\b^{-1}q^{1-n}\end{array}\bigg|\,q;\frac{q}{a}\right) \qquad (0.6.24)$$
$$=\frac{(a^{-1}c;q)_n}{(c;q)_n}a^n\,_2\phi_1\left(\begin{array}{c}q^{-n},a\\ac^{-1}q^{1-n}\end{array}\bigg|\,q;\frac{bq}{c}\right),\ n=0,1,2,\ldots. \qquad (0.6.25)$$

Limit cases of these formulas are

$$_2\phi_0\left(\begin{array}{c}q^{-n},b\\-\end{array}\bigg|\,q;z\right)=b^{-n}\,_3\phi_2\left(\begin{array}{c}q^{-n},b,bzq^{-n}\\0,0\end{array}\bigg|\,q;q\right),\ n=0,1,2,\ldots, \qquad (0.6.26)$$

or equivalently

$$_3\phi_2\left(\begin{array}{c}q^{-n},a,b\\0,0\end{array}\bigg|q;q\right) = (b;q)_n a^n {}_2\phi_1\left(\begin{array}{c}q^{-n},0\\b^{-1}q^{1-n}\end{array}\bigg|q;\frac{q}{a}\right) \qquad (0.6.27)$$

$$= a^n {}_2\phi_0\left(\begin{array}{c}q^{-n},a\\-\end{array}\bigg|q;\frac{bq^n}{a}\right), \; n=0,1,2,\ldots. \qquad (0.6.28)$$

On the next level we have Sears' transformation formula for a terminating balanced $_4\phi_3$ series :

$$_4\phi_3\left(\begin{array}{c}q^{-n},a,b,c\\d,e,f\end{array}\bigg|q;q\right)$$
$$= \frac{(a^{-1}e,a^{-1}f;q)_n}{(e,f;q)_n}a^n {}_4\phi_3\left(\begin{array}{c}q^{-n},a,b^{-1}d,c^{-1}d\\d,ae^{-1}q^{1-n},af^{-1}q^{1-n}\end{array}\bigg|q;q\right) \qquad (0.6.29)$$
$$= \frac{(a,a^{-1}b^{-1}ef,a^{-1}c^{-1}ef;q)_n}{(e,f,a^{-1}b^{-1}c^{-1}ef;q)_n} \times$$
$$_4\phi_3\left(\begin{array}{c}q^{-n},a^{-1}e,a^{-1}f,a^{-1}b^{-1}c^{-1}ef\\a^{-1}b^{-1}ef,a^{-1}c^{-1}ef,a^{-1}q^{1-n}\end{array}\bigg|q;q\right), \; def=abcq^{1-n}. \qquad (0.6.30)$$

Sears' transformation formula is a $q$-analogue of Whipple's transformation formula for a terminating balanced $_4F_3$ series :

$$_4F_3\left(\begin{array}{c}-n,a,b,c\\d,e,f\end{array}\bigg|1\right) = \frac{(e-a)_n(f-a)_n}{(e)_n(f)_n} \times$$
$$_4F_3\left(\begin{array}{c}-n,a,d-b,d-c\\d,a-e-n+1,a-f-n+1\end{array}\bigg|1\right), \; a+b+c+1=d+e+f+n. \qquad (0.6.31)$$

Whipple's formula can be used to show that the Wilson polynomials defined by (1.1.1) are symmetric in their parameters in the sense that the following 24 different forms are all equal :

$$W_n(x^2;a,b,c,d) = W_n(x^2;a,b,d,c) = W_n(x^2;a,c,b,d) = \cdots = W_n(x^2;d,c,b,a).$$

Sears' transformation formula can be used to derive similar symmetry relations for the Askey-Wilson polynomials defined by (3.1.1) :

$$p_n(x;a,b,c,d) = p_n(x;a,b,d,c) = p_n(x;a,c,b,d) = \cdots = p_n(x;d,c,b,a).$$

Finally, we mention a quadratic transformation formula which is due to Singh :

$$_4\phi_3\left(\begin{array}{c}a^2,b^2,c,d\\abq^{\frac{1}{2}},-abq^{\frac{1}{2}},-cd\end{array}\bigg|q;q\right) = {}_4\phi_3\left(\begin{array}{c}a^2,b^2,c^2,d^2\\a^2b^2q,-cd,-cdq\end{array}\bigg|q^2;q^2\right), \qquad (0.6.32)$$

which is valid when both sides terminate.

If we apply Singh's formula (0.6.32) to the continuous $q$-Jacobi polynomials defined by (3.10.1) and (3.10.2) and use Sears' transformation formula (0.6.29), formula (0.2.10) twice and also formula (0.2.18), then we find the quadratic transformation

$$P_n^{(\alpha,\beta)}(x|q^2) = \frac{(-q;q)_n}{(-q^{\alpha+\beta+1};q)_n}q^{n\alpha}P_n^{(\alpha,\beta)}(x;q).$$

## 0.7 Some special functions and their $q$-analogues

The classical exponential function $\exp(z)$ and the trigonometric functions $\sin(z)$ and $\cos(z)$ can be expressed in terms of hypergeometric functions as

$$\exp(z) = e^z = {}_0F_0\left(\begin{array}{c}-\\-\end{array}\bigg|z\right), \qquad (0.7.1)$$

$$\sin(z) = z \, _0F_1 \left( \begin{array}{c} - \\ \frac{3}{2} \end{array} \middle| -\frac{z^2}{4} \right) \tag{0.7.2}$$

and

$$\cos(z) = {}_0F_1 \left( \begin{array}{c} - \\ \frac{1}{2} \end{array} \middle| -\frac{z^2}{4} \right). \tag{0.7.3}$$

Further we have the well-known Bessel function $J_\nu(z)$ which can be defined by

$$J_\nu(z) := \frac{\left(\frac{z}{2}\right)^\nu}{\Gamma(\nu+1)} {}_0F_1 \left( \begin{array}{c} - \\ \nu+1 \end{array} \middle| -\frac{z^2}{4} \right). \tag{0.7.4}$$

Applying this formula to the generating function (1.11.6) of the Laguerre polynomials we obtain :

$$(xt)^{-\frac{\alpha}{2}} e^t J_\alpha(2\sqrt{xt}) = \frac{1}{\Gamma(\alpha+1)} \sum_{n=0}^{\infty} \frac{L_n^{(\alpha)}(x)}{(\alpha+1)_n} t^n.$$

These functions all have several $q$-analogues. The exponential function for instance has two different natural $q$-extensions, denoted by $e_q(z)$ and $E_q(z)$ defined by

$$e_q(z) := {}_1\phi_0 \left( \begin{array}{c} 0 \\ - \end{array} \middle| q; z \right) = \sum_{n=0}^{\infty} \frac{z^n}{(q;q)_n} \tag{0.7.5}$$

and

$$E_q(z) := {}_0\phi_0 \left( \begin{array}{c} - \\ - \end{array} \middle| q; -z \right) = \sum_{n=0}^{\infty} \frac{q^{\binom{n}{2}}}{(q;q)_n} z^n. \tag{0.7.6}$$

These $q$-analogues of the exponential function are related by

$$e_q(z) E_q(-z) = 1.$$

They are $q$-extensions of the exponential function since

$$\lim_{q \uparrow 1} e_q((1-q)z) = \lim_{q \uparrow 1} E_q((1-q)z) = e^z.$$

If we set $a = 0$ in the $q$-binomial theorem we find for the $q$-exponential functions :

$$e_q(z) = {}_1\phi_0 \left( \begin{array}{c} 0 \\ - \end{array} \middle| q; z \right) = \frac{1}{(z;q)_\infty}, \ |z| < 1. \tag{0.7.7}$$

Further we have

$$E_q(z) = {}_0\phi_0 \left( \begin{array}{c} - \\ - \end{array} \middle| q; -z \right) = (-z;q)_\infty. \tag{0.7.8}$$

For instance, these formulas can be used to obtain other versions of a generating function for several sets of orthogonal polynomials mentioned in this report.

If we assume that $|z| < 1$ we may define

$$\sin_q(z) := \frac{e_q(iz) - e_q(-iz)}{2i} = \sum_{n=0}^{\infty} \frac{(-1)^n z^{2n+1}}{(q;q)_{2n+1}} \tag{0.7.9}$$

and

$$\cos_q(z) := \frac{e_q(iz) + e_q(-iz)}{2} = \sum_{n=0}^{\infty} \frac{(-1)^n z^{2n}}{(q;q)_{2n}}. \tag{0.7.10}$$

These are $q$-analogues of the trigonometric functions $\sin(z)$ and $\cos(z)$. On the other hand we may define

$$\mathrm{Sin}_q(z) := \frac{E_q(iz) - E_q(-iz)}{2i} \tag{0.7.11}$$

and
$$\text{Cos}_q(z) := \frac{E_q(iz) + E_q(-iz)}{2}. \tag{0.7.12}$$

Then it is not very difficult to verify that

$$e_q(iz) = \cos_q(z) + i\sin_q(z) \text{ and } E_q(iz) = \text{Cos}_q(z) + i\,\text{Sin}_q(z).$$

Further we have

$$\begin{cases} \sin_q(z)\text{Sin}_q(z) + \cos_q(z)\text{Cos}_q(z) = 1 \\ \sin_q(z)\text{Cos}_q(z) - \text{Sin}_q(z)\cos_q(z) = 0. \end{cases}$$

The $q$-analogues of the trigonometric functions can be used to find different forms of formulas appearing in this report, although we will not use them.

Some $q$-analogues of the Bessel functions are given by

$$J_\nu^{(1)}(z;q) := \frac{(q^{\nu+1};q)_\infty}{(q;q)_\infty} \left(\frac{z}{2}\right)^\nu {}_2\phi_1\left(\begin{array}{c} 0,0 \\ q^{\nu+1} \end{array} \bigg| q; -\frac{z^2}{4}\right) \tag{0.7.13}$$

and

$$J_\nu^{(2)}(z;q) := \frac{(q^{\nu+1};q)_\infty}{(q;q)_\infty} \left(\frac{z}{2}\right)^\nu {}_0\phi_1\left(\begin{array}{c} - \\ q^{\nu+1} \end{array} \bigg| q; -\frac{q^{\nu+1}z^2}{4}\right). \tag{0.7.14}$$

These $q$-Bessel functions are connected by

$$J_\nu^{(2)}(z;q) = (-\frac{z^2}{4};q)_\infty \cdot J_\nu^{(1)}(z;q),\ |z| < 2.$$

They are $q$-analogues of the Bessel function since

$$\lim_{q\uparrow 1} J_\nu^{(k)}((1-q)z;q) = J_\nu(z),\ k=1,2.$$

These $q$-Bessel functions were introduced by F.H. Jackson in 1905. They are therefore referred to as Jackson $q$-Bessel functions. Other $q$-analogues of the Bessel function are the so-called Hahn-Exton $q$-Bessel functions.

As an example we remark that the generating function (3.20.5) for the little $q$-Laguerre (or Wall) polynomials can also be written as

$$\frac{(-t;q)_\infty(q;q)_\infty}{(q^{\alpha+1};q)_\infty}(xt)^{-\frac{\alpha}{2}} J_\alpha^{(1)}(2\sqrt{xt};q) = \sum_{n=0}^\infty \frac{q^{\binom{n}{2}}}{(q;q)_n} p_n(x;q^\alpha|q)t^n$$

or as

$$\frac{(q;q)_\infty}{(q^{\alpha+1};q)_\infty}(xt)^{-\frac{\alpha}{2}} E_q(t) J_\alpha^{(1)}(2\sqrt{xt};q) = \sum_{n=0}^\infty \frac{q^{\binom{n}{2}}}{(q;q)_n} p_n(x;q^\alpha|q)t^n.$$

## 0.8 The $q$-derivative and the $q$-integral

The $q$-derivative operator $\mathcal{D}_q$ is defined by

$$\mathcal{D}_q f(z) := \begin{cases} \dfrac{f(z) - f(qz)}{(1-q)z}, & z \neq 0 \\ f'(0), & z = 0. \end{cases} \tag{0.8.1}$$

Further we define
$$\mathcal{D}_q^0 f := f \text{ and } \mathcal{D}_q^n f := \mathcal{D}_q\left(\mathcal{D}_q^{n-1}f\right),\ n=1,2,3,\ldots \tag{0.8.2}$$

It is not very difficult to see that
$$\lim_{q \uparrow 1} \mathcal{D}_q f(z) = f'(z)$$
if the function $f$ is differentiable at $z$.

An easy consequence of this definition is
$$\mathcal{D}_q \left[ f(\gamma x) \right] = \gamma \left( \mathcal{D}_q f \right)(\gamma x) \tag{0.8.3}$$
for all real $\gamma$ or more general
$$\mathcal{D}_q^n \left[ f(\gamma x) \right] = \gamma^n \left( \mathcal{D}_q^n f \right)(\gamma x), \ n = 0, 1, 2, \ldots. \tag{0.8.4}$$

Further we have
$$\mathcal{D}_q \left[ f(x) g(x) \right] = f(qx) \mathcal{D}_q g(x) + g(x) \mathcal{D}_q f(x) \tag{0.8.5}$$
which is often referred to as the $q$-product rule. This can be generalized to a $q$-analogue of Leibniz' rule :
$$\mathcal{D}_q^n \left[ f(x) g(x) \right] = \sum_{k=0}^{n} \begin{bmatrix} n \\ k \end{bmatrix}_q \left( \mathcal{D}_q^{n-k} f \right) (q^k x) \left( \mathcal{D}_q^k g \right)(x), \ n = 0, 1, 2, \ldots. \tag{0.8.6}$$

As an example we note that the $q$-difference equation (3.21.5) of the $q$-Laguerre polynomials can also be written in terms of this $q$-derivative operator as
$$(1-q)^2 x \mathcal{D}_q^2 y(x) + (1-q) \left[ 1 - q^{\alpha+1} - q^{\alpha+2} x \right] (\mathcal{D}_q y)(qx) + (1-q^n) q^{\alpha+1} y(qx) = 0, \ y(x) = L_n^{(\alpha)}(x; q).$$

The $q$-integral is defined by
$$\int_0^z f(t) d_q t := z(1-q) \sum_{n=0}^{\infty} f(zq^n) q^n. \tag{0.8.7}$$

This definition is due to J. Thomae and F.H. Jackson. Jackson also defined a $q$-integral on $(0, \infty)$ by
$$\int_0^{\infty} f(t) d_q t := (1-q) \sum_{n=-\infty}^{\infty} f(q^n) q^n. \tag{0.8.8}$$

If the function $f$ is continuous on $[0, z]$ we have
$$\lim_{q \uparrow 1} \int_0^z f(t) d_q t = \int_0^z f(t) dt.$$

For instance, the orthogonality relation (3.12.2) for the little $q$-Jacobi polynomials can also be written in terms of a $q$-integral as :
$$\int_0^1 \frac{(qx; q)_{\infty}}{(q^{\beta+1} x; q)_{\infty}} x^{\alpha} p_m(x; q^{\alpha}, q^{\beta} | q) p_n(x; q^{\alpha}, q^{\beta} | q) d_q x$$
$$= (1-q) \frac{(q, q^{\alpha+\beta+2}; q)_{\infty}}{(q^{\alpha+1}, q^{\beta+1}; q)_{\infty}} \frac{(1-q^{\alpha+\beta+1})}{(1-q^{2n+\alpha+\beta+1})} \frac{(q, q^{\beta+1}; q)_n}{(q^{\alpha+1}, q^{\alpha+\beta+1}; q)_n} q^{n(\alpha+1)} \delta_{mn}, \ \alpha > -1, \ \beta > -1.$$

# ASKEY-SCHEME
## OF
# HYPERGEOMETRIC
# ORTHOGONAL POLYNOMIALS

$_4F_3(4)$ — Wilson, Racah

$_3F_2(3)$ — Continuous dual Hahn, Continuous Hahn, Hahn, Dual Hahn

$_2F_1(2)$ — Meixner-Pollaczek, Jacobi, Meixner, Krawtchouk

$_1F_1(1)/_2F_0(1)$ — Laguerre, Charlier

$_2F_0(0)$ — Hermite

# Chapter 1

# Hypergeometric orthogonal polynomials

## 1.1 Wilson

**Definition.**
$$\frac{W_n(x^2;a,b,c,d)}{(a+b)_n(a+c)_n(a+d)_n} = {}_4F_3\left(\begin{array}{c}-n, n+a+b+c+d-1, a+ix, a-ix\\a+b, a+c, a+d\end{array}\bigg| 1\right). \quad (1.1.1)$$

**Orthogonality.** When $Re(a,b,c,d) > 0$ and non-real parameters occur in conjugate pairs, then

$$\frac{1}{2\pi}\int_0^\infty \left|\frac{\Gamma(a+ix)\Gamma(b+ix)\Gamma(c+ix)\Gamma(d+ix)}{\Gamma(2ix)}\right|^2 W_m(x^2;a,b,c,d)W_n(x^2;a,b,c,d)dx$$
$$= (n+a+b+c+d-1)_n n! \frac{\Gamma(n+a+b)\cdots\Gamma(n+c+d)}{\Gamma(2n+a+b+c+d)}\delta_{mn}, \quad (1.1.2)$$

where

$$\Gamma(n+a+b)\cdots\Gamma(n+c+d)$$
$$= \Gamma(n+a+b)\Gamma(n+a+c)\Gamma(n+a+d)\Gamma(n+b+c)\Gamma(n+b+d)\Gamma(n+c+d).$$

If $a < 0$ and $a+b$, $a+c$, $a+d$ are positive or a pair of complex conjugates occur with positive real parts, then

$$\frac{1}{2\pi}\int_0^\infty \left|\frac{\Gamma(a+ix)\Gamma(b+ix)\Gamma(c+ix)\Gamma(d+ix)}{\Gamma(2ix)}\right|^2 W_m(x^2;a,b,c,d)W_n(x^2;a,b,c,d)dx +$$
$$+ \frac{\Gamma(a+b)\Gamma(a+c)\Gamma(a+d)\Gamma(b-a)\Gamma(c-a)\Gamma(d-a)}{\Gamma(-2a)} \times$$
$$\times \sum_{\substack{k=0,1,2\ldots\\a+k<0}} \frac{(2a)_k(a+1)_k(a+b)_k(a+c)_k(a+d)_k}{(k!)(a)_k(a-b+1)_k(a-c+1)_k(a-d+1)_k}W_m(-(a+k)^2)W_n(-(a+k)^2)$$
$$= (n+a+b+c+d-1)_n n! \frac{\Gamma(n+a+b)\cdots\Gamma(n+c+d)}{\Gamma(2n+a+b+c+d)}\delta_{mn}, \quad (1.1.3)$$

where
$$W_m(-(a+k)^2)W_n(-(a+k)^2) = W_m(-(a+k)^2; a, b, c, d)W_n(-(a+k)^2; a, b, c, d).$$

**Recurrence relation.**
$$-\left(a^2 + x^2\right)\tilde{W}_n(x^2) = A_n\tilde{W}_{n+1}(x^2) - (A_n + C_n)\tilde{W}_n(x^2) + C_n\tilde{W}_{n-1}(x^2), \quad (1.1.4)$$

where
$$\tilde{W}_n(x^2) := \tilde{W}_n(x^2; a, b, c, d) = \frac{W_n(x^2; a, b, c, d)}{(a+b)_n(a+c)_n(a+d)_n}$$

and
$$\begin{cases} A_n = \dfrac{(n+a+b+c+d-1)(n+a+b)(n+a+c)(n+a+d)}{(2n+a+b+c+d-1)(2n+a+b+c+d)} \\ \\ C_n = \dfrac{n(n+b+c-1)(n+b+d-1)(n+c+d-1)}{(2n+a+b+c+d-2)(2n+a+b+c+d-1)}. \end{cases}$$

**Difference equation.**
$$n(n+a+b+c+d-1)y(x) = B(x)y(x+i) - [B(x)+D(x)]y(x) + D(x)y(x-i), \quad (1.1.5)$$

where
$$y(x) = W_n(x^2; a, b, c, d)$$

and
$$\begin{cases} B(x) = \dfrac{(a-ix)(b-ix)(c-ix)(d-ix)}{2ix(2ix-1)} \\ \\ D(x) = \dfrac{(a+ix)(b+ix)(c+ix)(d+ix)}{2ix(2ix+1)}. \end{cases}$$

**Generating functions.**
$$_2F_1\left(\begin{array}{c} a+ix, b+ix \\ a+b \end{array}\bigg| t\right) {}_2F_1\left(\begin{array}{c} c-ix, d-ix \\ c+d \end{array}\bigg| t\right) = \sum_{n=0}^{\infty}\frac{W_n(x^2; a, b, c, d)t^n}{(a+b)_n(c+d)_n n!}. \quad (1.1.6)$$

$$_2F_1\left(\begin{array}{c} a+ix, c+ix \\ a+c \end{array}\bigg| t\right) {}_2F_1\left(\begin{array}{c} b-ix, d-ix \\ b+d \end{array}\bigg| t\right) = \sum_{n=0}^{\infty}\frac{W_n(x^2; a, b, c, d)t^n}{(a+c)_n(b+d)_n n!}. \quad (1.1.7)$$

$$_2F_1\left(\begin{array}{c} a+ix, d+ix \\ a+d \end{array}\bigg| t\right) {}_2F_1\left(\begin{array}{c} b-ix, c-ix \\ b+c \end{array}\bigg| t\right) = \sum_{n=0}^{\infty}\frac{W_n(x^2; a, b, c, d)t^n}{(a+d)_n(b+c)_n n!}. \quad (1.1.8)$$

$$(1-t)^{1-a-b-c-d} {}_4F_3\left(\begin{array}{c} \frac{1}{2}(a+b+c+d-1), \frac{1}{2}(a+b+c+d), a+ix, a-ix \\ a+b, a+c, a+d \end{array}\bigg| -\frac{4t}{(1-t)^2}\right)$$
$$= \sum_{n=0}^{\infty}\frac{(a+b+c+d-1)_n}{(a+b)_n(a+c)_n(a+d)_n n!}W_n(x^2; a, b, c, d)t^n. \quad (1.1.9)$$

**Remark.** If we set
$$a = \frac{1}{2}(\gamma + \delta + 1) \ ; \ b = \frac{1}{2}(2\alpha - \gamma - \delta + 1)$$
$$c = \frac{1}{2}(2\beta - \gamma + \delta + 1) \ ; \ d = \frac{1}{2}(\gamma - \delta + 1)$$

and
$$ix \to x + \frac{1}{2}(\gamma + \delta + 1)$$

in
$$\tilde{W}_n(x^2;a,b,c,d) = \frac{W_n(x^2;a,b,c,d)}{(a+b)_n(a+c)_n(a+d)_n},$$

defined by (1.1.1) and take

$\alpha + 1 = -N$ or $\beta + \delta + 1 = -N$ or $\gamma + 1 = -N$, with $N$ a nonnegative integer

we obtain the Racah polynomials defined by (1.2.1).

**References.** [31], [44], [45], [132], [133], [155], [156], [171], [175], [227], [228].

## 1.2 Racah

**Definition.**
$$R_n(\lambda(x);\alpha,\beta,\gamma,\delta) = {}_4\tilde{F}_3\left(\begin{array}{c}-n, n+\alpha+\beta+1, -x, x+\gamma+\delta+1 \\ \alpha+1, \beta+\delta+1, \gamma+1\end{array}\bigg| 1\right), \; n=0,1,2,\ldots,N, \quad (1.2.1)$$

where
$$\lambda(x) = x(x+\gamma+\delta+1)$$

and

$\alpha + 1 = -N$ or $\beta + \delta + 1 = -N$ or $\gamma + 1 = -N$, with $N$ a nonnegative integer.

**Orthogonality.**

$$\sum_{x=0}^{N} \frac{(\gamma+\delta+1)_x((\gamma+\delta+3)/2)_x(\alpha+1)_x(\beta+\delta+1)_x(\gamma+1)_x}{(x!)((\gamma+\delta+1)/2)_x(\gamma+\delta-\alpha+1)_x(\gamma-\beta+1)_x(\delta+1)_x} R_m(\lambda(x))R_n(\lambda(x))$$
$$= M \frac{(n+\alpha+\beta+1)_n(\beta+1)_n(\alpha-\delta+1)_n(\alpha+\beta-\gamma+1)_n n!}{(\alpha+\beta+2)_{2n}(\alpha+1)_n(\beta+\delta+1)_n(\gamma+1)_n} \delta_{mn}, \quad (1.2.2)$$

where
$$R_n(\lambda(x)) := R_n(\lambda(x);\alpha,\beta,\gamma,\delta)$$

and

$$M = \begin{cases} \dfrac{(\gamma+\delta+2)_N(-\beta)_N}{(\gamma-\beta+1)_N(\delta+1)_N} & \text{if} \quad \alpha+1 = -N \\[2ex] \dfrac{(\gamma+\delta+2)_N(\delta-\alpha)_N}{(\gamma+\delta-\alpha+1)_N(\delta+1)_N} & \text{if} \quad \beta+\delta+1 = -N \\[2ex] \dfrac{(-\delta)_N(\alpha+\beta+2)_N}{(\alpha-\delta+1)_N(\beta+1)_N} & \text{if} \quad \gamma+1 = -N. \end{cases}$$

**Recurrence relation.**
$$\lambda(x)R_n(\lambda(x)) = A_n R_{n+1}(\lambda(x)) - (A_n + C_n) R_n(\lambda(x)) + C_n R_{n-1}(\lambda(x)), \quad (1.2.3)$$

where
$$R_n(\lambda(x)) := R_n(\lambda(x);\alpha,\beta,\gamma,\delta)$$

and
$$\begin{cases} A_n = \dfrac{(n+\alpha+\beta+1)(n+\alpha+1)(n+\beta+\delta+1)(n+\gamma+1)}{(2n+\alpha+\beta+1)(2n+\alpha+\beta+2)} \\[2ex] C_n = \dfrac{n(n+\beta)(n+\alpha+\beta-\gamma)(n+\alpha-\delta)}{(2n+\alpha+\beta)(2n+\alpha+\beta+1)}, \end{cases}$$

hence

$$A_n = \begin{cases} \dfrac{(n+\beta-N)(n-N)(n+\beta+\delta+1)(n+\gamma+1)}{(2n+\beta-N)(2n+\beta-N+1)} & \text{if} \quad \alpha+1 = -N \\[2ex] \dfrac{(n+\alpha+\beta+1)(n+\alpha+1)(n-N)(n+\gamma+1)}{(2n+\alpha+\beta+1)(2n+\alpha+\beta+2)} & \text{if} \quad \beta+\delta+1 = -N \\[2ex] \dfrac{(n+\alpha+\beta+1)(n+\alpha+1)(n+\beta+\delta+1)(n-N)}{(2n+\alpha+\beta+1)(2n+\alpha+\beta+2)} & \text{if} \quad \gamma+1 = -N \end{cases}$$

and

$$C_n = \begin{cases} \dfrac{n(n+\beta)(n+\beta-\gamma-N-1)(n-\delta-N-1)}{(2n+\beta-N-1)(2n+\beta-N)} & \text{if} \quad \alpha+1 = -N \\[2ex] \dfrac{n(n+\beta)(n+\alpha+\beta-\gamma)(n+\alpha+\beta+N+1)}{(2n+\alpha+\beta)(2n+\alpha+\beta+1)} & \text{if} \quad \beta+\delta+1 = -N \\[2ex] \dfrac{n(n+\beta)(n+\alpha+\beta+N+1)(n+\alpha-\delta)}{(2n+\alpha+\beta)(2n+\alpha+\beta+1)} & \text{if} \quad \gamma+1 = -N. \end{cases}$$

**Difference equation.**

$$n(n+\alpha+\beta+1)y(x) = B(x)y(x+1) - [B(x)+D(x)]y(x) + D(x)y(x-1), \quad (1.2.4)$$

where

$$y(x) = R_n(\lambda(x); \alpha, \beta, \gamma, \delta)$$

and

$$\begin{cases} B(x) = \dfrac{(x+\alpha+1)(x+\beta+\delta+1)(x+\gamma+1)(x+\gamma+\delta+1)}{(2x+\gamma+\delta+1)(2x+\gamma+\delta+2)} \\[2ex] D(x) = \dfrac{x(x+\delta)(x-\beta+\gamma)(x-\alpha+\gamma+\delta)}{(2x+\gamma+\delta)(2x+\gamma+\delta+1)}. \end{cases}$$

**Generating functions.**

$$_2\tilde{F}_1\left(\begin{array}{c} x+\alpha+1, x+\gamma+\delta+1 \\ \alpha+1 \end{array}\bigg| t\right) {}_2F_1\left(\begin{array}{c} -x+\beta-\gamma, -x-\delta \\ \beta+1 \end{array}\bigg| t\right)$$
$$\simeq \sum_{n=0}^{N} \frac{(\beta+\delta+1)_n(\gamma+1)_n}{(\beta+1)_n n!} R_n(\lambda(x); \alpha, \beta, \gamma, \delta) t^n. \qquad (1.2.5)$$

$$_2\tilde{F}_1\left(\begin{array}{c} x+\beta+\delta+1, x+\gamma+\delta+1 \\ \beta+\delta+1 \end{array}\bigg| t\right) {}_2F_1\left(\begin{array}{c} -x+\alpha-\gamma-\delta, -x-\delta \\ \alpha-\delta+1 \end{array}\bigg| t\right)$$
$$\simeq \sum_{n=0}^{N} \frac{(\alpha+1)_n(\gamma+1)_n}{(\alpha-\delta+1)_n n!} R_n(\lambda(x); \alpha, \beta, \gamma, \delta) t^n. \qquad (1.2.6)$$

$$_2\tilde{F}_1\left(\begin{array}{c} x+\gamma+1, x+\gamma+\delta+1 \\ \gamma+1 \end{array}\bigg| t\right) {}_2F_1\left(\begin{array}{c} -x+\alpha-\gamma-\delta, -x+\beta-\gamma \\ \alpha+\beta-\gamma+1 \end{array}\bigg| t\right)$$
$$\simeq \sum_{n=0}^{N} \frac{(\alpha+1)_n(\beta+\delta+1)_n}{(\alpha+\beta-\gamma+1)_n n!} R_n(\lambda(x); \alpha, \beta, \gamma, \delta) t^n. \qquad (1.2.7)$$

$$(1-t)^{-\alpha-\beta-1} {}_4\tilde{F}_3 \left( \begin{array}{c} \frac{1}{2}(\alpha+\beta+1), \frac{1}{2}(\alpha+\beta+2), -x, x+\gamma+\delta+1 \\ \alpha+1, \beta+\delta+1, \gamma+1 \end{array} \bigg| -\frac{4t}{(1-t)^2} \right)$$
$$\simeq \sum_{n=0}^{N} \frac{(\alpha+\beta+1)_n}{n!} R_n(\lambda(x); \alpha, \beta, \gamma, \delta) t^n. \qquad (1.2.8)$$

**Remark.** If we set $\alpha = a+b-1$, $\beta = c+d-1$, $\gamma = a+d-1$, $\delta = a-d$ and $x \to -a+ix$ in the definition (1.2.1) of the Racah polynomials we obtain the Wilson polynomials defined by (1.1.1):

$$R_n(\lambda(-a+ix); a+b-1, c+d-1, a+d-1, a-d)$$
$$= \tilde{W}_n(x^2; a, b, c, d) = \frac{W_n(x^2; a, b, c, d)}{(a+b)_n(a+c)_n(a+d)_n}.$$

**References.** [31], [43], [45], [47], [90], [156], [180], [183], [192], [194], [227].

## 1.3 Continuous dual Hahn

**Definition.**
$$\frac{S_n(x^2; a, b, c)}{(a+b)_n(a+c)_n} = {}_3F_2 \left( \begin{array}{c} -n, a+ix, a-ix \\ a+b, a+c \end{array} \bigg| 1 \right). \qquad (1.3.1)$$

**Orthogonality.** When $a,b$ and $c$ are positive except possibly for a pair of complex conjugates with positive real parts, then

$$\frac{1}{2\pi} \int_0^\infty \left| \frac{\Gamma(a+ix)\Gamma(b+ix)\Gamma(c+ix)}{\Gamma(2ix)} \right|^2 S_m(x^2; a,b,c) S_n(x^2; a,b,c) dx$$
$$= \Gamma(n+a+b)\Gamma(n+a+c)\Gamma(n+b+c) n! \delta_{mn}. \qquad (1.3.2)$$

If $a < 0$ and $a+b$, $a+c$ are positive or a pair of complex conjugates with positive real parts, then

$$\frac{1}{2\pi} \int_0^\infty \left| \frac{\Gamma(a+ix)\Gamma(b+ix)\Gamma(c+ix)}{\Gamma(2ix)} \right|^2 S_m(x^2; a,b,c) S_n(x^2; a,b,c) dx +$$
$$+ \frac{\Gamma(a+b)\Gamma(a+c)\Gamma(b-a)\Gamma(c-a)}{\Gamma(-2a)} \times$$
$$\times \sum_{\substack{k=0,1,2\ldots \\ a+k<0}} \frac{(2a)_k(a+1)_k(a+b)_k(a+c)_k}{(k!)(a)_k(a-b+1)_k(a-c+1)_k} (-1)^k S_m(-(a+k)^2) S_n(-(a+k)^2)$$
$$= \Gamma(n+a+b)\Gamma(n+a+c)\Gamma(n+b+c) n! \delta_{mn}, \qquad (1.3.3)$$

where
$$S_m(-(a+k)^2) S_n(-(a+k)^2) = S_m(-(a+k)^2; a, b, c) S_n(-(a+k)^2; a, b, c).$$

**Recurrence relation.**
$$-(a^2 + x^2) \tilde{S}_n(x^2) = A_n \tilde{S}_{n+1}(x^2) - (A_n + C_n) \tilde{S}_n(x^2) + C_n \tilde{S}_{n-1}(x^2), \qquad (1.3.4)$$

where
$$\tilde{S}_n(x^2) := \tilde{S}_n(x^2; a, b, c) = \frac{S_n(x^2; a, b, c)}{(a+b)_n(a+c)_n}$$

and
$$\begin{cases} A_n = (n+a+b)(n+a+c) \\ C_n = n(n+b+c-1). \end{cases}$$

**Difference equation.**

$$ny(x) = B(x)y(x+i) - [B(x) + D(x)]y(x) + D(x)y(x-i), \ y(x) = S_n(x^2;a,b,c), \quad (1.3.5)$$

where
$$\begin{cases} B(x) = \dfrac{(a-ix)(b-ix)(c-ix)}{2ix(2ix-1)} \\ D(x) = \dfrac{(a+ix)(b+ix)(c+ix)}{2ix(2ix+1)}. \end{cases}$$

**Generating functions.**

$$(1-t)^{-c+ix} {}_2F_1\left(\begin{array}{c} a+ix, b+ix \\ a+b \end{array} \bigg| t\right) = \sum_{n=0}^{\infty} \frac{S_n(x^2;a,b,c)}{(a+b)_n n!} t^n. \quad (1.3.6)$$

$$(1-t)^{-b+ix} {}_2F_1\left(\begin{array}{c} a+ix, c+ix \\ a+c \end{array} \bigg| t\right) = \sum_{n=0}^{\infty} \frac{S_n(x^2;a,b,c)}{(a+c)_n n!} t^n. \quad (1.3.7)$$

$$(1-t)^{-a+ix} {}_2F_1\left(\begin{array}{c} b+ix, c+ix \\ b+c \end{array} \bigg| t\right) = \sum_{n=0}^{\infty} \frac{S_n(x^2;a,b,c)}{(b+c)_n n!} t^n. \quad (1.3.8)$$

$$e^t {}_2F_2\left(\begin{array}{c} a+ix, a-ix \\ a+b, a+c \end{array} \bigg| -t\right) = \sum_{n=0}^{\infty} \frac{S_n(x^2;a,b,c)}{(a+b)_n(a+c)_n n!} t^n. \quad (1.3.9)$$

**References.** [45], [130], [155], [156], [168], [169].

## 1.4 Continuous Hahn

**Definition.**

$$p_n(x;a,b,c,d) = i^n \frac{(a+c)_n (a+d)_n}{n!} {}_3F_2\left(\begin{array}{c} -n, n+a+b+c+d-1, a+ix \\ a+c, a+d \end{array} \bigg| 1\right). \quad (1.4.1)$$

**Orthogonality.**

$$\frac{1}{2\pi} \int_{-\infty}^{\infty} \Gamma(a+ix)\Gamma(b+ix)\Gamma(c-ix)\Gamma(d-ix) p_m(x;a,b,c,d) p_n(x;a,b,c,d) dx$$
$$= \frac{\Gamma(n+a+c)\Gamma(n+a+d)\Gamma(n+b+c)\Gamma(n+b+d)}{(2n+a+b+c+d-1)\Gamma(n+a+b+c+d-1)n!} \delta_{mn}, \quad (1.4.2)$$

where
$$Re(a,b,c,d) > 0, \ c = \bar{a} \ \text{and} \ d = \bar{b}.$$

**Recurrence relation.**

$$(a+ix)\tilde{p}_n(x) = A_n \tilde{p}_{n+1}(x) - (A_n + C_n)\tilde{p}_n(x) + C_n \tilde{p}_{n-1}(x), \quad (1.4.3)$$

where
$$\tilde{p}_n(x) := \tilde{p}_n(x;a,b,c,d) = \frac{n!}{i^n (a+c)_n (a+d)_n} p_n(x;a,b,c,d)$$

and

$$\begin{cases} A_n = -\dfrac{(n+a+b+c+d-1)(n+a+c)(n+a+d)}{(2n+a+b+c+d-1)(2n+a+b+c+d)} \\ C_n = \dfrac{n(n+b+c-1)(n+b+d-1)}{(2n+a+b+c+d-2)(2n+a+b+c+d-1)}. \end{cases}$$

**Difference equation.**

$$n(n+a+b+c+d-1)y(x) = B(x)y(x+i) - [B(x)+D(x)]\,y(x) + D(x)y(x-i), \quad (1.4.4)$$

where

$$y(x) = p_n(x;a,b,c,d)$$

and

$$\begin{cases} B(x) = (c-ix)(d-ix) \\ D(x) = (a+ix)(b+ix). \end{cases}$$

**Generating functions.**

$$_2F_0\left(\begin{matrix} a+ix, b+ix \\ - \end{matrix} \middle| -it\right) {}_2F_0\left(\begin{matrix} c-ix, d-ix \\ - \end{matrix} \middle| it\right) \sim \sum_{n=0}^{\infty} p_n(x;a,b,c,d)t^n. \quad (1.4.5)$$

$$_1F_1\left(\begin{matrix} a+ix \\ a+c \end{matrix} \middle| -it\right) {}_1F_1\left(\begin{matrix} d-ix \\ b+d \end{matrix} \middle| it\right) = \sum_{n=0}^{\infty} \frac{p_n(x;a,b,c,d)}{(a+c)_n(b+d)_n} t^n. \quad (1.4.6)$$

$$_1F_1\left(\begin{matrix} a+ix \\ a+d \end{matrix} \middle| -it\right) {}_1F_1\left(\begin{matrix} c-ix \\ b+c \end{matrix} \middle| it\right) = \sum_{n=0}^{\infty} \frac{p_n(x;a,b,c,d)}{(a+d)_n(b+c)_n} t^n. \quad (1.4.7)$$

$$(1-t)^{1-a-b-c-d}\,{}_3F_2\left(\begin{matrix} \tfrac{1}{2}(a+b+c+d-1), \tfrac{1}{2}(a+b+c+d), a+ix \\ a+c, a+d \end{matrix} \middle| -\frac{4t}{(1-t)^2}\right)$$

$$= \sum_{n=0}^{\infty} \frac{(a+b+c+d-1)_n}{(a+c)_n(a+d)_n i^n} p_n(x;a,b,c,d)t^n. \quad (1.4.8)$$

**Remark.** Since the generating function (1.4.5) is divergent this relation must be seen as an equality in terms of formal power series.

**References.** [29], [31], [46], [51], [120], [156].

## 1.5 Hahn

**Definition.**

$$Q_n(x;\alpha,\beta,N) = {}_3\tilde{F}_2\left(\begin{matrix} -n, n+\alpha+\beta+1, -x \\ \alpha+1, -N \end{matrix} \middle| 1\right), \; n=0,1,2,\ldots,N. \quad (1.5.1)$$

**Orthogonality.**

$$\sum_{x=0}^{N} \binom{\alpha+x}{x}\binom{N+\beta-x}{N-x} Q_m(x;\alpha,\beta,N) Q_n(x;\alpha,\beta,N)$$

$$= \frac{(-1)^n n!(\beta+1)_n(n+\alpha+\beta+1)_{N+1}}{(N!)(2n+\alpha+\beta+1)(-N)_n(\alpha+1)_n}\delta_{mn}. \quad (1.5.2)$$

**Recurrence relation.**

$$-xQ_n(x) = A_n Q_{n+1}(x) - (A_n + C_n) Q_n(x) + C_n Q_{n-1}(x), \qquad (1.5.3)$$

where

$$Q_n(x) := Q_n(x; \alpha, \beta, N)$$

and

$$\begin{cases} A_n = \dfrac{(n+\alpha+\beta+1)(n+\alpha+1)(N-n)}{(2n+\alpha+\beta+1)(2n+\alpha+\beta+2)} \\ \\ C_n = \dfrac{n(n+\beta)(n+\alpha+\beta+N+1)}{(2n+\alpha+\beta)(2n+\alpha+\beta+1)}. \end{cases}$$

**Difference equation.**

$$n(n+\alpha+\beta+1)y(x) = B(x)y(x+1) - [B(x) + D(x)]y(x) + D(x)y(x-1), \qquad (1.5.4)$$

where

$$y(x) = Q_n(x; \alpha, \beta, N)$$

and

$$\begin{cases} B(x) = (x-N)(x+\alpha+1) \\ \\ D(x) = x(x-\beta-N-1). \end{cases}$$

**Generating functions.**

$${}_1F_1\left(\begin{matrix} x-N \\ \beta+1 \end{matrix} \middle| t\right) {}_1F_1\left(\begin{matrix} -x \\ \alpha+1 \end{matrix} \middle| -t\right) \simeq \sum_{n=0}^{N} \frac{(-N)_n}{(\beta+1)_n n!} Q_n(x; \alpha, \beta, N) t^n. \qquad (1.5.5)$$

$${}_1\tilde{F}_1\left(\begin{matrix} x-N \\ -N \end{matrix} \middle| t\right) {}_1F_1\left(\begin{matrix} N-x+\beta+1 \\ \alpha+\beta+N+2 \end{matrix} \middle| -t\right) \simeq \sum_{n=0}^{N} \frac{(\alpha+1)_n}{(\alpha+\beta+N+2)_n n!} Q_n(x; \alpha, \beta, N) t^n. \qquad (1.5.6)$$

$$(1-t)^{-\alpha-\beta-1} {}_3\tilde{F}_2\left(\begin{matrix} \frac{1}{2}(\alpha+\beta+1), \frac{1}{2}(\alpha+\beta+2), -x \\ \alpha+1, -N \end{matrix} \middle| -\frac{4t}{(1-t)^2}\right)$$

$$\simeq \sum_{n=0}^{N} \frac{(\alpha+\beta+1)_n}{n!} Q_n(x; \alpha, \beta, N) t^n. \qquad (1.5.7)$$

**Remarks.** If we interchange the role of $x$ and $n$ in (1.5.1) we obtain the dual Hahn polynomials defined by (1.6.1).

Since

$$Q_n(x; \alpha, \beta, N) = R_x(\lambda(n); \alpha, \beta, N)$$

we obtain the dual orthogonality relation for the Hahn polynomials from the orthogonality relation (1.6.2) of the dual Hahn polynomials :

$$\sum_{n=0}^{N} \frac{(N!)(-N)_n(\alpha+1)_n(2n+\alpha+\beta+1)}{(-1)^n n!(\beta+1)_n(n+\alpha+\beta+1)_{N+1}} Q_n(x; \alpha, \beta, N) Q_n(y; \alpha, \beta, N)$$

$$= \frac{\delta_{xy}}{\binom{\alpha+x}{x}\binom{N+\beta-x}{N-x}}, \ x, y \in \{0, 1, 2, \ldots, N\}.$$

For $x = 0, 1, 2, \ldots, N$ the generating function (1.5.5) can also be written as :

$${}_1F_1\left(\begin{matrix} x-N \\ \beta+1 \end{matrix} \middle| t\right) {}_1F_1\left(\begin{matrix} -x \\ \alpha+1 \end{matrix} \middle| -t\right) = \sum_{n=0}^{N} \frac{(-N)_n}{(\beta+1)_n n!} Q_n(x; \alpha, \beta, N) t^n.$$

**References.** [10], [25], [27], [31], [45], [47], [77], [80], [82], [87], [89], [104], [105], [124], [144], [153], [156], [165], [167], [180], [187], [189], [190], [195], [217], [218], [227].

## 1.6 Dual Hahn

**Definition.**

$$R_n(\lambda(x); \gamma, \delta, N) = {}_3\tilde{F}_2\left(\begin{array}{c} -n, -x, x+\gamma+\delta+1 \\ \gamma+1, -N \end{array}\bigg| 1\right), \quad n = 0, 1, 2, \ldots, N, \tag{1.6.1}$$

where

$$\lambda(x) = x(x+\gamma+\delta+1).$$

**Orthogonality.**

$$\sum_{x=0}^{N} \frac{(N!)(-N)_x(\gamma+1)_x(2x+\gamma+\delta+1)}{(-1)^x(x!)(\delta+1)_x(x+\gamma+\delta+1)_{N+1}} R_m(\lambda(x); \gamma, \delta, N) R_n(\lambda(x); \gamma, \delta, N)$$

$$= \frac{\delta_{mn}}{\binom{\gamma+n}{n}\binom{N+\delta-n}{N-n}}. \tag{1.6.2}$$

**Recurrence relation.**

$$\lambda(x) R_n(\lambda(x)) = A_n R_{n+1}(\lambda(x)) - (A_n + C_n) R_n(\lambda(x)) + C_n R_{n-1}(\lambda(x)), \tag{1.6.3}$$

where

$$R_n(\lambda(x)) := R_n(\lambda(x); \gamma, \delta, N)$$

and

$$\begin{cases} A_n = (n-N)(n+\gamma+1) \\ C_n = n(n-\delta-N-1). \end{cases}$$

**Difference equation.**

$$-ny(x) = B(x)y(x+1) - [B(x) + D(x)]y(x) + D(x)y(x-1), \quad y(x) = R_n(\lambda(x); \gamma, \delta, N), \tag{1.6.4}$$

where

$$\begin{cases} B(x) = \dfrac{(x+\gamma+\delta+1)(x+\gamma+1)(N-x)}{(2x+\gamma+\delta+1)(2x+\gamma+\delta+2)} \\ D(x) = \dfrac{x(x+\delta)(x+\gamma+\delta+N+1)}{(2x+\gamma+\delta)(2x+\gamma+\delta+1)}. \end{cases}$$

**Generating functions.**

$$(1-t)^{x+\delta} {}_2\tilde{F}_1\left(\begin{array}{c} x-N, x+\gamma+\delta+1 \\ -N \end{array}\bigg| t\right) \simeq \sum_{n=0}^{N} \frac{(\gamma+1)_n}{n!} R_n(\lambda(x); \gamma, \delta, N) t^n. \tag{1.6.5}$$

$$(1-t)^{N-x} {}_2F_1\left(\begin{array}{c} -x, -x-\delta \\ \gamma+1 \end{array}\bigg| t\right) \simeq \sum_{n=0}^{N} \frac{(-N)_n}{n!} R_n(\lambda(x); \gamma, \delta, N) t^n. \tag{1.6.6}$$

$$(1-t)^x {}_2F_1\left(\begin{array}{c} x-N, x+\gamma+1 \\ -\delta-N \end{array}\bigg| t\right) \simeq \sum_{n=0}^{N} \frac{(-N)_n(\gamma+1)_n}{(-\delta-N)_n n!} R_n(\lambda(x); \gamma, \delta, N) t^n. \tag{1.6.7}$$

$$e^t {}_2\tilde{F}_2\left(\begin{array}{c} -x, x+\gamma+\delta+1 \\ \gamma+1, -N \end{array}\bigg| -t\right) \simeq \sum_{n=0}^{N} \frac{R_n(\lambda(x); \gamma, \delta, N)}{n!} t^n. \tag{1.6.8}$$

**Remarks.** If we interchange the role of $x$ and $n$ in the definition (1.6.1) of the dual Hahn polynomials we obtain the Hahn polynomials defined by (1.5.1).

Since
$$R_n(\lambda(x); \gamma, \delta, N) = Q_x(n; \gamma, \delta, N)$$
we obtain the dual orthogonality relation for the dual Hahn polynomials from the orthogonality relation (1.5.2) for the Hahn polynomials :

$$\sum_{n=0}^{N} \binom{\gamma+n}{n}\binom{N+\delta-n}{N-n} R_n(\lambda(x); \gamma, \delta, N) R_n(\lambda(y); \gamma, \delta, N)$$
$$= \frac{(-1)^x (x!)(\delta+1)_x (x+\gamma+\delta+1)_{N+1}}{(N!)(-N)_x (\gamma+1)_x (2x+\gamma+\delta+1)} \delta_{xy}, \ x, y \in \{0, 1, 2, \ldots, N\}.$$

For $x = 0, 1, 2, \ldots, N$ the generating function (1.6.6) can also be written as :

$$(1-t)^{N-x} {}_2F_1\left(\begin{matrix}-x, -x-\delta \\ \gamma+1\end{matrix}\bigg| t\right) = \sum_{n=0}^{N} \frac{(-N)_n}{n!} R_n(\lambda(x); \gamma, \delta, N) t^n.$$

For $x = 0, 1, 2, \ldots, N$ the generating function (1.6.7) can also be written as :

$$(1-t)^x {}_2F_1\left(\begin{matrix}x-N, x+\gamma+1 \\ -\delta-N\end{matrix}\bigg| t\right) = \sum_{n=0}^{N} \frac{(-N)_n (\gamma+1)_n}{(-\delta-N)_n n!} R_n(\lambda(x); \gamma, \delta, N) t^n.$$

**References.** [45], [47], [144], [153], [156], [180], [194], [217], [227].

## 1.7 Meixner-Pollaczek

**Definition.**
$$P_n^{(\lambda)}(x; \phi) = \frac{(2\lambda)_n}{n!} e^{in\phi} {}_2F_1\left(\begin{matrix}-n, \lambda+ix \\ 2\lambda\end{matrix}\bigg| 1-e^{-2i\phi}\right). \tag{1.7.1}$$

**Orthogonality.**
$$\frac{1}{2\pi} \int_{-\infty}^{\infty} e^{(2\phi-\pi)x} |\Gamma(\lambda+ix)|^2 P_m^{(\lambda)}(x;\phi) P_n^{(\lambda)}(x;\phi) dx$$
$$= \frac{\Gamma(n+2\lambda)}{(2\sin\phi)^{2\lambda} n!} \delta_{mn}, \ \lambda > 0 \ \text{and} \ 0 < \phi < \pi. \tag{1.7.2}$$

**Recurrence relation.**
$$(n+1)P_{n+1}^{(\lambda)}(x;\phi) - 2\left[x\sin\phi + (n+\lambda)\cos\phi\right] P_n^{(\lambda)}(x;\phi) + (n+2\lambda-1)P_{n-1}^{(\lambda)}(x;\phi) = 0. \tag{1.7.3}$$

**Difference equation.**
$$e^{i\phi}(\lambda-ix)y(x+i) + 2i\left[x\cos\phi - (n+\lambda)\sin\phi\right]y(x) - e^{-i\phi}(\lambda+ix)y(x-i) = 0, \tag{1.7.4}$$
where
$$y(x) = P_n^{(\lambda)}(x;\phi).$$

**Generating functions.**
$$(1-e^{i\phi}t)^{-\lambda+ix}(1-e^{-i\phi}t)^{-\lambda-ix} = \sum_{n=0}^{\infty} P_n^{(\lambda)}(x;\phi) t^n. \tag{1.7.5}$$

$$e^t {}_1F_1\left(\begin{matrix}\lambda+ix \\ 2\lambda\end{matrix}\bigg| (e^{-2i\phi}-1)t\right) = \sum_{n=0}^{\infty} \frac{P_n^{(\lambda)}(x;\phi)}{(2\lambda)_n e^{in\phi}} t^n. \tag{1.7.6}$$

**References.** [10], [15], [25], [31], [45], [47], [72], [77], [126], [156], [174], [180], [189], [227], [230].

## 1.8 Jacobi

**Definition.**
$$P_n^{(\alpha,\beta)}(x) = \frac{(\alpha+1)_n}{n!} {}_2F_1\left(\begin{array}{c} -n, n+\alpha+\beta+1 \\ \alpha+1 \end{array} \bigg| \frac{1-x}{2}\right). \tag{1.8.1}$$

**Orthogonality.**

$$\int_{-1}^{1} (1-x)^\alpha (1+x)^\beta P_m^{(\alpha,\beta)}(x) P_n^{(\alpha,\beta)}(x) dx$$
$$= \frac{2^{\alpha+\beta+1}}{2n+\alpha+\beta+1} \frac{\Gamma(n+\alpha+1)\Gamma(n+\beta+1)}{n!\,\Gamma(n+\alpha+\beta+1)} \delta_{mn}, \quad \alpha > -1 \text{ and } \beta > -1. \tag{1.8.2}$$

**Recurrence relation.**

$$xP_n^{(\alpha,\beta)}(x) = \frac{2(n+1)(n+\alpha+\beta+1)}{(2n+\alpha+\beta+1)(2n+\alpha+\beta+2)} P_{n+1}^{(\alpha,\beta)}(x) +$$
$$+ \frac{\beta^2 - \alpha^2}{(2n+\alpha+\beta)(2n+\alpha+\beta+2)} P_n^{(\alpha,\beta)}(x) + \frac{2(n+\alpha)(n+\beta)}{(2n+\alpha+\beta)(2n+\alpha+\beta+1)} P_{n-1}^{(\alpha,\beta)}(x). \tag{1.8.3}$$

**Differential equation.**

$$(1-x^2)y''(x) + [\beta - \alpha - (\alpha+\beta+2)x] y'(x) + n(n+\alpha+\beta+1)y(x) = 0, \ y(x) = P_n^{(\alpha,\beta)}(x). \tag{1.8.4}$$

**Generating functions.**

$$\frac{2^{\alpha+\beta}}{R(1+R-t)^\alpha (1+R+t)^\beta} = \sum_{n=0}^{\infty} P_n^{(\alpha,\beta)}(x) t^n, \ R = \sqrt{1-2xt+t^2}. \tag{1.8.5}$$

$${}_0F_1\left(\begin{array}{c} - \\ \alpha+1 \end{array} \bigg| \frac{(x-1)t}{2}\right) {}_0F_1\left(\begin{array}{c} - \\ \beta+1 \end{array} \bigg| \frac{(x+1)t}{2}\right) = \sum_{n=0}^{\infty} \frac{P_n^{(\alpha,\beta)}(x)}{(\alpha+1)_n(\beta+1)_n} t^n. \tag{1.8.6}$$

$$(1-t)^{-\alpha-\beta-1} {}_2F_1\left(\begin{array}{c} \frac{1}{2}(\alpha+\beta+1), \frac{1}{2}(\alpha+\beta+2) \\ \alpha+1 \end{array} \bigg| \frac{2(x-1)t}{(1-t)^2}\right)$$
$$= \sum_{n=0}^{\infty} \frac{(\alpha+\beta+1)_n}{(\alpha+1)_n} P_n^{(\alpha,\beta)}(x) t^n. \tag{1.8.7}$$

$$(1+t)^{-\alpha-\beta-1} {}_2F_1\left(\begin{array}{c} \frac{1}{2}(\alpha+\beta+1), \frac{1}{2}(\alpha+\beta+2) \\ \beta+1 \end{array} \bigg| \frac{2(x+1)t}{(1+t)^2}\right)$$
$$= \sum_{n=0}^{\infty} \frac{(\alpha+\beta+1)_n}{(\beta+1)_n} P_n^{(\alpha,\beta)}(x) t^n. \tag{1.8.8}$$

$${}_2F_1\left(\begin{array}{c} \gamma, \alpha+\beta+1-\gamma \\ \alpha+1 \end{array} \bigg| \frac{1-R-t}{2}\right) {}_2F_1\left(\begin{array}{c} \gamma, \alpha+\beta+1-\gamma \\ \beta+1 \end{array} \bigg| \frac{1-R+t}{2}\right)$$
$$= \sum_{n=0}^{\infty} \frac{(\gamma)_n (\alpha+\beta+1-\gamma)_n}{(\alpha+1)_n (\beta+1)_n} P_n^{(\alpha,\beta)}(x) t^n, \ R = \sqrt{1-2xt+t^2}, \ \gamma \text{ arbitrary.} \tag{1.8.9}$$

**Remarks.** The Jacobi polynomials defined by (1.8.1) and the Meixner polynomials given by (1.9.1) are related in the following way :

$$\frac{(\beta)_n}{n!} M_n(x;\beta,c) = P_n^{(\beta-1,-n-\beta-x)}\left(\frac{2-c}{c}\right).$$

The Jacobi polynomials are also related to the Gegenbauer (or ultraspherical) polynomials defined by (1.8.10) by the quadratic transformations :

$$C_{2n}^{(\lambda)}(x) = \frac{(\lambda)_n}{(\frac{1}{2})_n} P_n^{(\lambda-\frac{1}{2},-\frac{1}{2})}(2x^2 - 1)$$

and

$$C_{2n+1}^{(\lambda)}(x) = \frac{(\lambda)_{n+1}}{(\frac{1}{2})_{n+1}} x P_n^{(\lambda-\frac{1}{2},\frac{1}{2})}(2x^2 - 1).$$

**References.** [2], [3], [7], [9], [14], [25], [26], [27], [28], [31], [33], [34], [35], [42], [45], [50], [55], [61], [68], [69], [77], [80], [82], [90], [93], [94], [98], [99], [100], [101], [102], [103], [104], [105], [106], [115], [119], [122], [123], [124], [134], [146], [147], [148], [149], [151], [152], [155], [156], [165], [170], [173], [176], [178], [180], [182], [184], [185], [186], [188], [193], [203], [205], [211], [213], [214], [218], [220], [229], [232].

# Special cases

## 1.8.1 Gegenbauer / Ultraspherical

**Definition.** The Gegenbauer (or ultraspherical) polynomials are Jacobi polynomials with $\alpha = \beta = \lambda - \frac{1}{2}$ and another normalization :

$$C_n^{(\lambda)}(x) = \frac{(2\lambda)_n}{(\lambda+\frac{1}{2})_n} P_n^{(\lambda-\frac{1}{2},\lambda-\frac{1}{2})}(x) = \frac{(2\lambda)_n}{n!} {}_2F_1\left(\begin{array}{c} -n, n+2\lambda \\ \lambda + \frac{1}{2} \end{array} \middle| \frac{1-x}{2}\right), \ \lambda \neq 0. \tag{1.8.10}$$

**Orthogonality.**

$$\int_{-1}^{1}(1-x^2)^{\lambda-\frac{1}{2}} C_m^{(\lambda)}(x) C_n^{(\lambda)}(x) dx = \frac{\pi \Gamma(n+2\lambda) 2^{1-2\lambda}}{\{\Gamma(\lambda)\}^2 (n+\lambda)n!} \delta_{mn}, \ \lambda > -\frac{1}{2} \ \text{and} \ \lambda \neq 0. \tag{1.8.11}$$

**Recurrence relation.**

$$2(n+\lambda)x C_n^{(\lambda)}(x) = (n+1) C_{n+1}^{(\lambda)}(x) + (n+2\lambda-1) C_{n-1}^{(\lambda)}(x). \tag{1.8.12}$$

**Differential equation.**

$$(1-x^2)y''(x) - (2\lambda+1)xy'(x) + n(n+2\lambda)y(x) = 0, \ y(x) = C_n^{(\lambda)}(x). \tag{1.8.13}$$

**Generating functions.**

$$(1 - 2xt + t^2)^{-\lambda} = \sum_{n=0}^{\infty} C_n^{(\lambda)}(x) t^n. \tag{1.8.14}$$

$$R^{-1}\left(\frac{1+R-xt}{2}\right)^{\frac{1}{2}-\lambda} = \sum_{n=0}^{\infty} \frac{(\lambda+\frac{1}{2})_n}{(2\lambda)_n} C_n^{(\lambda)}(x) t^n, \ R = \sqrt{1-2xt+t^2}. \tag{1.8.15}$$

$${}_0F_1\left(\begin{array}{c} - \\ \lambda + \frac{1}{2} \end{array} \middle| \frac{(x-1)t}{2}\right) {}_0F_1\left(\begin{array}{c} - \\ \lambda + \frac{1}{2} \end{array} \middle| \frac{(x+1)t}{2}\right) = \sum_{n=0}^{\infty} \frac{C_n^{(\lambda)}(x)}{(2\lambda)_n (\lambda+\frac{1}{2})_n} t^n. \tag{1.8.16}$$

$$e^{xt} {}_0F_1\left(\begin{array}{c} - \\ \lambda + \frac{1}{2} \end{array} \middle| \frac{(x^2-1)t^2}{4}\right) = \sum_{n=0}^{\infty} \frac{C_n^{(\lambda)}(x)}{(2\lambda)_n} t^n. \tag{1.8.17}$$

$$\displaylines{{}_2F_1\left(\begin{array}{c}\gamma,2\lambda-\gamma\\ \lambda+\frac{1}{2}\end{array}\bigg|\frac{1-R-t}{2}\right){}_2F_1\left(\begin{array}{c}\gamma,2\lambda-\gamma\\ \lambda+\frac{1}{2}\end{array}\bigg|\frac{1-R+t}{2}\right)\hfill\cr\hfill=\sum_{n=0}^{\infty}\frac{(\gamma)_n(2\lambda-\gamma)_n}{(2\lambda)_n(\lambda+\frac{1}{2})_n}C_n^{(\lambda)}(x)t^n,\ R=\sqrt{1-2xt+t^2},\ \gamma\text{ arbitrary}.\quad(1.8.18)}$$

$$(1-xt)^{-\gamma}{}_2F_1\left(\begin{array}{c}\frac{1}{2}\gamma,\frac{1}{2}\gamma+\frac{1}{2}\\ \lambda+\frac{1}{2}\end{array}\bigg|\frac{(x^2-1)t^2}{(1-xt)^2}\right)=\sum_{n=0}^{\infty}\frac{(\gamma)_n}{(2\lambda)_n}C_n^{(\lambda)}(x)t^n,\ \gamma\text{ arbitrary}.\quad(1.8.19)$$

**Remarks.** The case $\lambda=0$ needs another normalization. In that case we have the Chebyshev polynomials of the first kind described in the next subsection.

The Gegenbauer (or ultraspherical) polynomials defined by (1.8.10) and the Jacobi polynomials given by (1.8.1) are related by the quadratic transformations :

$$C_{2n}^{(\lambda)}(x)=\frac{(\lambda)_n}{(\frac{1}{2})_n}P_n^{(\lambda-\frac{1}{2},-\frac{1}{2})}(2x^2-1)$$

and

$$C_{2n+1}^{(\lambda)}(x)=\frac{(\lambda)_{n+1}}{(\frac{1}{2})_{n+1}}xP_n^{(\lambda-\frac{1}{2},\frac{1}{2})}(2x^2-1).$$

**References.** [2], [4], [26], [27], [31], [33], [41], [54], [55], [56], [61], [62], [63], [64], [65], [66], [68], [77], [83], [92], [94], [103], [108], [119], [123], [156], [164], [180], [203], [205], [220], [223], [232].

### 1.8.2 Chebyshev

**Definitions.** The Chebyshev polynomials of the first kind can be obtained from the Jacobi polynomials by taking $\alpha=\beta=-\frac{1}{2}$ :

$$T_n(x)=\frac{P_n^{(-\frac{1}{2},-\frac{1}{2})}(x)}{P_n^{(-\frac{1}{2},-\frac{1}{2})}(1)}={}_2F_1\left(\begin{array}{c}-n,n\\ \frac{1}{2}\end{array}\bigg|\frac{1-x}{2}\right)\quad(1.8.20)$$

and the Chebyshev polynomials of the second kind can be obtained from the Jacobi polynomials by taking $\alpha=\beta=\frac{1}{2}$ :

$$U_n(x)=(n+1)\frac{P_n^{(\frac{1}{2},\frac{1}{2})}(x)}{P_n^{(\frac{1}{2},\frac{1}{2})}(1)}=(n+1){}_2F_1\left(\begin{array}{c}-n,n+2\\ \frac{3}{2}\end{array}\bigg|\frac{1-x}{2}\right).\quad(1.8.21)$$

**Orthogonality.**

$$\int_{-1}^{1}(1-x^2)^{-\frac{1}{2}}T_m(x)T_n(x)dx=\begin{cases}\frac{\pi}{2}\delta_{mn},&n\neq 0\\ \pi\delta_{mn},&n=0.\end{cases}\quad(1.8.22)$$

$$\int_{-1}^{1}(1-x^2)^{\frac{1}{2}}U_m(x)U_n(x)dx=\frac{\pi}{2}\delta_{mn}.\quad(1.8.23)$$

**Recurrence relations.**

$$2xT_n(x)=T_{n+1}(x)+T_{n-1}(x),\ T_0(x)=1\text{ and }T_1(x)=x.\quad(1.8.24)$$

$$2xU_n(x) = U_{n+1}(x) + U_{n-1}(x). \tag{1.8.25}$$

**Differential equations.**

$$(1-x^2)y''(x) - xy'(x) + n^2 y(x) = 0, \ y(x) = T_n(x). \tag{1.8.26}$$

$$(1-x^2)y''(x) - 3xy'(x) + n(n+2)y(x) = 0, \ y(x) = U_n(x). \tag{1.8.27}$$

**Generating functions.**

$$\frac{1-xt}{1-2xt+t^2} = \sum_{n=0}^{\infty} T_n(x)t^n. \tag{1.8.28}$$

$$R^{-1}\sqrt{\frac{1}{2}(1+R-xt)} = \sum_{n=0}^{\infty} \frac{\left(\frac{1}{2}\right)_n}{n!} T_n(x)t^n, \ R = \sqrt{1-2xt+t^2}. \tag{1.8.29}$$

$$_0F_1\left(\begin{array}{c}-\\ \frac{1}{2}\end{array}\bigg|\ \frac{(x-1)t}{2}\right) {}_0F_1\left(\begin{array}{c}-\\ \frac{1}{2}\end{array}\bigg|\ \frac{(x+1)t}{2}\right) = \sum_{n=0}^{\infty} \frac{T_n(x)}{\left(\frac{1}{2}\right)_n n!} t^n. \tag{1.8.30}$$

$$e^{xt}{}_0F_1\left(\begin{array}{c}-\\ \frac{1}{2}\end{array}\bigg|\ \frac{(x^2-1)t^2}{4}\right) = \sum_{n=0}^{\infty} \frac{T_n(x)}{n!} t^n. \tag{1.8.31}$$

$$_2F_1\left(\begin{array}{c}\gamma, -\gamma\\ \frac{1}{2}\end{array}\bigg|\ \frac{1-R-t}{2}\right) {}_2F_1\left(\begin{array}{c}\gamma, -\gamma\\ \frac{1}{2}\end{array}\bigg|\ \frac{1-R+t}{2}\right)$$
$$= \sum_{n=0}^{\infty} \frac{(\gamma)_n(-\gamma)_n}{\left(\frac{1}{2}\right)_n n!} T_n(x)t^n, \ R = \sqrt{1-2xt+t^2}, \ \gamma \text{ arbitrary.} \tag{1.8.32}$$

$$(1-xt)^{-\gamma} {}_2F_1\left(\begin{array}{c}\frac{1}{2}\gamma, \frac{1}{2}\gamma + \frac{1}{2}\\ \frac{1}{2}\end{array}\bigg|\ \frac{(x^2-1)t^2}{(1-xt)^2}\right) = \sum_{n=0}^{\infty} \frac{(\gamma)_n}{n!} T_n(x)t^n, \ \gamma \text{ arbitrary.} \tag{1.8.33}$$

$$\frac{1}{1-2xt+t^2} = \sum_{n=0}^{\infty} U_n(x)t^n. \tag{1.8.34}$$

$$\frac{1}{R\sqrt{\frac{1}{2}(1+R-xt)}} = \sum_{n=0}^{\infty} \frac{\left(\frac{3}{2}\right)_n}{(n+1)!} U_n(x)t^n, \ R = \sqrt{1-2xt+t^2}. \tag{1.8.35}$$

$$_0F_1\left(\begin{array}{c}-\\ \frac{3}{2}\end{array}\bigg|\ \frac{(x-1)t}{2}\right) {}_0F_1\left(\begin{array}{c}-\\ \frac{3}{2}\end{array}\bigg|\ \frac{(x+1)t}{2}\right) = \sum_{n=0}^{\infty} \frac{U_n(x)}{\left(\frac{3}{2}\right)_n (n+1)!} t^n. \tag{1.8.36}$$

$$e^{xt}{}_0F_1\left(\begin{array}{c}-\\ \frac{3}{2}\end{array}\bigg|\ \frac{(x^2-1)t^2}{4}\right) = \sum_{n=0}^{\infty} \frac{U_n(x)}{(n+1)!} t^n. \tag{1.8.37}$$

$$_2F_1\left(\begin{array}{c}\gamma, 2-\gamma\\ \frac{3}{2}\end{array}\bigg|\ \frac{1-R-t}{2}\right) {}_2F_1\left(\begin{array}{c}\gamma, 2-\gamma\\ \frac{3}{2}\end{array}\bigg|\ \frac{1-R+t}{2}\right)$$
$$= \sum_{n=0}^{\infty} \frac{(\gamma)_n(2-\gamma)_n}{\left(\frac{3}{2}\right)_n (n+1)!} U_n(x)t^n, \ R = \sqrt{1-2xt+t^2}, \ \gamma \text{ arbitrary.} \tag{1.8.38}$$

$$(1-xt)^{-\gamma} {}_2F_1\left(\begin{array}{c}\frac{1}{2}\gamma, \frac{1}{2}\gamma + \frac{1}{2}\\ \frac{3}{2}\end{array}\bigg|\ \frac{(x^2-1)t^2}{(1-xt)^2}\right) = \sum_{n=0}^{\infty} \frac{(\gamma)_n}{(n+1)!} U_n(x)t^n, \ \gamma \text{ arbitrary.} \tag{1.8.39}$$

**Remarks.** The Chebyshev polynomials can also be written as :

$$T_n(x) = \cos(n \arccos x)$$

and

$$U_n(x) = \frac{\sin(n+1)\theta}{\sin \theta}, \ x = \cos \theta.$$

Further we have

$$U_n(x) = C_n^{(1)}(x)$$

where $C_n^{(\lambda)}(x)$ denotes the Gegenbauer (or ultraspherical) polynomial defined by (1.8.10) in the preceding subsection.

**References.** [2], [33], [36], [77], [83], [94], [119], [170], [180], [205], [206], [210], [220], [232].

### 1.8.3 Legendre / Spherical

**Definition.** The Legendre (or spherical) polynomials are Jacobi polynomials with $\alpha = \beta = 0$ :

$$P_n(x) = P_n^{(0,0)}(x) = {}_2F_1\left(\begin{array}{c}-n, n+1\\1\end{array}\bigg|\frac{1-x}{2}\right). \tag{1.8.40}$$

**Orthogonality.**

$$\int_{-1}^{1} P_m(x)P_n(x)dx = \frac{2}{2n+1}\delta_{mn}. \tag{1.8.41}$$

**Recurrence relation.**

$$(2n+1)xP_n(x) = (n+1)P_{n+1}(x) + nP_{n-1}(x). \tag{1.8.42}$$

**Differential equation.**

$$(1-x^2)y''(x) - 2xy'(x) + n(n+1)y(x) = 0, \ y(x) = P_n(x). \tag{1.8.43}$$

**Generating functions.**

$$\frac{1}{\sqrt{1-2xt+t^2}} = \sum_{n=0}^{\infty} P_n(x)t^n. \tag{1.8.44}$$

$${}_0F_1\left(\begin{array}{c}-\\1\end{array}\bigg|\frac{(x-1)t}{2}\right) {}_0F_1\left(\begin{array}{c}-\\1\end{array}\bigg|\frac{(x+1)t}{2}\right) = \sum_{n=0}^{\infty} \frac{P_n(x)}{n!^2}t^n. \tag{1.8.45}$$

$$e^{xt} {}_0F_1\left(\begin{array}{c}-\\1\end{array}\bigg|\frac{(x^2-1)t^2}{4}\right) = \sum_{n=0}^{\infty} \frac{P_n(x)}{n!}t^n. \tag{1.8.46}$$

$${}_2F_1\left(\begin{array}{c}\gamma, 1-\gamma\\1\end{array}\bigg|\frac{1-R-t}{2}\right) {}_2F_1\left(\begin{array}{c}\gamma, 1-\gamma\\1\end{array}\bigg|\frac{1-R+t}{2}\right)$$
$$= \sum_{n=0}^{\infty} \frac{(\gamma)_n(1-\gamma)_n}{n!^2}P_n(x)t^n, \ R = \sqrt{1-2xt+t^2}, \ \gamma \text{ arbitrary.} \tag{1.8.47}$$

$$(1-xt)^{-\gamma} {}_2F_1\left(\begin{array}{c}\frac{1}{2}\gamma, \frac{1}{2}\gamma+\frac{1}{2}\\1\end{array}\bigg|\frac{(x^2-1)t^2}{(1-xt)^2}\right) = \sum_{n=0}^{\infty} \frac{(\gamma)_n}{n!}P_n(x)t^n, \ \gamma \text{ arbitrary.} \tag{1.8.48}$$

**References.** [2], [5], [10], [53], [55], [77], [83], [94], [119], [180], [182], [205], [220], [232].

## 1.9 Meixner

**Definition.**
$$M_n(x;\beta,c) = {}_2F_1\left(\begin{array}{c}-n,-x\\ \beta\end{array}\bigg|\, 1-\frac{1}{c}\right). \tag{1.9.1}$$

**Orthogonality.**
$$\sum_{x=0}^{\infty}\frac{(\beta)_x}{x!}c^x M_m(x;\beta,c)M_n(x;\beta,c) = \frac{c^{-n}n!}{(\beta)_n(1-c)^{\beta}}\delta_{mn},\ \beta>0\ \text{and}\ 0<c<1. \tag{1.9.2}$$

**Recurrence relation.**
$$(c-1)xM_n(x;\beta,c) = c(n+\beta)M_{n+1}(x;\beta,c) + \\ -[n+(n+\beta)c]\,M_n(x;\beta,c) + nM_{n-1}(x;\beta,c). \tag{1.9.3}$$

**Difference equation.**
$$n(c-1)y(x) = c(x+\beta)y(x+1) - [x+(x+\beta)c]\,y(x) + xy(x-1),\ y(x)=M_n(x;\beta,c). \tag{1.9.4}$$

**Generating functions.**
$$\left(1-\frac{t}{c}\right)^x (1-t)^{-x-\beta} = \sum_{n=0}^{\infty}\frac{(\beta)_n}{n!}M_n(x;\beta,c)t^n. \tag{1.9.5}$$

$$e^t\,{}_1F_1\left(\begin{array}{c}-x\\ \beta\end{array}\bigg|\,\left(\frac{1-c}{c}\right)t\right) = \sum_{n=0}^{\infty}\frac{M_n(x;\beta,c)}{n!}t^n. \tag{1.9.6}$$

**Remarks.** The Meixner polynomials defined by (1.9.1) and the Jacobi polynomials given by (1.8.1) are related in the following way :
$$\frac{(\beta)_n}{n!}M_n(x;\beta,c) = P_n^{(\beta-1,-n-\beta-x)}\left(\frac{2-c}{c}\right).$$

The Meixner polynomials are also related to the Krawtchouk polynomials defined by (1.10.1) in the following way :
$$K_n(x;p,N) = M_n\left(x;-N,\frac{p}{p-1}\right).$$

**References.** [7], [10], [15], [17], [25], [27], [31], [36], [45], [47], [77], [82], [94], [97], [104], [105], [129], [135], [143], [156], [174], [180], [189], [222], [233].

## 1.10 Krawtchouk

**Definition.**
$$K_n(x;p,N) = {}_2\tilde{F}_1\left(\begin{array}{c}-n,-x\\ -N\end{array}\bigg|\,\frac{1}{p}\right),\ n=0,1,2,\ldots,N. \tag{1.10.1}$$

**Orthogonality.**
$$\sum_{x=0}^{N}\binom{N}{x}p^x(1-p)^{N-x}K_m(x;p,N)K_n(x;p,N) = \frac{(-1)^n n!}{(-N)_n}\left(\frac{1-p}{p}\right)^n\delta_{mn},\ 0<p<1. \tag{1.10.2}$$

**Recurrence relation.**

$$-xK_n(x;p,N) = p(N-n)K_{n+1}(x;p,N) +$$
$$- [p(N-n) + n(1-p)] K_n(x;p,N) + n(1-p)K_{n-1}(x;p,N). \qquad (1.10.3)$$

**Difference equation.**

$$-ny(x) = p(N-x)y(x+1) - [p(N-x) + x(1-p)] y(x) + x(1-p)y(x-1), \qquad (1.10.4)$$

where

$$y(x) = K_n(x;p,N).$$

**Generating functions.**

$$\left(1 - \frac{(1-p)}{p}t\right)^x (1+t)^{N-x} \simeq \sum_{n=0}^{N} \binom{N}{n} K_n(x;p,N)t^n. \qquad (1.10.5)$$

$$e^t {}_1\tilde{F}_1 \left( \begin{matrix} -x \\ -N \end{matrix} \middle| -\frac{t}{p} \right) \simeq \sum_{n=0}^{N} \frac{K_n(x;p,N)}{n!}t^n. \qquad (1.10.6)$$

**Remarks.** The Krawtchouk polynomials are self-dual, which means that

$$K_n(x;p,N) = K_x(n;p,N), \; n, x \in \{0,1,2,\ldots,N\}.$$

By using this relation we easily obtain the so-called dual orthogonality relation from the orthogonality relation (1.10.2):

$$\sum_{n=0}^{N} \binom{N}{n} p^n (1-p)^{N-n} K_n(x;p,N) K_n(y;p,N) = \frac{\left(\frac{1-p}{p}\right)^x}{\binom{N}{x}} \delta_{xy},$$

where $0 < p < 1$ and $x, y \in \{0,1,2,\ldots,N\}$.

The Krawtchouk polynomials are related to the Meixner polynomials defined by (1.9.1) in the following way:

$$K_n(x;p,N) = M_n \left( x; -N, \frac{p}{p-1} \right).$$

For $x = 0, 1, 2, \ldots, N$ the generating function (1.10.5) can also be written as:

$$\left(1 - \frac{(1-p)}{p}t\right)^x (1+t)^{N-x} = \sum_{n=0}^{N} \binom{N}{n} K_n(x;p,N)t^n.$$

**References.** [10], [25], [27], [31], [45], [77], [87], [90], [91], [94], [104], [105], [143], [154], [156], [167], [180], [189], [191], [217], [218], [220], [233].

## 1.11 Laguerre

**Definition.**

$$L_n^{(\alpha)}(x) = \frac{(\alpha+1)_n}{n!} {}_1F_1 \left( \begin{matrix} -n \\ \alpha+1 \end{matrix} \middle| x \right). \qquad (1.11.1)$$

**Orthogonality.**

$$\int_0^\infty x^\alpha e^{-x} L_m^{(\alpha)}(x) L_n^{(\alpha)}(x) dx = \frac{\Gamma(n+\alpha+1)}{n!} \delta_{mn}, \; \alpha > -1. \qquad (1.11.2)$$

**Recurrence relation.**

$$(n+1)L_{n+1}^{(\alpha)}(x) - (2n+\alpha+1-x)L_n^{(\alpha)}(x) + (n+\alpha)L_{n-1}^{(\alpha)}(x) = 0. \tag{1.11.3}$$

**Differential equation.**

$$xy''(x) + (\alpha+1-x)y'(x) + ny(x) = 0, \ y(x) = L_n^{(\alpha)}(x). \tag{1.11.4}$$

**Generating functions.**

$$(1-t)^{-\alpha-1}\exp\left(\frac{xt}{t-1}\right) = \sum_{n=0}^{\infty} L_n^{(\alpha)}(x)t^n. \tag{1.11.5}$$

$$e^t {}_0F_1\left(\begin{array}{c}-\\ \alpha+1\end{array}\bigg|-xt\right) = \sum_{n=0}^{\infty} \frac{L_n^{(\alpha)}(x)}{(\alpha+1)_n} t^n. \tag{1.11.6}$$

$$(1-t)^{-\gamma} {}_1F_1\left(\begin{array}{c}\gamma\\ \alpha+1\end{array}\bigg|\frac{xt}{t-1}\right) = \sum_{n=0}^{\infty} \frac{(\gamma)_n}{(\alpha+1)_n} L_n^{(\alpha)}(x) t^n, \ \gamma \text{ arbitrary}. \tag{1.11.7}$$

**Remarks.** The definition (1.11.1) of the Laguerre polynomials can also be written as :

$$L_n^{(\alpha)}(x) = \frac{1}{n!}\sum_{k=0}^{n} \frac{(-n)_k}{k!} (\alpha+k+1)_{n-k} x^k.$$

In this way the Laguerre polynomials can be defined for all $\alpha$. Then we have the following connection with the Charlier polynomials defined by (1.12.1) :

$$\frac{(-a)^n}{n!} C_n(x;a) = L_n^{(x-n)}(a).$$

The Laguerre polynomials defined by (1.11.1) and the Hermite polynomials defined by (1.13.1) are connected by the following quadratic transformations :

$$H_{2n}(x) = (-1)^n n! 2^{2n} L_n^{(-\frac{1}{2})}(x^2)$$

and

$$H_{2n+1}(x) = (-1)^n n! 2^{2n+1} x L_n^{(\frac{1}{2})}(x^2).$$

In combinatorics the Laguerre polynomials with $\alpha = 0$ are often called Rook polynomials.

**References.** [1], [2], [3], [7], [9], [10], [14], [15], [25], [27], [31], [35], [36], [40], [45], [52], [55], [57], [58], [61], [65], [66], [67], [68], [70], [75], [77], [81], [82], [94], [95], [106], [116], [117], [119], [123], [124], [129], [134], [135], [139], [141], [143], [150], [152], [155], [156], [174], [180], [182], [205], [210], [211], [220], [222].

## 1.12 Charlier

**Definition.**

$$C_n(x;a) = {}_2F_0\left(\begin{array}{c}-n,-x\\ -\end{array}\bigg|-\frac{1}{a}\right). \tag{1.12.1}$$

**Orthogonality.**

$$\sum_{x=0}^{\infty} \frac{a^x}{x!} C_m(x;a) C_n(x;a) = n! a^{-n} e^a \delta_{mn}, \ a > 0. \tag{1.12.2}$$

**Recurrence relation.**

$$-xC_n(x;a) = aC_{n+1}(x;a) - (n+a)C_n(x;a) + nC_{n-1}(x;a). \tag{1.12.3}$$

**Difference equation.**

$$-ny(x) = ay(x+1) - (x+a)y(x) + xy(x-1), \ y(x) = C_n(x;a). \tag{1.12.4}$$

**Generating function.**

$$e^t \left(1 - \frac{t}{a}\right)^x = \sum_{n=0}^{\infty} \frac{C_n(x;a)}{n!} t^n. \tag{1.12.5}$$

**Remark.** The definition (1.11.1) of the Laguerre polynomials can also be written as :

$$L_n^{(\alpha)}(x) = \frac{1}{n!} \sum_{k=0}^{n} \frac{(-n)_k}{k!} (\alpha + k + 1)_{n-k} x^k.$$

In this way the Laguerre polynomials can be defined for all $\alpha$. Then we have the following connection with the Charlier polynomials defined by (1.12.1) :

$$\frac{(-a)^n}{n!} C_n(x;a) = L_n^{(x-n)}(a).$$

**References.** [7], [10], [15], [17], [25], [27], [45], [77], [78], [87], [94], [104], [105], [129], [156], [167], [174], [180], [220], [222], [233].

## 1.13 Hermite

**Definition.**

$$H_n(x) = (2x)^n {}_2F_0\left(\begin{array}{c} -n/2, -(n-1)/2 \\ - \end{array} \middle| -\frac{1}{x^2}\right). \tag{1.13.1}$$

**Orthogonality.**

$$\int_{-\infty}^{\infty} e^{-x^2} H_m(x) H_n(x) dx = 2^n n! \sqrt{\pi} \delta_{mn}. \tag{1.13.2}$$

**Recurrence relation.**

$$H_{n+1}(x) - 2xH_n(x) + 2nH_{n-1}(x) = 0. \tag{1.13.3}$$

**Differential equation.**

$$y''(x) - 2xy'(x) + 2ny(x) = 0, \ y(x) = H_n(x). \tag{1.13.4}$$

**Generating functions.**

$$\exp\left(2xt - t^2\right) = \sum_{n=0}^{\infty} \frac{H_n(x)}{n!} t^n. \tag{1.13.5}$$

$$\begin{cases} e^t \cos(2x\sqrt{t}) = \sum_{n=0}^{\infty} \frac{(-1)^n}{(2n)!} H_{2n}(x) t^n \\ \frac{e^t}{\sqrt{t}} \sin(2x\sqrt{t}) = \sum_{n=0}^{\infty} \frac{(-1)^n}{(2n+1)!} H_{2n+1}(x) t^n. \end{cases} \tag{1.13.6}$$

$$\begin{cases} e^{-t^2}\cosh(2xt) = \sum_{n=0}^{\infty} \frac{H_{2n}(x)}{(2n)!} t^{2n} \\ \\ e^{-t^2}\sinh(2xt) = \sum_{n=0}^{\infty} \frac{H_{2n+1}(x)}{(2n+1)!} t^{2n+1}. \end{cases} \quad (1.13.7)$$

$$\begin{cases} (1+t^2)^{-\gamma} {}_1F_1\left(\begin{array}{c}\gamma\\ \frac{1}{2}\end{array}\bigg|\frac{x^2t^2}{1+t^2}\right) = \sum_{n=0}^{\infty} \frac{(\gamma)_n}{(2n)!} H_{2n}(x) t^{2n} \\ \\ \frac{xt}{\sqrt{1+t^2}} {}_1F_1\left(\begin{array}{c}\gamma+\frac{1}{2}\\ \frac{3}{2}\end{array}\bigg|\frac{x^2t^2}{1+t^2}\right) = \sum_{n=0}^{\infty} \frac{(\gamma+\frac{1}{2})_n}{(2n+1)!} H_{2n+1}(x) t^{2n+1}. \end{cases} \quad (1.13.8)$$

$$\frac{1+2xt+4t^2}{(1+4t^2)^{\frac{3}{2}}} \exp\left(\frac{4x^2t^2}{1+4t^2}\right) = \sum_{n=0}^{\infty} \frac{H_n(x)}{[n/2]!} t^n, \quad (1.13.9)$$

where $[\alpha]$ denotes the largest integer smaller than or equal to $\alpha$.

**Remarks.** The Hermite polynomials can also be written as :

$$\frac{H_n(x)}{n!} = \sum_{k=0}^{[n/2]} \frac{(-1)^k (2x)^{n-2k}}{k!(n-2k)!},$$

where $[\alpha]$ denotes the largest integer smaller than or equal to $\alpha$.

The Laguerre polynomials defined by (1.11.1) and the Hermite polynomials defined by (1.13.1) are connected by the following quadratic transformations :

$$H_{2n}(x) = (-1)^n n! 2^{2n} L_n^{(-\frac{1}{2})}(x^2)$$

and

$$H_{2n+1}(x) = (-1)^n n! 2^{2n+1} x L_n^{(\frac{1}{2})}(x^2).$$

**References.** [2], [7], [10], [14], [15], [25], [27], [31], [35], [45], [49], [55], [57], [58], [71], [77], [81], [83], [94], [118], [119], [123], [156], [174], [180], [182], [205], [210], [220], [222], [225], [231].

# Chapter 2

# Limit relations between hypergeometric orthogonal polynomials

## 2.1 Wilson $\to$ Continuous dual Hahn

The continuous dual Hahn polynomials can be found from the Wilson polynomials defined by (1.1.1) by dividing by $(a+d)_n$ and letting $d \to \infty$ :

$$\lim_{d \to \infty} \frac{W_n(x^2;a,b,c,d)}{(a+d)_n} = S_n(x^2;a,b,c),$$

where $S_n(x^2;a,b,c)$ is defined by (1.3.1).

## 2.2 Wilson $\to$ Continuous Hahn

The continuous Hahn polynomials defined by (1.4.1) are obtained from the Wilson polynomials by the substitution $a \to a - it$, $b \to b - it$, $c \to c + it$, $d \to d + it$ and $x \to x + t$ in the definition (1.1.1) of the Wilson polynomials and the limit $t \to \infty$ in the following way :

$$\lim_{t \to \infty} \frac{W_n\left((x+t)^2;a-it,b-it,c+it,d+it\right)}{(-2t)^n n!} = p_n(x;a,b,c,d).$$

## 2.3 Wilson $\to$ Jacobi

The Jacobi polynomials given by (1.8.1) can be found from the Wilson polynomials by substituting $a = b = \frac{1}{2}(\alpha+1)$, $c = \frac{1}{2}(\beta+1) + it$, $d = \frac{1}{2}(\beta+1) - it$ and $x \to t\sqrt{\frac{1}{2}(1-x)}$ in the definition (1.1.1) of the Wilson polynomials and taking the limit $t \to \infty$. In fact we have

$$\lim_{t \to \infty} \frac{W_n\left(\frac{1}{2}(1-x)t^2;\frac{1}{2}(\alpha+1),\frac{1}{2}(\alpha+1),\frac{1}{2}(\beta+1)+it,\frac{1}{2}(\beta+1)-it\right)}{t^{2n} n!} = P_n^{(\alpha,\beta)}(x).$$

## 2.4 Racah $\to$ Hahn

If we take $\gamma + 1 = -N$ and let $\delta \to \infty$ in the definition (1.2.1) of the Racah polynomials, we obtain the Hahn polynomials defined by (1.5.1). Hence

$$\lim_{\delta \to \infty} R_n(\lambda(x);\alpha,\beta,-N-1,\delta) = Q_n(x;\alpha,\beta,N).$$

The Hahn polynomials can also be obtained from the Racah polynomials by taking $\delta = -\beta - N - 1$ in the definition (1.2.1) and letting $\gamma \to \infty$ :

$$\lim_{\gamma \to \infty} R_n(\lambda(x); \alpha, \beta, \gamma, -\beta - N - 1) = Q_n(x; \alpha, \beta, N).$$

Another way to do this is to take $\alpha + 1 = -N$ and $\beta \to \beta + \gamma + N + 1$ in the definition (1.2.1) of the Racah polynomials and then take the limit $\delta \to \infty$. In that case we obtain the Hahn polynomials given by (1.5.1) in the following way :

$$\lim_{\delta \to \infty} R_n(\lambda(x); -N - 1, \beta + \gamma + N + 1, \gamma, \delta) = Q_n(x; \gamma, \beta, N).$$

## 2.5 Racah → Dual Hahn

If we take $\alpha + 1 = -N$ and let $\beta \to \infty$ in (1.2.1), then we obtain the dual Hahn polynomials from the Racah polynomials. So we have

$$\lim_{\beta \to \infty} R_n(\lambda(x); -N - 1, \beta, \gamma, \delta) = R_n(\lambda(x); \gamma, \delta, N).$$

And if we take $\beta = -\delta - N - 1$ and let $\alpha \to \infty$ in (1.2.1), then we also obtain the dual Hahn polynomials :

$$\lim_{\alpha \to \infty} R_n(\lambda(x); \alpha, -\delta - N - 1, \gamma, \delta) = R_n(\lambda(x); \gamma, \delta, N).$$

Finally, if we take $\gamma + 1 = -N$ and $\delta \to \alpha + \delta + N + 1$ in the definition (1.2.1) of the Racah polynomials and take the limit $\beta \to \infty$ we find the dual Hahn polynomials given by (1.6.1) in the following way :

$$\lim_{\beta \to \infty} R_n(\lambda(x); \alpha, \beta, -N - 1, \alpha + \delta + N + 1) = R_n(\lambda(x); \alpha, \delta, N).$$

## 2.6 Continuous dual Hahn → Meixner-Pollaczek

The Meixner-Pollaczek polynomials given by (1.7.1) can be obtained from the continuous dual Hahn polynomials by the substitutions $x \to x - t$, $a = \lambda + it$, $b = \lambda - it$ and $c = t \cot \phi$ in the definition (1.3.1) and the limit $t \to \infty$ :

$$\lim_{t \to \infty} \frac{S_n\left((x-t)^2; \lambda + it, \lambda - it, t \cot \phi\right)}{\left(\frac{t}{\sin \phi}\right)_n n!} = P_n^{(\lambda)}(x; \phi).$$

## 2.7 Continuous Hahn → Meixner-Pollaczek

By taking $x \to x - t$, $a = \lambda + it$, $c = \lambda - it$ and $b = d = -t \tan \phi$ in the definition (1.4.1) of the continuous Hahn polynomials and taking the limit $t \to \infty$ we obtain the Meixner-Pollaczek polynomials defined by (1.7.1) :

$$\lim_{t \to \infty} \frac{p_n\left(x - t; \lambda + it, -t \tan \phi, \lambda - it, -t \tan \phi\right)}{\left(\frac{it}{\cos \phi}\right)_n i^n} = P_n^{(\lambda)}(x; \phi).$$

## 2.8 Continuous Hahn → Jacobi

The Jacobi polynomials defined by (1.8.1) follow from the continuous Hahn polynomials by the substitution $x \to -\frac{1}{2}xt$, $a = \frac{1}{2}(\alpha + 1 + it)$, $b = \frac{1}{2}(\beta + 1 - it)$, $c = \frac{1}{2}(\alpha + 1 - it)$ and $d = \frac{1}{2}(\beta + 1 + it)$ in (1.4.1), division by $(-1)^n t^n$ and the limit $t \to \infty$ :

$$\lim_{t \to \infty} \frac{p_n\left(-\frac{1}{2}xt; \frac{1}{2}(\alpha + 1 + it), \frac{1}{2}(\beta + 1 - it), \frac{1}{2}(\alpha + 1 - it), \frac{1}{2}(\beta + 1 + it)\right)}{(-1)^n t^n} = P_n^{(\alpha, \beta)}(x).$$

## 2.9 Hahn → Jacobi

To find the Jacobi polynomials from the Hahn polynomials we take $x \to Nx$ in (1.5.1) and let $N \to \infty$. We have

$$\lim_{N \to \infty} Q_n(Nx; \alpha, \beta, N) = \frac{P_n^{(\alpha,\beta)}(1-2x)}{P_n^{(\alpha,\beta)}(1)}.$$

## 2.10 Hahn → Meixner

If we take $\alpha = b - 1$, $\beta = N(1-c)c^{-1}$ in the definition (1.5.1) of the Hahn polynomials and let $N \to \infty$ we find the Meixner polynomials given by (1.9.1) :

$$\lim_{N \to \infty} Q_n\left(x; b-1, N\frac{1-c}{c}, N\right) = M_n(x; b, c).$$

## 2.11 Hahn → Krawtchouk

If we take $\alpha = pt$ and $\beta = (1-p)t$ in the definition (1.5.1) of the Hahn polynomials and let $t \to \infty$ we obtain the Krawtchouk polynomials defined by (1.10.1) :

$$\lim_{t \to \infty} Q_n\left(x; pt, (1-p)t, N\right) = K_n(x; p, N).$$

## 2.12 Dual Hahn → Meixner

To obtain the Meixner polynomials from the dual Hahn polynomials we have to take $\gamma = \beta - 1$ and $\delta = N(1-c)c^{-1}$ in the definition (1.6.1) of the dual Hahn polynomials and let $N \to \infty$ :

$$\lim_{N \to \infty} R_n\left(\lambda(x); \beta-1, N\frac{1-c}{c}, N\right) = M_n(x; \beta, c).$$

## 2.13 Dual Hahn → Krawtchouk

In the same way we find the Krawtchouk polynomials from the dual Hahn polynomials by setting $\gamma = pt$, $\delta = (1-p)t$ in (1.6.1) and let $t \to \infty$ :

$$\lim_{t \to \infty} R_n\left(\lambda(x); pt, (1-p)t, N\right) = K_n(x; p, N).$$

## 2.14 Meixner-Pollaczek → Laguerre

The Laguerre polynomials can be obtained from the Meixner-Pollaczek polynomials defined by (1.7.1) by the substitution $\lambda = \frac{1}{2}(\alpha + 1)$, $x \to -\frac{1}{2}\phi^{-1}x$ and letting $\phi \to 0$ :

$$\lim_{\phi \to 0} P_n^{(\frac{1}{2}\alpha + \frac{1}{2})}\left(-\frac{x}{2\phi}; \phi\right) = L_n^{(\alpha)}(x).$$

## 2.15 Meixner-Pollaczek → Hermite

If we substitute $x \to (\sin\phi)^{-1}(x\sqrt{\lambda} - \lambda\cos\phi)$ in the definition (1.7.1) of the Meixner-Pollaczek polynomials and then let $\lambda \to \infty$ we obtain the Hermite polynomials :

$$\lim_{\lambda \to \infty} \lambda^{-\frac{n}{2}} P_n^{(\lambda)}\left(\frac{x\sqrt{\lambda} - \lambda\cos\phi}{\sin\phi}; \phi\right) = \frac{H_n(x)}{n!}.$$

## 2.16 Jacobi → Laguerre

The Laguerre polynomials can be obtained from the Jacobi polynomials defined by (1.8.1) by letting $x \to 1 - 2\beta^{-1}x$ and then $\beta \to \infty$ :

$$\lim_{\beta \to \infty} P_n^{(\alpha,\beta)}\left(1 - \frac{2x}{\beta}\right) = L_n^{(\alpha)}(x).$$

## 2.17 Jacobi → Hermite

The Hermite polynomials given by (1.13.1) follow from the Jacobi polynomials defined by (1.8.1) by taking $\beta = \alpha$ and letting $\alpha \to \infty$ in the following way :

$$\lim_{\alpha \to \infty} \alpha^{-\frac{n}{2}} P_n^{(\alpha,\alpha)}\left(\frac{x}{\alpha^{\frac{1}{2}}}\right) = \frac{H_n(x)}{2^n n!}.$$

## 2.18 Meixner → Laguerre

If we take $\beta = \alpha + 1$ and $x \to (1-c)^{-1}x$ in the definition (1.9.1) of the Meixner polynomials and let $c \to 1$ we obtain the Laguerre polynomials :

$$\lim_{c \to 1} M_n\left(\frac{x}{1-c}; \alpha + 1, c\right) = \frac{L_n^{(\alpha)}(x)}{L_n^{(\alpha)}(0)}.$$

## 2.19 Meixner → Charlier

If we take $c = (a+\beta)^{-1}a$ in the definition (1.9.1) of the Meixner polynomials and let $\beta \to \infty$ we find the Charlier polynomials :

$$\lim_{\beta \to \infty} M_n\left(x; \beta, \frac{a}{a+\beta}\right) = C_n(x;a).$$

## 2.20 Krawtchouk → Charlier

The Charlier polynomials given by (1.12.1) can be found from the Krawtchouk polynomials defined by (1.10.1) by taking $p = N^{-1}a$ and let $N \to \infty$ :

$$\lim_{N \to \infty} K_n\left(x; \frac{a}{N}, N\right) = C_n(x;a).$$

## 2.21 Krawtchouk → Hermite

The Hermite polynomials follow from the Krawtchouk polynomials defined by (1.10.1) by setting $x \to pN + x\sqrt{2p(1-p)N}$ and then letting $N \to \infty$ :

$$\lim_{N \to \infty} \sqrt{\binom{N}{n}} K_n\left(pN + x\sqrt{2p(1-p)N}; p, N\right) = \frac{(-1)^n H_n(x)}{\sqrt{2^n (n!) \left(\frac{p}{1-p}\right)^n}}.$$

## 2.22 Laguerre → Hermite

The Hermite polynomials defined by (1.13.1) can be obtained from the Laguerre polynomials given by (1.11.1) by taking the limit $\alpha \to \infty$ in the following way :

$$\lim_{\alpha \to \infty} \left(\frac{2}{\alpha}\right)^{\frac{n}{2}} L_n^{(\alpha)}\left((2\alpha)^{\frac{1}{2}}x + \alpha\right) = \frac{(-1)^n}{n!} H_n(x).$$

## 2.23 Charlier → Hermite

If we set $x \to (2a)^{1/2}x + a$ in the definition (1.12.1) of the Charlier polynomials and let $a \to \infty$ we find the Hermite polynomials defined by (1.13.1). In fact we have

$$\lim_{a \to \infty} (2a)^{\frac{n}{2}} C_n\left((2a)^{\frac{1}{2}}x + a; a\right) = (-1)^n H_n(x).$$

# SCHEME
## OF
# BASIC HYPERGEOMETRIC
# ORTHOGONAL POLYNOMIALS

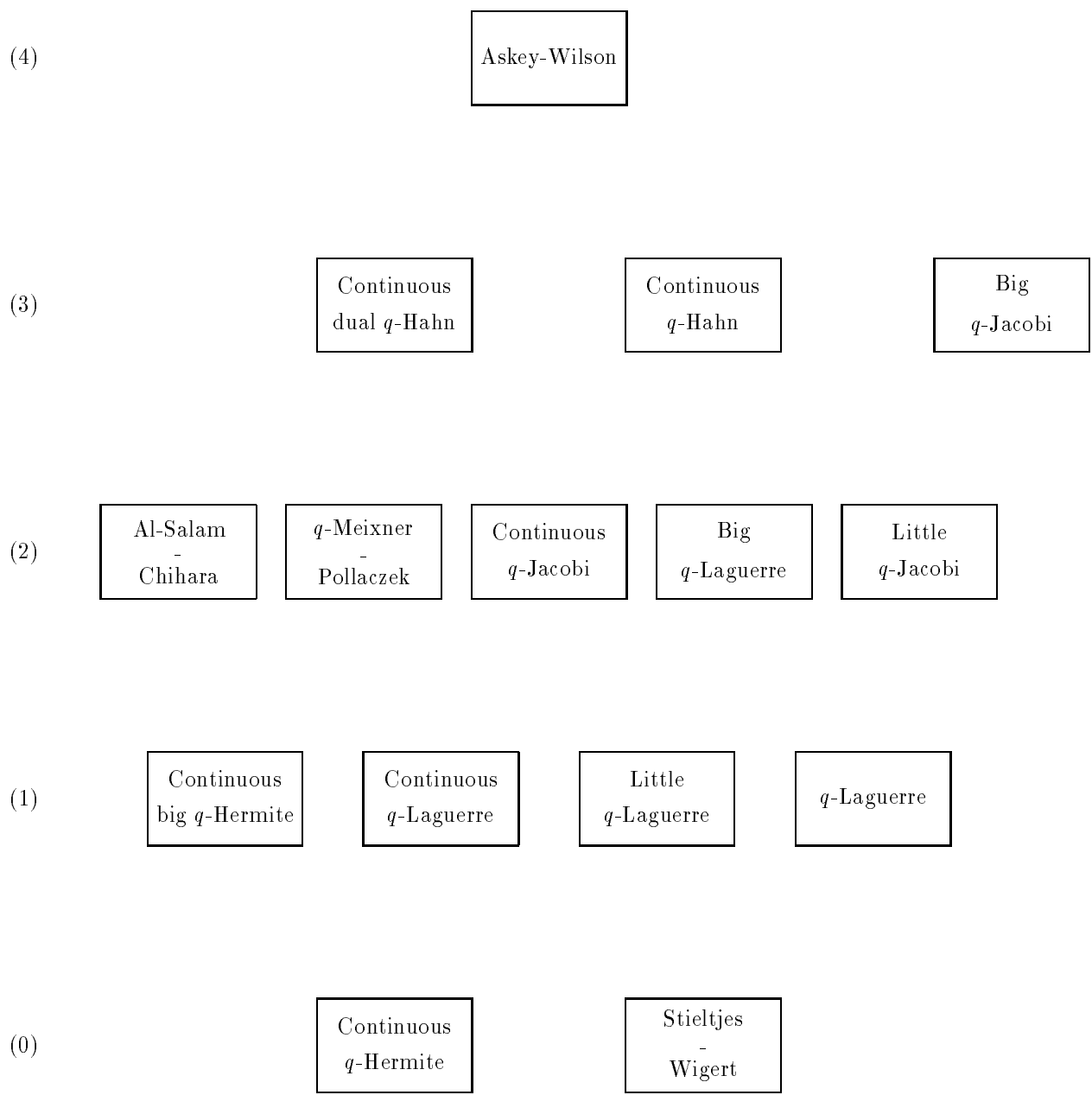

# SCHEME
OF
# BASIC HYPERGEOMETRIC
# ORTHOGONAL POLYNOMIALS

| | | | | |
|---|---|---|---|---|
| | | $q$-Racah | | (4) |

| Big $q$-Jacobi | | $q$-Hahn | | Dual $q$-Hahn | | (3) |

| $q$-Meixner | Quantum $q$-Krawtchouk | $q$-Krawtchouk | Affine $q$-Krawtchouk | Dual $q$-Krawtchouk | (2) |

| Alternative $q$-Charlier | $q$-Charlier | Al-Salam - Carlitz I | Al-Salam - Carlitz II | (1) |

| Discrete $q$-Hermite I | Discrete $q$-Hermite II | (0) |

# Chapter 3

# Basic hypergeometric orthogonal polynomials

## 3.1 Askey-Wilson

**Definition.**
$$\frac{a^n p_n(x;a,b,c,d|q)}{(ab,ac,ad;q)_n} = {}_4\phi_3\left(\begin{array}{c} q^{-n}, abcdq^{n-1}, ae^{i\theta}, ae^{-i\theta} \\ ab, ac, ad \end{array} \bigg| q; q\right), \quad x = \cos\theta. \qquad (3.1.1)$$

The Askey-Wilson polynomials are $q$-analogues of the Wilson polynomials given by (1.1.1).

**Orthogonality.** When $a, b, c, d$ are real, or occur in complex conjugate pairs if complex, and $\max(|a|, |b|, |c|, |d|) < 1$, then we have the following orthogonality relation

$$\frac{1}{2\pi} \int_{-1}^{1} \frac{w(x)}{\sqrt{1-x^2}} p_m(x;a,b,c,d|q) p_n(x;a,b,c,d|q) dx = h_n \delta_{mn}, \qquad (3.1.2)$$

where

$$w(x) := w(x;a,b,c,d|q) = \left|\frac{(e^{2i\theta};q)_\infty}{(ae^{i\theta}, be^{i\theta}, ce^{i\theta}, de^{i\theta};q)_\infty}\right|^2 = \frac{h(x,1)h(x,-1)h(x,q^{\frac{1}{2}})h(x,-q^{\frac{1}{2}})}{h(x,a)h(x,b)h(x,c)h(x,d)},$$

with

$$h(x,\alpha) := \prod_{k=0}^{\infty} \left[1 - 2\alpha x q^k + \alpha^2 q^{2k}\right] = \left(\alpha e^{i\theta}, \alpha e^{-i\theta}; q\right)_\infty, \quad x = \cos\theta$$

and

$$h_n = \frac{(abcdq^{n-1};q)_n (abcdq^{2n};q)_\infty}{(q^{n+1}, abq^n, acq^n, adq^n, bcq^n, bdq^n, cdq^n;q)_\infty}.$$

If $a > 1$ and $b, c, d$ are real or one is real and the other two are complex conjugates, $\max(|b|, |c|, |d|) < 1$ and the pairwise products of $a, b, c$ and $d$ have absolute value less than one, then we have another orthogonality relation given by :

$$\frac{1}{2\pi} \int_{-1}^{1} \frac{w(x)}{\sqrt{1-x^2}} p_m(x;a,b,c,d|q) p_n(x;a,b,c,d|q) dx +$$
$$+ \sum_{\substack{k \\ 1 < aq^k \leq a}} w_k p_m(x_k;a,b,c,d|q) p_n(x_k;a,b,c,d|q) = h_n \delta_{mn}, \qquad (3.1.3)$$

where $w(x)$ and $h_n$ are as before,

$$x_k = \frac{aq^k + (aq^k)^{-1}}{2}$$

and

$$w_k = \frac{(a^{-2};q)_\infty}{(q,ab,ac,ad,a^{-1}b,a^{-1}c,a^{-1}d;q)_\infty} \frac{(1-a^2q^{2k})(a^2,ab,ac,ad;q)_k}{(1-a^2)(q,ab^{-1}q,ac^{-1}q,ad^{-1}q;q)_k} \left(\frac{q}{abcd}\right)^k.$$

**Recurrence relation.**

$$2x\tilde{p}_n(x) = A_n\tilde{p}_{n+1}(x) + \left[a + a^{-1} - (A_n + C_n)\right]\tilde{p}_n(x) + C_n\tilde{p}_{n-1}(x), \qquad (3.1.4)$$

where

$$\tilde{p}_n(x) := \frac{a^n p_n(x;a,b,c,d|q)}{(ab,ac,ad;q)_n}$$

and

$$\begin{cases} A_n = \dfrac{(1-abq^n)(1-acq^n)(1-adq^n)(1-abcdq^{n-1})}{a(1-abcdq^{2n-1})(1-abcdq^{2n})} \\[2ex] C_n = \dfrac{a(1-q^n)(1-bcq^{n-1})(1-bdq^{n-1})(1-cdq^{n-1})}{(1-abcdq^{2n-2})(1-abcdq^{2n-1})}. \end{cases}$$

**$q$-Difference equation.**

$$(1-q)^2 D_q\left[\tilde{w}(x;aq^{\frac{1}{2}},bq^{\frac{1}{2}},cq^{\frac{1}{2}},dq^{\frac{1}{2}}|q)D_q y(x)\right] + $$
$$+ \lambda_n \tilde{w}(x;a,b,c,d|q)y(x) = 0, \ y(x) = p_n(x;a,b,c,d|q), \qquad (3.1.5)$$

where

$$\tilde{w}(x;a,b,c,d|q) := \frac{w(x;a,b,c,d|q)}{\sqrt{1-x^2}},$$

$$\lambda_n = 4q^{-n+1}(1-q^n)(1-abcdq^{n-1})$$

and

$$D_q f(x) := \frac{\delta_q f(x)}{\delta_q x} \text{ with } \delta_q f(e^{i\theta}) = f(q^{\frac{1}{2}}e^{i\theta}) - f(q^{-\frac{1}{2}}e^{i\theta}), \ x = \cos\theta.$$

If we define

$$P_n(z) := \frac{(ab,ac,ad;q)_n}{a^n} {}_4\phi_3\left(\begin{array}{c} q^{-n},abcdq^{n-1},az,az^{-1} \\ ab,ac,ad \end{array} \bigg| q;q\right)$$

then the $q$-difference equation can also be written in the form

$$q^{-n}(1-q^n)(1-abcdq^{n-1})P_n(z)$$
$$= A(z)P_n(qz) - \left[A(z) + A(z^{-1})\right]P_n(z) + A(z^{-1})P_n(q^{-1}z), \qquad (3.1.6)$$

where

$$A(z) = \frac{(1-az)(1-bz)(1-cz)(1-dz)}{(1-z^2)(1-qz^2)}.$$

**Generating functions.**

$${}_2\phi_1\left(\begin{array}{c} ae^{i\theta},be^{i\theta} \\ ab \end{array} \bigg| q;e^{-i\theta}t\right) {}_2\phi_1\left(\begin{array}{c} ce^{-i\theta},de^{-i\theta} \\ cd \end{array} \bigg| q;e^{i\theta}t\right) = \sum_{n=0}^{\infty} \frac{p_n(x;a,b,c,d|q)}{(ab,cd,q;q)_n}t^n, \ x=\cos\theta. \quad (3.1.7)$$

$${}_2\phi_1\left(\begin{array}{c} ae^{i\theta},ce^{i\theta} \\ ac \end{array} \bigg| q;e^{-i\theta}t\right) {}_2\phi_1\left(\begin{array}{c} be^{-i\theta},de^{-i\theta} \\ bd \end{array} \bigg| q;e^{i\theta}t\right) = \sum_{n=0}^{\infty} \frac{p_n(x;a,b,c,d|q)}{(ac,bd,q;q)_n}t^n, \ x=\cos\theta. \quad (3.1.8)$$

$$_2\phi_1\left(\begin{matrix}ae^{i\theta},de^{i\theta}\\ad\end{matrix}\bigg|q;e^{-i\theta}t\right){}_2\phi_1\left(\begin{matrix}be^{-i\theta},ce^{-i\theta}\\bc\end{matrix}\bigg|q;e^{i\theta}t\right)=\sum_{n=0}^{\infty}\frac{p_n(x;a,b,c,d|q)}{(ad,bc,q;q)_n}t^n,\ x=\cos\theta. \quad(3.1.9)$$

**Remark.** The $q$-Racah polynomials defined by (3.2.1) and the Askey-Wilson polynomials given by (3.1.1) are related in the following way. If we substitute $a^2=\gamma\delta q$, $b^2=\alpha^2\gamma^{-1}\delta^{-1}q$, $c^2=\beta^2\gamma^{-1}\delta q$, $d^2=\gamma\delta^{-1}q$ and $e^{2i\theta}=\gamma\delta q^{2x+1}$ in the definition (3.1.1) of the Askey-Wilson polynomials we find :

$$R_n(\mu(x);\alpha,\beta,\gamma,\delta|q)=\frac{(\gamma\delta q)^{\frac{1}{2}n}p_n(\nu(x);\gamma^{\frac{1}{2}}\delta^{\frac{1}{2}}q^{\frac{1}{2}},\alpha\gamma^{-\frac{1}{2}}\delta^{-\frac{1}{2}}q^{\frac{1}{2}},\beta\gamma^{-\frac{1}{2}}\delta^{\frac{1}{2}}q^{\frac{1}{2}},\gamma^{\frac{1}{2}}\delta^{-\frac{1}{2}}q^{\frac{1}{2}}|q)}{(\alpha q,\beta\delta q,\gamma q;q)_n},$$

where

$$\nu(x)=\frac{1}{2}\gamma^{\frac{1}{2}}\delta^{\frac{1}{2}}q^{x+\frac{1}{2}}+\frac{1}{2}\gamma^{-\frac{1}{2}}\delta^{-\frac{1}{2}}q^{-x-\frac{1}{2}}.$$

**References.** [10], [25], [31], [45], [47], [48], [73], [112], [114], [127], [131], [133], [134], [136], [137], [140], [162], [163], [166], [176], [179], [180], [197], [198], [200], [202], [228].

## 3.2 $q$-Racah

**Definition.**
$$R_n(\mu(x);\alpha,\beta,\gamma,\delta|q)={}_4\tilde{\phi}_3\left(\begin{matrix}q^{-n},\alpha\beta q^{n+1},q^{-x},\gamma\delta q^{x+1}\\\alpha q,\beta\delta q,\gamma q\end{matrix}\bigg|q;q\right),\ n=0,1,2,\ldots,N,\quad(3.2.1)$$

where
$$\mu(x):=q^{-x}+\gamma\delta q^{x+1}$$

and
$$\alpha q=q^{-N}\text{ or }\beta\delta q=q^{-N}\text{ or }\gamma q=q^{-N},\text{ with }N\text{ a nonnegative integer.}$$

Since
$$(q^{-x},\gamma\delta q^{x+1};q)_k=\prod_{j=0}^{k-1}\left(1-\mu(x)q^j+\gamma\delta q^{2j+1}\right),$$

it is clear that $R_n(\mu(x);\alpha,\beta,\gamma,\delta|q)$ is a polynomial of degree $n$ in $\mu(x)$.

**Orthogonality.**
$$\sum_{x=0}^{N}\frac{(\gamma\delta q,\alpha q,\beta\delta q,\gamma q;q)_x}{(q,\alpha^{-1}\gamma\delta q,\beta^{-1}\gamma q,\delta q;q)_x}\frac{(1-\gamma\delta q^{2x+1})}{(\alpha\beta q)^x(1-\gamma\delta q)}R_m(\mu(x))R_n(\mu(x))=h_n\delta_{mn},\quad(3.2.2)$$

where
$$R_n(\mu(x)):=R_n(\mu(x);\alpha,\beta,\gamma,\delta|q)$$

and
$$h_n=\frac{(\gamma\delta q^2,\alpha^{-1}\beta^{-1}\gamma,\alpha^{-1}\delta,\beta^{-1};q)_\infty}{(\alpha^{-1}\gamma\delta q,\beta^{-1}\gamma q,\delta q,\alpha^{-1}\beta^{-1}q^{-1};q)_\infty}\frac{(1-\alpha\beta q)(\gamma\delta q)^n}{(1-\alpha\beta q^{2n+1})}\frac{(q,\beta q,\alpha\delta^{-1}q,\alpha\beta\gamma^{-1}q;q)_n}{(\alpha\beta q,\alpha q,\beta\delta q,\gamma q;q)_n}.$$

This implies

$$h_n=\begin{cases}\dfrac{(\gamma\delta q^2,\beta^{-1};q)_N}{(\beta^{-1}\gamma q,\delta q;q)_N}\dfrac{(1-\beta q^{-N})(\gamma\delta q)^n}{(1-\beta q^{2n-N})}\dfrac{(q,\beta q,\delta^{-1}q^{-N},\beta\gamma^{-1}q^{-N};q)_n}{(\beta q^{-N},q^{-N},\beta\delta q,\gamma q;q)_n}&\text{if }\alpha q=q^{-N}\\[2ex]\dfrac{(\beta\gamma^{-1},\alpha\beta q^2;q)_N}{(\alpha\beta\gamma^{-1}q,\beta q;q)_N}\dfrac{(1-\alpha\beta q)(\beta^{-1}\gamma q^{-N})^n}{(1-\alpha\beta q^{2n+1})}\dfrac{(q,\beta q,\alpha\beta q^{N+2},\alpha\beta\gamma^{-1}q;q)_n}{(\alpha\beta q,\alpha q,q^{-N},\gamma q;q)_n}&\text{if }\beta\delta q=q^{-N}\\[2ex]\dfrac{(\alpha\beta q^2,\delta^{-1};q)_N}{(\beta q,\alpha\delta^{-1}q;q)_N}\dfrac{(1-\alpha\beta q)(\delta q^{-N})^n}{(1-\alpha\beta q^{2n+1})}\dfrac{(q,\beta q,\alpha\delta^{-1}q,\alpha\beta q^{N+2};q)_n}{(\alpha\beta q,\alpha q,\beta\delta q,q^{-N};q)_n}&\text{if }\gamma q=q^{-N}.\end{cases}$$

**Recurrence relation.**

$$-\left(1 - q^{-x}\right)\left(1 - \gamma\delta q^{x+1}\right) R_n(\mu(x))$$
$$= A_n R_{n+1}(\mu(x)) - (A_n + C_n) R_n(\mu(x)) + C_n R_{n-1}(\mu(x)), \qquad (3.2.3)$$

where

$$\begin{cases} A_n = \dfrac{(1 - \alpha q^{n+1})(1 - \gamma q^{n+1})(1 - \alpha\beta q^{n+1})(1 - \beta\delta q^{n+1})}{(1 - \alpha\beta q^{2n+1})(1 - \alpha\beta q^{2n+2})} \\[2mm] C_n = \dfrac{q(1 - q^n)(1 - \beta q^n)(\delta - \alpha q^n)(\gamma - \alpha\beta q^n)}{(1 - \alpha\beta q^{2n})(1 - \alpha\beta q^{2n+1})}. \end{cases}$$

**$q$-Difference equation.**

$$\Delta\left[w(x-1)B(x-1)\Delta y(x-1)\right] + $$
$$- q^{-n}(1 - q^n)(1 - \alpha\beta q^{n+1})w(x)y(x) = 0, \; y(x) = R_n(\mu(x); \alpha, \beta, \gamma, \delta|q), \quad (3.2.4)$$

where

$$\Delta f(x) := f(x+1) - f(x),$$

$$w(x) = \frac{(\gamma\delta q, \alpha q, \beta\delta q, \gamma q; q)_x}{(q, \alpha^{-1}\gamma\delta q, \beta^{-1}\gamma q, \delta q; q)_x} \frac{(1 - \gamma\delta q^{2x+1})}{(\alpha\beta q)^x(1 - \gamma\delta q)}$$

and $B(x)$ as below. This $q$-difference equation can also be written in the form

$$q^{-n}(1 - q^n)(1 - \alpha\beta q^{n+1})y(x) = B(x)y(x+1) - [B(x) + D(x)]y(x) + D(x)y(x-1), \quad (3.2.5)$$

where

$$y(x) = R_n(\mu(x); \alpha, \beta, \gamma, \delta|q)$$

and

$$\begin{cases} B(x) = \dfrac{(1 - \alpha q^{x+1})(1 - \beta\delta q^{x+1})(1 - \gamma q^{x+1})(1 - \gamma\delta q^{x+1})}{(1 - \gamma\delta q^{2x+1})(1 - \gamma\delta q^{2x+2})} \\[2mm] D(x) = \dfrac{q(1 - q^x)(1 - \delta q^x)(\beta - \gamma q^x)(\alpha - \gamma\delta q^x)}{(1 - \gamma\delta q^{2x})(1 - \gamma\delta q^{2x+1})}. \end{cases}$$

**Generating functions.**

$$_2\tilde{\phi}_1\left(\begin{matrix} \alpha q^{x+1}, \gamma\delta q^{x+1} \\ \alpha q \end{matrix} \bigg| q; q^{-x}t\right) {}_2\phi_1\left(\begin{matrix} \beta\gamma^{-1}q^{-x}, \delta^{-1}q^{-x} \\ \beta q \end{matrix} \bigg| q; \gamma\delta q^{x+1}t\right)$$
$$\simeq \sum_{n=0}^{N} \frac{(\beta\delta q, \gamma q; q)_n}{(\beta q, q; q)_n} R_n(\mu(x); \alpha, \beta, \gamma, \delta|q)t^n. \qquad (3.2.6)$$

$$_2\tilde{\phi}_1\left(\begin{matrix} \beta\delta q^{x+1}, \gamma\delta q^{x+1} \\ \beta\delta q \end{matrix} \bigg| q; q^{-x}t\right) {}_2\phi_1\left(\begin{matrix} \alpha\gamma^{-1}\delta^{-1}q^{-x}, \delta^{-1}q^{-x} \\ \alpha\delta^{-1}q \end{matrix} \bigg| q; \gamma\delta q^{x+1}t\right)$$
$$\simeq \sum_{n=0}^{N} \frac{(\alpha q, \gamma q; q)_n}{(\alpha\delta^{-1}q, q; q)_n} R_n(\mu(x); \alpha, \beta, \gamma, \delta|q)t^n. \qquad (3.2.7)$$

$$_2\tilde{\phi}_1\left(\begin{matrix} \gamma q^{x+1}, \gamma\delta q^{x+1} \\ \gamma q \end{matrix} \bigg| q; q^{-x}t\right) {}_2\phi_1\left(\begin{matrix} \alpha\gamma^{-1}\delta^{-1}q^{-x}, \beta\gamma^{-1}q^{-x} \\ \alpha\beta\gamma^{-1}q \end{matrix} \bigg| q; \gamma\delta q^{x+1}t\right)$$
$$\simeq \sum_{n=0}^{N} \frac{(\alpha q, \beta\delta q; q)_n}{(\alpha\beta\gamma^{-1}q, q; q)_n} R_n(\mu(x); \alpha, \beta, \gamma, \delta|q)t^n. \qquad (3.2.8)$$

**Remark.** The Askey-Wilson polynomials defined by (3.1.1) and the $q$-Racah polynomials given by (3.2.1) are related in the following way. If we substitute $\alpha = abq^{-1}$, $\beta = cdq^{-1}$, $\gamma = adq^{-1}$, $\delta = ad^{-1}$ and $q^x = a^{-1}e^{-i\theta}$ in the definition (3.2.1) of the $q$-Racah polynomials we find :

$$\mu(x) = 2a\cos\theta$$

and

$$R_n\left(2a\cos\theta; abq^{-1}, cdq^{-1}, adq^{-1}, ad^{-1}|q\right) = \frac{a^n p_n(x;a,b,c,d|q)}{(ab,ac,ad;q)_n}.$$

**References.** [10], [22], [25], [43], [45], [111], [114], [127], [160], [180], [183], [197].

## 3.3 Continuous dual $q$-Hahn

**Definition.**

$$\frac{a^n p_n(x;a,b,c|q)}{(ab,ac;q)_n} = {}_3\phi_2\left(\begin{array}{c}q^{-n}, ae^{i\theta}, ae^{-i\theta}\\ab,ac\end{array}\bigg| q;q\right), \quad x = \cos\theta. \tag{3.3.1}$$

**Orthogonality.** When $a, b, c$ are real, or one is real and the other two are complex conjugates, and $\max(|a|,|b|,|c|) < 1$, then we have the following orthogonality relation

$$\frac{1}{2\pi}\int_{-1}^{1}\frac{w(x)}{\sqrt{1-x^2}}p_m(x;a,b,c|q)p_n(x;a,b,c|q)dx = h_n\delta_{mn}, \tag{3.3.2}$$

where

$$w(x) := w(x;a,b,c|q) = \left|\frac{(e^{2i\theta};q)_\infty}{(ae^{i\theta},be^{i\theta},ce^{i\theta};q)_\infty}\right|^2 = \frac{h(x,1)h(x,-1)h(x,q^{\frac{1}{2}})h(x,-q^{\frac{1}{2}})}{h(x,a)h(x,b)h(x,c)},$$

with

$$h(x,\alpha) := \prod_{k=0}^{\infty}\left[1 - 2\alpha x q^k + \alpha^2 q^{2k}\right] = \left(\alpha e^{i\theta}, \alpha e^{-i\theta}; q\right)_\infty, \quad x = \cos\theta$$

and

$$h_n = \frac{1}{(q^{n+1}, abq^n, acq^n, bcq^n; q)_\infty}.$$

If $a > 1$ and $b$ and $c$ are real or complex conjugates, $\max(|b|,|c|) < 1$ and the pairwise products of $a, b$ and $c$ have absolute value less than one, then we have another orthogonality relation given by :

$$\frac{1}{2\pi}\int_{-1}^{1}\frac{w(x)}{\sqrt{1-x^2}}p_m(x;a,b,c|q)p_n(x;a,b,c|q)dx +$$
$$+ \sum_{\substack{k \\ 1 < aq^k \leq a}} w_k p_m(x_k;a,b,c|q)p_n(x_k;a,b,c|q) = h_n\delta_{mn}, \tag{3.3.3}$$

where $w(x)$ and $h_n$ are as before,

$$x_k = \frac{aq^k + (aq^k)^{-1}}{2}$$

and

$$w_k = \frac{(a^{-2};q)_\infty}{(q,ab,ac,a^{-1}b,a^{-1}c;q)_\infty}\frac{(1-a^2q^{2k})(a^2,ab,ac;q)_k}{(1-a^2)(q,ab^{-1}q,ac^{-1}q;q)_k}(-1)^k q^{-\binom{k}{2}}\left(\frac{1}{a^2bc}\right)^k.$$

**Recurrence relation.**

$$2x\tilde{p}_n(x) = A_n\tilde{p}_{n+1}(x) + \left[a + a^{-1} - (A_n + C_n)\right]\tilde{p}_n(x) + C_n\tilde{p}_{n-1}(x), \qquad (3.3.4)$$

where

$$\tilde{p}_n(x) := \frac{a^n p_n(x;a,b,c|q)}{(ab,ac;q)_n}$$

and

$$\begin{cases} A_n = a^{-1}(1 - abq^n)(1 - acq^n) \\ C_n = a(1 - q^n)(1 - bcq^{n-1}). \end{cases}$$

**$q$-Difference equation.**

$$(1-q)^2 D_q\left[\tilde{w}(x;aq^{\frac{1}{2}},bq^{\frac{1}{2}},cq^{\frac{1}{2}}|q)D_q y(x)\right] + \\ + 4q^{-n+1}(1-q^n)\tilde{w}(x;a,b,c|q)y(x) = 0, \ y(x) = p_n(x;a,b,c|q), \quad (3.3.5)$$

where

$$\tilde{w}(x;a,b,c|q) := \frac{w(x;a,b,c|q)}{\sqrt{1-x^2}}$$

and

$$D_q f(x) := \frac{\delta_q f(x)}{\delta_q x} \text{ with } \delta_q f(e^{i\theta}) = f(q^{\frac{1}{2}}e^{i\theta}) - f(q^{-\frac{1}{2}}e^{i\theta}), \ x = \cos\theta.$$

If we define

$$P_n(z) := \frac{(ab,ac;q)_n}{a^n}{}_3\phi_2\left(\begin{array}{c}q^{-n}, az, az^{-1} \\ ab, ac\end{array}\bigg| q;q\right)$$

then the $q$-difference equation can also be written in the form

$$q^{-n}(1-q^n)P_n(z) = A(z)P_n(qz) - \left[A(z) + A(z^{-1})\right]P_n(z) + A(z^{-1})P_n(q^{-1}z), \qquad (3.3.6)$$

where

$$A(z) = \frac{(1-az)(1-bz)(1-cz)}{(1-z^2)(1-qz^2)}.$$

**Generating functions.**

$$\frac{(ct;q)_\infty}{(e^{i\theta}t;q)_\infty}{}_2\phi_1\left(\begin{array}{c}ae^{i\theta}, be^{i\theta} \\ ab\end{array}\bigg| q; e^{-i\theta}t\right) = \sum_{n=0}^{\infty}\frac{p_n(x;a,b,c|q)}{(ab,q;q)_n}t^n, \ x = \cos\theta. \qquad (3.3.7)$$

$$\frac{(bt;q)_\infty}{(e^{i\theta}t;q)_\infty}{}_2\phi_1\left(\begin{array}{c}ae^{i\theta}, ce^{i\theta} \\ ac\end{array}\bigg| q; e^{-i\theta}t\right) = \sum_{n=0}^{\infty}\frac{p_n(x;a,b,c|q)}{(ac,q;q)_n}t^n, \ x = \cos\theta. \qquad (3.3.8)$$

$$\frac{(at;q)_\infty}{(e^{i\theta}t;q)_\infty}{}_2\phi_1\left(\begin{array}{c}be^{i\theta}, ce^{i\theta} \\ bc\end{array}\bigg| q; e^{-i\theta}t\right) = \sum_{n=0}^{\infty}\frac{p_n(x;a,b,c|q)}{(bc,q;q)_n}t^n, \ x = \cos\theta. \qquad (3.3.9)$$

**References.**

## 3.4 Continuous $q$-Hahn

**Definition.**

$$\frac{(ae^{i\phi})^n p_n(x;a,b,c,d;q)}{(abe^{2i\phi},ac,ad;q)_n} = {}_4\phi_3\left(\begin{array}{c}q^{-n}, abcdq^{n-1}, ae^{i(\theta+2\phi)}, ae^{-i\theta} \\ abe^{2i\phi}, ac, ad\end{array}\bigg| q;q\right), \ x = \cos(\theta+\phi). \quad (3.4.1)$$

**Orthogonality.** When $c = a$ and $d = b$ then we have, if $a$ and $b$ are real and $\max(|a|, |b|) < 1$ or if $b = \overline{a}$ and $|a| < 1$:

$$\frac{1}{4\pi} \int_{-\pi}^{\pi} w(\cos(\theta + \phi)) p_m(\cos(\theta + \phi); a, b, c, d; q) p_n(\cos(\theta + \phi); a, b, c, d; q) d\theta = h_n \delta_{mn}, \quad (3.4.2)$$

where

$$w(x) := w(x; a, b, c, d; q) = \left| \frac{(e^{2i(\theta+\phi)}; q)_\infty}{(ae^{i(\theta+2\phi)}, be^{i(\theta+2\phi)}, ce^{i\theta}, de^{i\theta}; q)_\infty} \right|^2$$

$$= \frac{h(x, 1) h(x, -1) h(x, q^{\frac{1}{2}}) h(x, -q^{\frac{1}{2}})}{h(x, ae^{i\phi}) h(x, be^{i\phi}) h(x, ce^{-i\phi}) h(x, de^{-i\phi})},$$

with

$$h(x, \alpha) := \prod_{k=0}^{\infty} \left[ 1 - 2\alpha x q^k + \alpha^2 q^{2k} \right] = \left( \alpha e^{i(\theta+\phi)}, \alpha e^{-i(\theta+\phi)}; q \right)_\infty, \quad x = \cos(\theta + \phi)$$

and

$$h_n = \frac{(abcdq^{n-1}; q)_n (abcdq^{2n}; q)_\infty}{(q^{n+1}, abq^n e^{2i\phi}, acq^n, adq^n, bcq^n, bdq^n, cdq^n e^{-2i\phi}; q)_\infty}.$$

**Recurrence relation.**

$$2x \tilde{p}_n(x) = A_n \tilde{p}_{n+1}(x) + \left[ ae^{i\phi} + a^{-1} e^{-i\phi} - (A_n + C_n) \right] \tilde{p}_n(x) + C_n \tilde{p}_{n-1}(x), \quad (3.4.3)$$

where

$$\tilde{p}_n(x) := \frac{(ae^{i\phi})^n p_n(x; a, b, c, d; q)}{(ac, ad, abe^{2i\phi}; q)_n}$$

and

$$\begin{cases} A_n = \dfrac{(1 - abe^{2i\phi} q^n)(1 - acq^n)(1 - adq^n)(1 - abcdq^{n-1})}{ae^{i\phi}(1 - abcdq^{2n-1})(1 - abcdq^{2n})} \\ C_n = \dfrac{ae^{i\phi}(1 - q^n)(1 - bcq^{n-1})(1 - bdq^{n-1})(1 - cde^{-2i\phi} q^{n-1})}{(1 - abcdq^{2n-2})(1 - abcdq^{2n-1})}. \end{cases}$$

**$q$-Difference equation.**

$$(1-q)^2 D_q \left[ \tilde{w}(x; aq^{\frac{1}{2}}, bq^{\frac{1}{2}}, cq^{\frac{1}{2}}, dq^{\frac{1}{2}}; q) D_q y(x) \right] +$$
$$+ \lambda_n \tilde{w}(x; a, b, c, d; q) y(x) = 0, \ y(x) = p_n(x; a, b, c, d; q), \quad (3.4.4)$$

where

$$\tilde{w}(x; a, b, c, d; q) := \frac{w(x; a, b, c, d; q)}{\sqrt{1 - x^2}},$$

$$\lambda_n = 4q^{-n+1}(1 - q^n)(1 - abcdq^{n-1})$$

and

$$D_q f(x) := \frac{\delta_q f(x)}{\delta_q x} \ \text{with} \ \delta_q f(e^{i(\theta+\phi)}) = f(q^{\frac{1}{2}} e^{i(\theta+\phi)}) - f(q^{-\frac{1}{2}} e^{i(\theta+\phi)}), \ x = \cos(\theta + \phi).$$

**Generating functions.**

$$_2\phi_1 \left( \begin{array}{c} ae^{i(\theta+2\phi)}, be^{i(\theta+2\phi)} \\ abe^{2i\phi} \end{array} \bigg| q; e^{-i(\theta+\phi)} t \right) {}_2\phi_1 \left( \begin{array}{c} ce^{-i(\theta+2\phi)}, de^{-i(\theta+2\phi)} \\ cde^{-2i\phi} \end{array} \bigg| q; e^{i(\theta+\phi)} t \right)$$

$$= \sum_{n=0}^{\infty} \frac{p_n(x; a, b, c, d; q) t^n}{(abe^{2i\phi}, cde^{-2i\phi}, q; q)_n}, \ x = \cos(\theta + \phi). \quad (3.4.5)$$

$$_2\phi_1\left(\begin{array}{c}ae^{i(\theta+2\phi)},ce^{i\theta}\\ac\end{array}\bigg|q;e^{-i(\theta+\phi)}t\right){}_2\phi_1\left(\begin{array}{c}be^{-i\theta},de^{-i(\theta+2\phi)}\\bd\end{array}\bigg|q;e^{i(\theta+\phi)}t\right)$$

$$=\sum_{n=0}^{\infty}\frac{p_n(x;a,b,c,d;q)}{(ac,bd,q;q)_n}t^n,\ x=\cos(\theta+\phi). \tag{3.4.6}$$

$$_2\phi_1\left(\begin{array}{c}ae^{i(\theta+2\phi)},de^{i\theta}\\ad\end{array}\bigg|q;e^{-i(\theta+\phi)}t\right){}_2\phi_1\left(\begin{array}{c}be^{-i\theta},ce^{-i(\theta+2\phi)}\\bc\end{array}\bigg|q;e^{i(\theta+\phi)}t\right)$$

$$=\sum_{n=0}^{\infty}\frac{p_n(x;a,b,c,d;q)}{(ad,bc,q;q)_n}t^n,\ x=\cos(\theta+\phi). \tag{3.4.7}$$

**References.** [25], [45], [114].

## 3.5 Big $q$-Jacobi

**Definition.**
$$P_n(x;a,b,c;q)={}_3\phi_2\left(\begin{array}{c}q^{-n},abq^{n+1},x\\aq,cq\end{array}\bigg|q;q\right). \tag{3.5.1}$$

**Orthogonality.**

$$\int_{cq}^{aq}\frac{(a^{-1}x,c^{-1}x;q)_\infty}{(x,bc^{-1}x;q)_\infty}P_m(x;a,b,c;q)P_n(x;a,b,c;q)d_qx$$
$$=\ aq(1-q)\frac{(q,a^{-1}c,ac^{-1}q,abq^2;q)_\infty}{(aq,bq,cq,abc^{-1}q;q)_\infty}\times$$
$$\times\frac{(1-abq)}{(1-abq^{2n+1})}\frac{(q,bq,abc^{-1}q;q)_n}{(abq,aq,cq;q)_n}(-acq^2)^nq^{\binom{n}{2}}\delta_{mn}. \tag{3.5.2}$$

**Recurrence relation.**

$$(x-1)P_n(x;a,b,c;q)$$
$$=\ A_nP_{n+1}(x;a,b,c;q)-(A_n+C_n)P_n(x;a,b,c;q)+C_nP_{n-1}(x;a,b,c;q), \tag{3.5.3}$$

where

$$\begin{cases}A_n=\dfrac{(1-aq^{n+1})(1-cq^{n+1})(1-abq^{n+1})}{(1-abq^{2n+1})(1-abq^{2n+2})}\\[2mm]C_n=-acq^{n+1}\dfrac{(1-q^n)(1-bq^n)(1-abc^{-1}q^n)}{(1-abq^{2n})(1-abq^{2n+1})}.\end{cases}$$

**$q$-Difference equation.**

$$q^{-n}(1-q^n)(1-abq^{n+1})x^2y(x)=B(x)y(qx)-[B(x)+D(x)]y(x)+D(x)y(q^{-1}x), \tag{3.5.4}$$

where
$$y(x)=P_n(x;a,b,c;q)$$

and
$$\begin{cases}B(x)=aq(x-1)(bx-c)\\D(x)=(x-aq)(x-cq).\end{cases}$$

**Generating functions.**

$$_2\phi_1\left(\begin{array}{c}aqx^{-1},0\\aq\end{array}\bigg|q;xt\right){}_1\phi_1\left(\begin{array}{c}bc^{-1}x\\bq\end{array}\bigg|q;cqt\right)=\sum_{n=0}^{\infty}\frac{(cq;q)_n}{(bq,q;q)_n}P_n(x;a,b,c;q)t^n. \qquad (3.5.5)$$

$$_2\phi_1\left(\begin{array}{c}cqx^{-1},0\\cq\end{array}\bigg|q;xt\right){}_1\phi_1\left(\begin{array}{c}bc^{-1}x\\abc^{-1}q\end{array}\bigg|q;aqt\right)=\sum_{n=0}^{\infty}\frac{(aq;q)_n}{(abc^{-1}q,q;q)_n}P_n(x;a,b,c;q)t^n. \qquad (3.5.6)$$

**Remarks.** The big $q$-Jacobi polynomials with $c=0$ and the little $q$-Jacobi polynomials defined by (3.12.1) are related in the following way :

$$P_n(x;a,b,0;q)=\frac{(bq;q)_n}{(aq;q)_n}(-1)^n a^n q^{n+\binom{n}{2}}p_n\left(\frac{x}{aq};b,a\bigg|q\right).$$

Sometimes the big $q$-Jacobi polynomials are defined in terms of four parameters instead of three. In fact the polynomials given by the definition

$$P_n(x;a,b,c,d;q)={}_3\phi_2\left(\begin{array}{c}q^{-n},abq^{n+1},ac^{-1}qx\\aq,-ac^{-1}dq\end{array}\bigg|q;q\right)$$

are orthogonal on the interval $[-d,c]$ with respect to the weight function

$$\frac{(c^{-1}qx,-d^{-1}qx;q)_\infty}{(ac^{-1}qx,-bd^{-1}qx;q)_\infty}d_qx.$$

These polynomials are not really different from those defined by (3.5.1) since we have

$$P_n(x;a,b,c,d;q)=P_n(ac^{-1}qx;a,b,-ac^{-1}d;q)$$

and

$$P_n(x;a,b,c;q)=P_n(x;a,b,aq,-cq;q).$$

**References.** [8], [10], [25], [114], [121], [127], [140], [142], [160], [163], [176], [180], [181], [212].

# Special case

## 3.5.1 Big $q$-Legendre

**Definition.** The big $q$-Legendre polynomials are big $q$-Jacobi polynomials with $a=b=1$ :

$$P_n(x;c;q)={}_3\phi_2\left(\begin{array}{c}q^{-n},q^{n+1},x\\q,cq\end{array}\bigg|q;q\right). \qquad (3.5.7)$$

**Orthogonality.**

$$\int_{cq}^{q}P_m(x;c;q)P_n(x;c;q)d_qx=q(1-c)\frac{(1-q)}{(1-q^{2n+1})}\frac{(c^{-1}q;q)_n}{(cq;q)_n}(-cq^2)^n q^{\binom{n}{2}}\delta_{mn}. \qquad (3.5.8)$$

**Recurrence relation.**

$$(x-1)P_n(x;c;q)=A_nP_{n+1}(x;c;q)-(A_n+C_n)P_n(x;c;q)+C_nP_{n-1}(x;c;q), \qquad (3.5.9)$$

where

$$\begin{cases} A_n = \dfrac{(1-q^{n+1})(1-cq^{n+1})}{(1+q^{n+1})(1-q^{2n+1})} \\ \\ C_n = -cq^{n+1}\dfrac{(1-q^n)(1-c^{-1}q^n)}{(1+q^n)(1-q^{2n+1})}. \end{cases}$$

**$q$-Difference equation.**

$$q^{-n}(1-q^n)(1-q^{n+1})x^2 y(x) = B(x)y(qx) - [B(x)+D(x)]y(x) + D(x)y(q^{-1}x), \qquad (3.5.10)$$

where

$$y(x) = P_n(x;c;q)$$

and

$$\begin{cases} B(x) = q(x-1)(x-c) \\ \\ D(x) = (x-q)(x-cq). \end{cases}$$

**Generating functions.**

$$_2\phi_1\left(\begin{matrix} qx^{-1},0 \\ q \end{matrix}\bigg| q;xt\right) {}_1\phi_1\left(\begin{matrix} c^{-1}x \\ q \end{matrix}\bigg| q;cqt\right) = \sum_{n=0}^{\infty} \frac{(cq;q)_n}{(q,q;q)_n} P_n(x;c;q)t^n. \qquad (3.5.11)$$

$$_2\phi_1\left(\begin{matrix} cqx^{-1},0 \\ cq \end{matrix}\bigg| q;xt\right) {}_1\phi_1\left(\begin{matrix} c^{-1}x \\ c^{-1}q \end{matrix}\bigg| q;qt\right) = \sum_{n=0}^{\infty} \frac{P_n(x;c;q)}{(c^{-1}q;q)_n} t^n. \qquad (3.5.12)$$

**References.** [160].

## 3.6 $q$-Hahn

**Definition.**

$$Q_n(q^{-x};\alpha,\beta,N|q) = {}_3\tilde{\phi}_2\left(\begin{matrix} q^{-n},\alpha\beta q^{n+1},q^{-x} \\ \alpha q, q^{-N} \end{matrix}\bigg| q;q\right), \; n=0,1,2,\ldots,N. \qquad (3.6.1)$$

**Orthogonality.**

$$\sum_{x=0}^{N} \frac{(\alpha q, q^{-N};q)_x}{(q,\beta^{-1}q^{-N};q)_x}(\alpha\beta q)^{-x} Q_m(q^{-x};\alpha,\beta,N|q)Q_n(q^{-x};\alpha,\beta,N|q)$$

$$= \frac{(\alpha\beta q^2;q)_N}{(\beta q;q)_N(\alpha q)^N} \frac{(q,\beta q,\alpha\beta q^{N+2};q)_n}{(\alpha\beta q,\alpha q,q^{-N};q)_n} \frac{(1-\alpha\beta q)(-\alpha q)^n}{(1-\alpha\beta q^{2n+1})} q^{\binom{n}{2}-Nn}\delta_{mn}. \qquad (3.6.2)$$

**Recurrence relation.**

$$-\left(1-q^{-x}\right)Q_n(q^{-x}) = A_n Q_{n+1}(q^{-x}) - (A_n+C_n)Q_n(q^{-x}) + C_n Q_{n-1}(q^{-x}), \qquad (3.6.3)$$

where

$$Q_n(q^{-x}) := Q_n(q^{-x};\alpha,\beta,N|q)$$

and

$$\begin{cases} A_n = \dfrac{(1-q^{n-N})(1-\alpha q^{n+1})(1-\alpha\beta q^{n+1})}{(1-\alpha\beta q^{2n+1})(1-\alpha\beta q^{2n+2})} \\ \\ C_n = -\dfrac{\alpha q^n(1-q^n)(1-\beta q^n)(q^{-N}-\alpha\beta q^{n+1})}{(1-\alpha\beta q^{2n})(1-\alpha\beta q^{2n+1})}. \end{cases}$$

**$q$-Difference equation.**

$$q^{-n}(1-q^n)(1-\alpha\beta q^{n+1})y(x) = B(x)y(x+1) - [B(x)+D(x)]y(x) + D(x)y(x-1), \quad (3.6.4)$$

where

$$y(x) = Q_n(q^{-x};\alpha,\beta,N|q)$$

and

$$\begin{cases} B(x) = (1-q^{x-N})(1-\alpha q^{x+1}) \\ D(x) = \alpha q(1-q^x)(\beta - q^{x-N-1}). \end{cases}$$

**Generating functions.**

$$_2\phi_1\left(\begin{matrix} q^{x-N},0 \\ \beta q \end{matrix} \bigg| q; q^{-x}t\right) {}_1\phi_1\left(\begin{matrix} q^{-x} \\ \alpha q \end{matrix} \bigg| q; \alpha qt\right) \simeq \sum_{n=0}^{N} \frac{(q^{-N};q)_n}{(\beta q,q;q)_n} Q_n(q^{-x};\alpha,\beta,N|q)t^n. \quad (3.6.5)$$

$$_2\tilde{\phi}_1\left(\begin{matrix} q^{x-N},0 \\ q^{-N} \end{matrix} \bigg| q; q^{-x}t\right) {}_1\phi_1\left(\begin{matrix} \beta q^{N+1-x} \\ \alpha\beta q^{N+2} \end{matrix} \bigg| q; \alpha qt\right)$$

$$\simeq \sum_{n=0}^{N} \frac{(\alpha q;q)_n}{(\alpha\beta q^{N+2},q;q)_n} Q_n(q^{-x};\alpha,\beta,N|q)t^n. \quad (3.6.6)$$

**Remarks.** The $q$-Hahn polynomials defined by (3.6.1) and the dual $q$-Hahn polynomials given by (3.7.1) are related in the following way :

$$Q_n(q^{-x};\alpha,\beta,N|q) = R_x(\mu(n);\alpha,\beta,N|q),$$

with

$$\mu(n) = q^{-n} + \alpha\beta q^{n+1}$$

or

$$R_n(\mu(x);\gamma,\delta,N|q) = Q_x(q^{-n};\gamma,\delta,N|q),$$

where

$$\mu(x) = q^{-x} + \gamma\delta q^{x+1}.$$

For $x = 0, 1, 2, \ldots, N$ the generating function (3.6.5) can also be written as :

$$_2\phi_1\left(\begin{matrix} q^{x-N},0 \\ \beta q \end{matrix} \bigg| q; q^{-x}t\right) {}_1\phi_1\left(\begin{matrix} q^{-x} \\ \alpha q \end{matrix} \bigg| q; \alpha qt\right) = \sum_{n=0}^{N} \frac{(q^{-N};q)_n}{(\beta q,q;q)_n} Q_n(q^{-x};\alpha,\beta,N|q)t^n.$$

**References.** [10], [25], [43], [45], [88], [111], [114], [121], [142], [145], [158], [160], [180], [197], [215], [217], [218].

## 3.7 Dual $q$-Hahn

**Definition.**

$$R_n(\mu(x);\gamma,\delta,N|q) = {}_3\tilde{\phi}_2\left(\begin{matrix} q^{-n},q^{-x},\gamma\delta q^{x+1} \\ \gamma q, q^{-N} \end{matrix} \bigg| q;q\right), \; n = 0,1,2,\ldots,N, \quad (3.7.1)$$

where

$$\mu(x) := q^{-x} + \gamma\delta q^{x+1}.$$

**Orthogonality.**

$$\sum_{x=0}^{N} \frac{(\gamma q, \gamma\delta q, q^{-N}; q)_x}{(q, \gamma\delta q^{N+2}, \delta q; q)_x} \frac{(1-\gamma\delta q^{2x+1})}{(1-\gamma\delta q)(-\gamma q)^x} q^{Nx-\binom{x}{2}} R_m(\mu(x); \gamma, \delta, N|q) R_n(\mu(x); \gamma, \delta, N|q)$$
$$= \frac{(\gamma\delta q^2; q)_N}{(\delta q; q)_N}(\gamma q)^{-N} \frac{(q, \delta^{-1}q^{-N}; q)_n}{(\gamma q, q^{-N}; q)_n}(\gamma\delta q)^n \delta_{mn}. \tag{3.7.2}$$

**Recurrence relation.**

$$-\left(1-q^{-x}\right)\left(1-\gamma\delta q^{x+1}\right) R_n(\mu(x))$$
$$= A_n R_{n+1}(\mu(x)) - (A_n + C_n) R_n(\mu(x)) + C_n R_{n-1}(\mu(x)), \tag{3.7.3}$$

where

$$R_n(\mu(x)) := R_n(\mu(x); \gamma, \delta, N|q)$$

and

$$\begin{cases} A_n = \left(1-q^{n-N}\right)\left(1-\gamma q^{n+1}\right) \\ C_n = \gamma q \left(1-q^n\right)\left(\delta - q^{n-N-1}\right). \end{cases}$$

**$q$-Difference equation.**

$$q^{-n}(1-q^n)y(x) = B(x)y(x+1) - [B(x)+D(x)]y(x) + D(x)y(x-1), \tag{3.7.4}$$

where

$$y(x) = R_n(\mu(x); \gamma, \delta, N|q)$$

and

$$\begin{cases} B(x) = \dfrac{(1-\gamma q^{x+1})(1-\gamma\delta q^{x+1})(1-q^{x-N})}{(1-\gamma\delta q^{2x+1})(1-\gamma\delta q^{2x+2})} \\ D(x) = -\dfrac{\gamma q^{x-N}(1-q^x)(1-\delta q^x)(1-\gamma\delta q^{x+N+1})}{(1-\gamma\delta q^{2x})(1-\gamma\delta q^{2x+1})}. \end{cases}$$

**Generating functions.**

$$\frac{(\gamma qt; q)_\infty}{(\gamma\delta q^{x+1}t; q)_\infty} {}_2\tilde{\phi}_1\left(\begin{array}{c} q^{x-N}, \gamma\delta q^{x+1} \\ q^{-N} \end{array} \bigg| q; q^{-x}t\right) \simeq \sum_{n=0}^{N} \frac{(\gamma q; q)_n}{(q; q)_n} R_n(\mu(x); \gamma, \delta, N|q) t^n. \tag{3.7.5}$$

$$\frac{(q^{-N}t; q)_\infty}{(q^{-x}t; q)_\infty} {}_2\phi_1\left(\begin{array}{c} q^{-x}, \delta^{-1}q^{-x} \\ \gamma q \end{array} \bigg| q; \gamma\delta q^{x+1}t\right) \simeq \sum_{n=0}^{N} \frac{(q^{-N}; q)_n}{(q; q)_n} R_n(\mu(x); \gamma, \delta, N|q) t^n. \tag{3.7.6}$$

$$\frac{(\gamma\delta qt; q)_\infty}{(\gamma\delta q^{x+1}t; q)_\infty} {}_2\phi_1\left(\begin{array}{c} q^{x-N}, \gamma q^{x+1} \\ \delta^{-1}q^{-N} \end{array} \bigg| q; q^{-x}t\right) \simeq \sum_{n=0}^{N} \frac{(q^{-N}, \gamma q; q)_n}{(\delta^{-1}q^{-N}, q; q)_n} R_n(\mu(x); \gamma, \delta, N|q) t^n. \tag{3.7.7}$$

**Remarks.** The dual $q$-Hahn polynomials defined by (3.7.1) and the $q$-Hahn polynomials given by (3.6.1) are related in the following way :

$$Q_n(q^{-x}; \alpha, \beta, N|q) = R_x(\mu(n); \alpha, \beta, N|q),$$

with

$$\mu(n) = q^{-n} + \alpha\beta q^{n+1}$$

or

$$R_n(\mu(x); \gamma, \delta, N|q) = Q_x(q^{-n}; \gamma, \delta, N|q),$$

where
$$\mu(x) = q^{-x} + \gamma\delta q^{x+1}.$$

For $x = 0, 1, 2, \ldots, N$ the generating function (3.7.6) can also be written as :

$$(q^{-N}t;q)_{N-x} \cdot {}_2\phi_1\left(\begin{array}{c}q^{-x},\delta^{-1}q^{-x}\\ \gamma q\end{array}\Big| q;\gamma\delta q^{x+1}t\right) = \sum_{n=0}^{N}\frac{(q^{-N};q)_n}{(q;q)_n}R_n(\mu(x);\gamma,\delta,N|q)t^n.$$

For $x = 0, 1, 2, \ldots, N$ the generating function (3.7.7) can also be written as :

$$(\gamma\delta qt;q)_x \cdot {}_2\phi_1\left(\begin{array}{c}q^{x-N},\gamma q^{x+1}\\ \delta^{-1}q^{-N}\end{array}\Big| q;q^{-x}t\right) = \sum_{n=0}^{N}\frac{(q^{-N},\gamma q;q)_n}{(\delta^{-1}q^{-N},q;q)_n}R_n(\mu(x);\gamma,\delta,N|q)t^n.$$

**References.** [25], [43], [45], [114], [145], [180], [217].

## 3.8 Al-Salam-Chihara

**Definition.**
$$\begin{aligned}Q_n(x;a,b|q) &= \frac{(ab;q)_n}{a^n}{}_3\phi_2\left(\begin{array}{c}q^{-n},ae^{i\theta},ae^{-i\theta}\\ ab,0\end{array}\Big| q;q\right) \quad (3.8.1)\\ &= (be^{-i\theta};q)_n e^{in\theta}{}_2\phi_1\left(\begin{array}{c}q^{-n},ae^{i\theta}\\ b^{-1}q^{1-n}e^{i\theta}\end{array}\Big| q;b^{-1}qe^{-i\theta}\right), \quad x=\cos\theta.\end{aligned}$$

**Orthogonality.** When $a$ and $b$ are real or complex conjugates and $\max(|a|,|b|) < 1$, then we have the following orthogonality relation

$$\frac{1}{2\pi}\int_{-1}^{1}\frac{w(x)}{\sqrt{1-x^2}}Q_m(x;a,b|q)Q_n(x;a,b|q)dx = \frac{\delta_{mn}}{(q^{n+1},abq^n;q)_\infty}, \quad (3.8.2)$$

where
$$w(x) := w(x;a,b|q) = \left|\frac{(e^{2i\theta};q)_\infty}{(ae^{i\theta},be^{i\theta};q)_\infty}\right|^2 = \frac{h(x,1)h(x,-1)h(x,q^{\frac{1}{2}})h(x,-q^{\frac{1}{2}})}{h(x,a)h(x,b)},$$

with
$$h(x,\alpha) := \prod_{k=0}^{\infty}\left[1 - 2\alpha x q^k + \alpha^2 q^{2k}\right] = \left(\alpha e^{i\theta},\alpha e^{-i\theta};q\right)_\infty, \quad x = \cos\theta.$$

If $a > 1$, $|b| < 1$ and $|ab| < 1$, then we have another orthogonality relation given by :

$$\begin{aligned}&\frac{1}{2\pi}\int_{-1}^{1}\frac{w(x)}{\sqrt{1-x^2}}Q_m(x;a,b|q)Q_n(x;a,b|q)dx + \\ &+ \sum_{\substack{k\\1<aq^k\leq a}} w_k Q_m(x_k;a,b|q)Q_n(x_k;a,b|q) = \frac{\delta_{mn}}{(q^{n+1},abq^n;q)_\infty},\end{aligned} \quad (3.8.3)$$

where $w(x)$ is as before,
$$x_k = \frac{aq^k + (aq^k)^{-1}}{2}$$

and
$$w_k = \frac{(a^{-2};q)_\infty}{(q,ab,a^{-1}b;q)_\infty} \frac{(1-a^2q^{2k})(a^2,ab;q)_k}{(1-a^2)(q,ab^{-1}q;q)_k} q^{-k^2} \left(\frac{1}{a^3b}\right)^k.$$

**Recurrence relation.**

$$2x\tilde{Q}_n(x) = A_n\tilde{Q}_{n+1}(x) + \left[a + a^{-1} - (A_n + C_n)\right]\tilde{Q}_n(x) + C_n\tilde{Q}_{n-1}(x), \qquad (3.8.4)$$

where

$$\tilde{Q}_n(x) := \frac{a^n Q_n(x;a,b|q)}{(ab;q)_n}$$

and

$$\begin{cases} A_n = a^{-1}(1-abq^n) \\ C_n = a(1-q^n). \end{cases}$$

**$q$-Difference equation.**

$$(1-q)^2 D_q\left[\tilde{w}(x;aq^{\frac{1}{2}},bq^{\frac{1}{2}}|q)D_q y(x)\right] + \\ + 4q^{-n+1}(1-q^n)\tilde{w}(x;a,b|q)y(x) = 0, \; y(x) = Q_n(x;a,b|q), \qquad (3.8.5)$$

where

$$\tilde{w}(x;a,b|q) := \frac{w(x;a,b|q)}{\sqrt{1-x^2}}$$

and

$$D_q f(x) := \frac{\delta_q f(x)}{\delta_q x} \text{ with } \delta_q f(e^{i\theta}) = f(q^{\frac{1}{2}}e^{i\theta}) - f(q^{-\frac{1}{2}}e^{i\theta}), \; x = \cos\theta.$$

If we define

$$P_n(z) := \frac{(ab;q)_n}{a^n} {}_3\phi_2\left(\begin{array}{c} q^{-n}, az, az^{-1} \\ ab, 0 \end{array}\bigg| q;q\right)$$

then the $q$-difference equation can also be written in the form

$$q^{-n}(1-q^n)P_n(z) = A(z)P_n(qz) - \left[A(z) + A(z^{-1})\right]P_n(z) + A(z^{-1})P_n(q^{-1}z), \qquad (3.8.6)$$

where

$$A(z) = \frac{(1-az)(1-bz)}{(1-z^2)(1-qz^2)}.$$

**Generating functions.**

$$\frac{1}{(e^{i\theta}t;q)_\infty} {}_2\phi_1\left(\begin{array}{c} ae^{i\theta}, be^{i\theta} \\ ab \end{array}\bigg| q; e^{-i\theta}t\right) = \sum_{n=0}^\infty \frac{Q_n(x;a,b|q)}{(ab,q;q)_n} t^n, \; x = \cos\theta. \qquad (3.8.7)$$

$$\frac{(at,bt;q)_\infty}{(e^{i\theta}t, e^{-i\theta}t;q)_\infty} = \sum_{n=0}^\infty \frac{Q_n(x;a,b|q)}{(q;q)_n} t^n, \; x = \cos\theta. \qquad (3.8.8)$$

**References.** [10], [15], [16], [39], [79], [84].

## 3.9 $q$-Meixner-Pollaczek

**Definition.**

$$P_n(x;a|q) = a^{-n}e^{-in\phi}\frac{(a^2;q)_n}{(q;q)_n}{}_3\phi_2\left(\begin{array}{c}q^{-n},ae^{i(\theta+2\phi)},ae^{-i\theta}\\a^2,0\end{array}\bigg|q;q\right) \quad (3.9.1)$$

$$= \frac{(ae^{-i\theta};q)_n}{(q;q)_n}e^{in(\theta+\phi)}{}_2\phi_1\left(\begin{array}{c}q^{-n},ae^{i\theta}\\a^{-1}q^{1-n}e^{i\theta}\end{array}\bigg|q;qa^{-1}e^{-i(\theta+2\phi)}\right),\ x=\cos(\theta+\phi).$$

**Orthogonality.**

$$\frac{1}{2\pi}\int_{-\pi}^{\pi}w(\cos(\theta+\phi);a|q)P_m(\cos(\theta+\phi);a|q)P_n(\cos(\theta+\phi);a|q)d\theta = \frac{\delta_{mn}}{(q;q)_n(q,a^2q^n;q)_\infty}, \quad (3.9.2)$$

where
$$0 < a < 1$$

and
$$w(x;a|q) = \left|\frac{(e^{2i(\theta+\phi)};q)_\infty}{(ae^{i(\theta+2\phi)},ae^{i\theta};q)_\infty}\right|^2 = \frac{h(x,1)h(x,-1)h(x,q^{\frac{1}{2}})h(x,-q^{\frac{1}{2}})}{h(x,ae^{i\phi})h(x,ae^{-i\phi})},$$

with
$$h(x,\alpha) := \prod_{k=0}^{\infty}\left[1-2\alpha xq^k+\alpha^2q^{2k}\right] = \left(\alpha e^{i(\theta+\phi)},\alpha e^{-i(\theta+\phi)};q\right)_\infty,\ x=\cos(\theta+\phi).$$

**Recurrence relation.**

$$2xP_n(x;a|q) = (1-q^{n+1})P_{n+1}(x;a|q) +$$
$$+ 2aq^n\cos\phi P_n(x;a|q) + (1-a^2q^{n-1})P_{n-1}(x;a|q). \quad (3.9.3)$$

**$q$-Difference equation.**

$$(1-q)^2D_q\left[\tilde{w}(x;aq^{\frac{1}{2}}|q)D_qy(x)\right] + 4q^{-n+1}(1-q^n)\tilde{w}(x;a|q)y(x) = 0,\ y(x) = P_n(x;a|q), \quad (3.9.4)$$

where
$$\tilde{w}(x;a|q) := \frac{w(x;a|q)}{\sqrt{1-x^2}}$$

and
$$D_qf(x) := \frac{\delta_qf(x)}{\delta_qx}\ \text{with}\ \delta_qf(e^{i(\theta+\phi)}) = f(q^{\frac{1}{2}}e^{i(\theta+\phi)}) - f(q^{-\frac{1}{2}}e^{i(\theta+\phi)}),\ x=\cos(\theta+\phi).$$

**Generating functions.**

$$\left|\frac{(ae^{i\phi}t;q)_\infty}{(e^{i(\theta+\phi)}t;q)_\infty}\right|^2 = \sum_{n=0}^{\infty}P_n(x;a|q)t^n,\ x=\cos(\theta+\phi). \quad (3.9.5)$$

$$\frac{1}{(e^{i(\theta+\phi)}t;q)_\infty}{}_2\phi_1\left(\begin{array}{c}ae^{i(\theta+2\phi)},ae^{i\theta}\\a^2\end{array}\bigg|q;e^{-i(\theta+\phi)}t\right) = \sum_{n=0}^{\infty}\frac{P_n(x;a|q)}{(a^2;q)_n}t^n,\ x=\cos(\theta+\phi). \quad (3.9.6)$$

**References.** [10], [16], [39], [45], [72], [126].

## 3.10 Continuous $q$-Jacobi

**Definitions.** If we take $a = q^{\frac{1}{2}\alpha+\frac{1}{4}}$, $b = q^{\frac{1}{2}\alpha+\frac{3}{4}}$, $c = -q^{\frac{1}{2}\beta+\frac{1}{4}}$ and $d = -q^{\frac{1}{2}\beta+\frac{3}{4}}$ in the definition (3.1.1) of the Askey-Wilson polynomials we find after renormalizing

$$P_n^{(\alpha,\beta)}(x|q) = \frac{(q^{\alpha+1};q)_n}{(q;q)_n} {}_4\phi_3\left(\begin{array}{c} q^{-n}, q^{n+\alpha+\beta+1}, q^{\frac{1}{2}\alpha+\frac{1}{4}}e^{i\theta}, q^{\frac{1}{2}\alpha+\frac{1}{4}}e^{-i\theta} \\ q^{\alpha+1}, -q^{\frac{1}{2}(\alpha+\beta+1)}, -q^{\frac{1}{2}(\alpha+\beta+2)} \end{array} \bigg| q;q \right), \; x = \cos\theta. \quad (3.10.1)$$

In [196] M. Rahman takes $a = q^{\frac{1}{2}}$, $b = q^{\alpha+\frac{1}{2}}$, $c = -q^{\beta+\frac{1}{2}}$ and $d = -q^{\frac{1}{2}}$ to obtain after renormalizing

$$P_n^{(\alpha,\beta)}(x;q) = \frac{(q^{\alpha+1}, -q^{\beta+1};q)_n}{(q, -q;q)_n} {}_4\phi_3\left(\begin{array}{c} q^{-n}, q^{n+\alpha+\beta+1}, q^{\frac{1}{2}}e^{i\theta}, q^{\frac{1}{2}}e^{-i\theta} \\ q^{\alpha+1}, -q^{\beta+1}, -q \end{array} \bigg| q;q \right), \; x = \cos\theta. \quad (3.10.2)$$

These two $q$-analogues of the Jacobi polynomials are not really different, since they are connected by the quadratic transformation :

$$P_n^{(\alpha,\beta)}(x|q^2) = \frac{(-q;q)_n}{(-q^{\alpha+\beta+1};q)_n} q^{n\alpha} P_n^{(\alpha,\beta)}(x;q).$$

**Orthogonality.** For $\alpha \geq -\frac{1}{2}$ and $\beta \geq -\frac{1}{2}$ the orthogonality relations are respectively

$$\frac{1}{2\pi} \int_{-1}^{1} \frac{w(x|q)}{\sqrt{1-x^2}} P_m^{(\alpha,\beta)}(x|q) P_n^{(\alpha,\beta)}(x|q) dx$$

$$= \frac{(q^{\frac{1}{2}(\alpha+\beta+2)}, q^{\frac{1}{2}(\alpha+\beta+3)};q)_\infty}{(q, q^{\alpha+1}, q^{\beta+1}, -q^{\frac{1}{2}(\alpha+\beta+1)}, -q^{\frac{1}{2}(\alpha+\beta+2)};q)_\infty} \times$$

$$\times \frac{(1-q^{\alpha+\beta+1})(q^{\alpha+1}, q^{\beta+1}, -q^{\frac{1}{2}(\alpha+\beta+3)};q)_n}{(1-q^{2n+\alpha+\beta+1})(q, q^{\alpha+1}, -q^{\frac{1}{2}(\alpha+\beta+1)};q)_n} q^{(\alpha+\frac{1}{2})n} \delta_{mn}, \quad (3.10.3)$$

where

$$w(x|q) := w(x; q^\alpha, q^\beta|q) = \left| \frac{(e^{2i\theta};q)_\infty}{(q^{\frac{1}{2}\alpha+\frac{1}{4}}e^{i\theta}, q^{\frac{1}{2}\alpha+\frac{3}{4}}e^{i\theta}, -q^{\frac{1}{2}\beta+\frac{1}{4}}e^{i\theta}, -q^{\frac{1}{2}\beta+\frac{3}{4}}e^{i\theta};q)_\infty} \right|^2$$

$$= \left| \frac{(e^{i\theta}, -e^{i\theta};q^{\frac{1}{2}})_\infty}{(q^{\frac{1}{2}\alpha+\frac{1}{4}}e^{i\theta}, -q^{\frac{1}{2}\beta+\frac{1}{4}}e^{i\theta};q^{\frac{1}{2}})_\infty} \right|^2$$

$$= \frac{h(x,1)h(x,-1)h(x,q^{\frac{1}{2}})h(x,-q^{\frac{1}{2}})}{h(x,q^{\frac{1}{2}\alpha+\frac{1}{4}})h(x,q^{\frac{1}{2}\alpha+\frac{3}{4}})h(x,-q^{\frac{1}{2}\beta+\frac{1}{4}})h(x,-q^{\frac{1}{2}\beta+\frac{3}{4}})},$$

with

$$h(x,\alpha) := \prod_{k=0}^{\infty} \left[1 - 2\alpha x q^k + \alpha^2 q^{2k}\right] = \left(\alpha e^{i\theta}, \alpha e^{-i\theta};q\right)_\infty, \; x = \cos\theta$$

and

$$\frac{1}{2\pi} \int_{-1}^{1} \frac{w(x;q)}{\sqrt{1-x^2}} P_m^{(\alpha,\beta)}(x;q) P_n^{(\alpha,\beta)}(x;q) dx$$

$$= \frac{(q^{\alpha+\beta+2};q)_\infty}{(q, -q, q^{\alpha+1}, -q^{\alpha+1}, q^{\beta+1}, -q^{\beta+1}, -q^{\alpha+\beta+1};q)_\infty} \times$$

$$\times \frac{(1-q^{\alpha+\beta+1})(q^{\alpha+1}, q^{\beta+1}, -q^{\alpha+1}, -q^{\beta+1}, -q^{\alpha+\beta+1};q)_n}{(1-q^{2n+\alpha+\beta+1})(q^{\alpha+\beta+1}, q, -q, -q;q)_n} q^n \delta_{mn}, \quad (3.10.4)$$

where

$$w(x;q) := w(x;q^\alpha, q^\beta; q) = \left| \frac{(e^{2i\theta};q)_\infty}{(q^{\alpha+\frac{1}{2}}e^{i\theta}, q^{\frac{1}{2}}e^{i\theta}, -q^{\beta+\frac{1}{2}}e^{i\theta}, -q^{\frac{1}{2}}e^{i\theta};q)_\infty} \right|^2$$

$$= \left| \frac{(e^{i\theta}, -e^{i\theta};q)_\infty}{(q^{\alpha+\frac{1}{2}}e^{i\theta}, -q^{\beta+\frac{1}{2}}e^{i\theta};q)_\infty} \right|^2 = \frac{h(x,1)h(x,-1)}{h(x,q^{\alpha+\frac{1}{2}})h(x,-q^{\beta+\frac{1}{2}})},$$

with

$$h(x,\alpha) := \prod_{k=0}^\infty \left[1 - 2\alpha x q^k + \alpha^2 q^{2k}\right] = \left(\alpha e^{i\theta}, \alpha e^{-i\theta};q\right)_\infty, \quad x = \cos\theta.$$

**Recurrence relations.**

$$2x\tilde{P}_n(x|q) = A_n \tilde{P}_{n+1}(x|q) + \left[q^{\frac{1}{2}\alpha+\frac{1}{4}} + q^{-\frac{1}{2}\alpha-\frac{1}{4}} - (A_n + C_n)\right]\tilde{P}_n(x|q) + C_n \tilde{P}_{n-1}(x|q), \quad (3.10.5)$$

where

$$\tilde{P}_n(x|q) := \frac{(q;q)_n}{(q^{\alpha+1};q)_n} P_n^{(\alpha,\beta)}(x|q)$$

and

$$\begin{cases} A_n = \dfrac{(1-q^{n+\alpha+1})(1-q^{n+\alpha+\beta+1})(1+q^{n+\frac{1}{2}(\alpha+\beta+1)})(1+q^{n+\frac{1}{2}(\alpha+\beta+2)})}{q^{\frac{1}{2}\alpha+\frac{1}{4}}(1-q^{2n+\alpha+\beta+1})(1-q^{2n+\alpha+\beta+2})} \\ \\ C_n = \dfrac{q^{\frac{1}{2}\alpha+\frac{1}{4}}(1-q^n)(1-q^{n+\beta})(1+q^{n+\frac{1}{2}(\alpha+\beta)})(1+q^{n+\frac{1}{2}(\alpha+\beta+1)})}{(1-q^{2n+\alpha+\beta})(1-q^{2n+\alpha+\beta+1})}. \end{cases}$$

$$2x\tilde{P}_n(x;q) = A_n \tilde{P}_{n+1}(x;q) + \left[q^{\frac{1}{2}} + q^{-\frac{1}{2}} - (A_n + C_n)\right]\tilde{P}_n(x;q) + C_n \tilde{P}_{n-1}(x;q), \quad (3.10.6)$$

where

$$\tilde{P}_n(x;q) := \frac{(q,-q;q)_n}{(q^{\alpha+1},-q^{\beta+1};q)_n} P_n^{(\alpha,\beta)}(x;q)$$

and

$$\begin{cases} A_n = \dfrac{(1-q^{n+\alpha+1})(1-q^{n+\alpha+\beta+1})(1+q^{n+1})(1+q^{n+\beta+1})}{q^{\frac{1}{2}}(1-q^{2n+\alpha+\beta+1})(1-q^{2n+\alpha+\beta+2})} \\ \\ C_n = \dfrac{q^{\frac{1}{2}}(1-q^n)(1-q^{n+\beta})(1+q^{n+\alpha})(1+q^{n+\alpha+\beta})}{(1-q^{2n+\alpha+\beta})(1-q^{2n+\alpha+\beta+1})}. \end{cases}$$

**$q$-Difference equations.**

$$(1-q)^2 D_q \left[\tilde{w}(x;q^{\alpha+\frac{1}{2}}, q^{\beta+\frac{1}{2}}|q) D_q y(x)\right] + \lambda_n \tilde{w}(x;q^\alpha, q^\beta|q) y(x) = 0, \; y(x) = P_n^{(\alpha,\beta)}(x|q), \quad (3.10.7)$$

where

$$\tilde{w}(x;q^\alpha, q^\beta|q) := \frac{w(x;q^\alpha, q^\beta|q)}{\sqrt{1-x^2}},$$

$$\lambda_n = 4q^{-n+1}(1-q^n)(1-q^{n+\alpha+\beta+1})$$

and

$$D_q f(x) := \frac{\delta_q f(x)}{\delta_q x} \;\text{ with }\; \delta_q f(e^{i\theta}) = f(q^{\frac{1}{2}}e^{i\theta}) - f(q^{-\frac{1}{2}}e^{i\theta}), \; x = \cos\theta.$$

$$(1-q)^2 D_q \left[\tilde{w}(x;q^{\alpha+\frac{1}{2}}, q^{\beta+\frac{1}{2}};q) D_q y(x)\right] + \lambda_n \tilde{w}(x;q^\alpha, q^\beta;q) y(x) = 0, \; y(x) = P_n^{(\alpha,\beta)}(x;q), \quad (3.10.8)$$

where
$$\tilde{w}(x;q^\alpha,q^\beta;q) := \frac{w(x;q^\alpha,q^\beta;q)}{\sqrt{1-x^2}},$$
$$\lambda_n = 4q^{-n+1}(1-q^n)(1-q^{n+\alpha+\beta+1})$$
and
$$D_q f(x) := \frac{\delta_q f(x)}{\delta_q x} \ \text{ with } \ \delta_q f(e^{i\theta}) = f(q^{\frac{1}{2}}e^{i\theta}) - f(q^{-\frac{1}{2}}e^{i\theta}), \ x = \cos\theta.$$

**Generating functions.**

$$_2\phi_1\left(\begin{array}{c} q^{\frac{1}{2}\alpha+\frac{1}{4}}e^{i\theta}, q^{\frac{1}{2}\alpha+\frac{3}{4}}e^{i\theta} \\ q^{\alpha+1} \end{array}\bigg|q; e^{-i\theta}t\right){}_2\phi_1\left(\begin{array}{c} -q^{\frac{1}{2}\beta+\frac{1}{4}}e^{-i\theta}, -q^{\frac{1}{2}\beta+\frac{3}{4}}e^{-i\theta} \\ q^{\beta+1} \end{array}\bigg|q; e^{i\theta}t\right)$$
$$= \sum_{n=0}^{\infty} \frac{(-q^{\frac{1}{2}(\alpha+\beta+1)}, -q^{\frac{1}{2}(\alpha+\beta+2)};q)_n}{(q^{\alpha+1}, q^{\beta+1};q)_n} \frac{P_n^{(\alpha,\beta)}(x|q)}{q^{(\frac{1}{2}\alpha+\frac{1}{4})n}} t^n, \ x=\cos\theta. \quad (3.10.9)$$

$$_2\phi_1\left(\begin{array}{c} q^{\frac{1}{2}\alpha+\frac{1}{4}}e^{i\theta}, -q^{\frac{1}{2}\beta+\frac{1}{4}}e^{i\theta} \\ -q^{\frac{1}{2}(\alpha+\beta+1)} \end{array}\bigg|q; e^{-i\theta}t\right){}_2\phi_1\left(\begin{array}{c} q^{\frac{1}{2}\alpha+\frac{3}{4}}e^{-i\theta}, -q^{\frac{1}{2}\beta+\frac{3}{4}}e^{-i\theta} \\ -q^{\frac{1}{2}(\alpha+\beta+3)} \end{array}\bigg|q; e^{i\theta}t\right)$$
$$= \sum_{n=0}^{\infty} \frac{(-q^{\frac{1}{2}(\alpha+\beta+2)};q)_n}{(-q^{\frac{1}{2}(\alpha+\beta+3)};q)_n} \frac{P_n^{(\alpha,\beta)}(x|q)}{q^{(\frac{1}{2}\alpha+\frac{1}{4})n}} t^n, \ x=\cos\theta. \quad (3.10.10)$$

$$_2\phi_1\left(\begin{array}{c} q^{\frac{1}{2}\alpha+\frac{1}{4}}e^{i\theta}, -q^{\frac{1}{2}\beta+\frac{3}{4}}e^{i\theta} \\ -q^{\frac{1}{2}(\alpha+\beta+2)} \end{array}\bigg|q; e^{-i\theta}t\right){}_2\phi_1\left(\begin{array}{c} q^{\frac{1}{2}\alpha+\frac{3}{4}}e^{-i\theta}, -q^{\frac{1}{2}\beta+\frac{1}{4}}e^{-i\theta} \\ -q^{\frac{1}{2}(\alpha+\beta+2)} \end{array}\bigg|q; e^{i\theta}t\right)$$
$$= \sum_{n=0}^{\infty} \frac{(-q^{\frac{1}{2}(\alpha+\beta+1)};q)_n}{(-q^{\frac{1}{2}(\alpha+\beta+2)};q)_n} \frac{P_n^{(\alpha,\beta)}(x|q)}{q^{(\frac{1}{2}\alpha+\frac{1}{4})n}} t^n, \ x=\cos\theta. \quad (3.10.11)$$

$$_2\phi_1\left(\begin{array}{c} q^{\frac{1}{2}}e^{i\theta}, q^{\alpha+\frac{1}{2}}e^{i\theta} \\ q^{\alpha+1} \end{array}\bigg|q; e^{-i\theta}t\right){}_2\phi_1\left(\begin{array}{c} -q^{\frac{1}{2}}e^{-i\theta}, -q^{\beta+\frac{1}{2}}e^{-i\theta} \\ q^{\beta+1} \end{array}\bigg|q; e^{i\theta}t\right)$$
$$= \sum_{n=0}^{\infty} \frac{(-q,-q;q)_n}{(q^{\alpha+1}, q^{\beta+1};q)_n} \frac{P_n^{(\alpha,\beta)}(x;q)}{q^{\frac{1}{2}n}} t^n, \ x=\cos\theta. \quad (3.10.12)$$

$$_2\phi_1\left(\begin{array}{c} q^{\frac{1}{2}}e^{i\theta}, -q^{\beta+\frac{1}{2}}e^{i\theta} \\ -q^{\beta+1} \end{array}\bigg|q; e^{-i\theta}t\right){}_2\phi_1\left(\begin{array}{c} -q^{\frac{1}{2}}e^{-i\theta}, q^{\alpha+\frac{1}{2}}e^{-i\theta} \\ -q^{\alpha+1} \end{array}\bigg|q; e^{i\theta}t\right)$$
$$= \sum_{n=0}^{\infty} \frac{(-q,-q;q)_n}{(-q^{\alpha+1}, -q^{\beta+1};q)_n} \frac{P_n^{(\alpha,\beta)}(x;q)}{q^{\frac{1}{2}n}} t^n, \ x=\cos\theta. \quad (3.10.13)$$

$$_2\phi_1\left(\begin{array}{c} q^{\frac{1}{2}}e^{i\theta}, -q^{\frac{1}{2}}e^{i\theta} \\ -q \end{array}\bigg|q; e^{-i\theta}t\right){}_2\phi_1\left(\begin{array}{c} q^{\alpha+\frac{1}{2}}e^{-i\theta}, -q^{\beta+\frac{1}{2}}e^{-i\theta} \\ -q^{\alpha+\beta+1} \end{array}\bigg|q; e^{i\theta}t\right)$$
$$= \sum_{n=0}^{\infty} \frac{(-q;q)_n}{(-q^{\alpha+\beta+1};q)_n} \frac{P_n^{(\alpha,\beta)}(x;q)}{q^{\frac{1}{2}n}} t^n, \ x=\cos\theta. \quad (3.10.14)$$

**Remark.** The continuous $q$-Jacobi polynomials given by (3.10.2) and the continuous $q$-ultraspherical (or Rogers) polynomials given by (3.10.15) are connected by the quadratic transformations :
$$C_{2n}(x;q^\lambda|q) = \frac{(q^\lambda, -q;q)_n}{(q^{\frac{1}{2}}, -q^{\frac{1}{2}};q)_n} q^{-\frac{1}{2}n} P_n^{(\lambda-\frac{1}{2},-\frac{1}{2})}(2x^2-1;q)$$

and
$$C_{2n+1}(x;q^\lambda|q) = \frac{(q^\lambda,-1;q)_{n+1}}{(q^{\frac{1}{2}},-q^{\frac{1}{2}};q)_{n+1}} q^{-\frac{1}{2}n} x P_n^{(\lambda-\frac{1}{2},\frac{1}{2})}(2x^2-1;q).$$

**References.** [45], [112], [114], [136], [179], [180], [196], [198], [199].

# Special cases

### 3.10.1 Continuous $q$-ultraspherical / Rogers

**Definition.** If we set $a = \beta^{\frac{1}{2}}$, $b = \beta^{\frac{1}{2}} q^{\frac{1}{2}}$, $c = -\beta^{\frac{1}{2}}$ and $d = -\beta^{\frac{1}{2}} q^{\frac{1}{2}}$ in the definition (3.1.1) of the Askey-Wilson polynomials and change the normalization we obtain the continuous $q$-ultraspherical (or Rogers) polynomials :

$$\begin{aligned} C_n(x;\beta|q) &= \frac{(\beta^2;q)_n}{(q;q)_n} \beta^{-\frac{1}{2}n} {}_4\phi_3\left(\begin{array}{c} q^{-n}, \beta^2 q^n, \beta^{\frac{1}{2}} e^{i\theta}, \beta^{\frac{1}{2}} e^{-i\theta} \\ \beta q^{\frac{1}{2}}, -\beta, -\beta q^{\frac{1}{2}} \end{array} \bigg| q;q\right) \quad (3.10.15) \\ &= \frac{(\beta^2;q)_n}{(q;q)_n} \beta^{-n} e^{-in\theta} {}_3\phi_2\left(\begin{array}{c} q^{-n}, \beta, \beta e^{2i\theta} \\ \beta^2, 0 \end{array} \bigg| q;q\right) \\ &= \frac{(\beta;q)_n}{(q;q)_n} e^{in\theta} {}_2\phi_1\left(\begin{array}{c} q^{-n}, \beta \\ \beta^{-1} q^{1-n} \end{array} \bigg| q; \beta^{-1} q e^{-2i\theta}\right), \ x = \cos\theta. \end{aligned}$$

**Orthogonality.**

$$\frac{1}{2\pi} \int_{-1}^{1} \frac{w(x)}{\sqrt{1-x^2}} C_m(x;\beta|q) C_n(x;\beta|q) dx = \frac{(\beta,\beta q;q)_\infty}{(\beta^2,q;q)_\infty} \frac{(\beta^2;q)_n}{(q;q)_n} \frac{(1-\beta)}{(1-\beta q^n)} \delta_{mn}, \ |\beta| < 1, \quad (3.10.16)$$

where

$$\begin{aligned} w(x) := w(x;\beta|q) &= \left|\frac{(e^{2i\theta};q)_\infty}{(\beta^{\frac{1}{2}} e^{i\theta}, \beta^{\frac{1}{2}} q^{\frac{1}{2}} e^{i\theta}, -\beta^{\frac{1}{2}} e^{i\theta}, -\beta^{\frac{1}{2}} q^{\frac{1}{2}} e^{i\theta};q)_\infty}\right|^2 = \left|\frac{(e^{2i\theta};q)_\infty}{(\beta e^{2i\theta};q)_\infty}\right|^2 \\ &= \frac{h(x,1) h(x,-1) h(x,q^{\frac{1}{2}}) h(x,-q^{\frac{1}{2}})}{h(x,\beta^{\frac{1}{2}}) h(x,\beta^{\frac{1}{2}} q^{\frac{1}{2}}) h(x,-\beta^{\frac{1}{2}}) h(x,-\beta^{\frac{1}{2}} q^{\frac{1}{2}})}, \end{aligned}$$

with

$$h(x,\alpha) := \prod_{k=0}^{\infty} \left[1 - 2\alpha x q^k + \alpha^2 q^{2k}\right] = \left(\alpha e^{i\theta}, \alpha e^{-i\theta}; q\right)_\infty, \ x = \cos\theta.$$

**Recurrence relation.**

$$2(1-\beta q^n) x C_n(x;\beta|q) = (1-q^{n+1}) C_{n+1}(x;\beta|q) + (1-\beta^2 q^{n-1}) C_{n-1}(x;\beta|q). \quad (3.10.17)$$

**$q$-Difference equation.**

$$(1-q)^2 D_q \left[\tilde{w}(x;\beta q^{\frac{1}{2}}|q) D_q y(x)\right] + \lambda_n \tilde{w}(x;\beta|q) y(x) = 0, \ y(x) = C_n(x;\beta|q), \quad (3.10.18)$$

where

$$\tilde{w}(x;\beta|q) := \frac{w(x;\beta|q)}{\sqrt{1-x^2}},$$

$$\lambda_n = 4 q^{-n+1} (1-q^n)(1-\beta^2 q^n)$$

and

$$D_q f(x) := \frac{\delta_q f(x)}{\delta_q x} \ \text{with} \ \delta_q f(e^{i\theta}) = f(q^{\frac{1}{2}} e^{i\theta}) - f(q^{-\frac{1}{2}} e^{i\theta}), \ x = \cos\theta.$$

**Generating functions.**

$$\frac{(\beta e^{i\theta}t, \beta e^{-i\theta}t; q)_\infty}{(e^{i\theta}t, e^{-i\theta}t; q)_\infty} = \sum_{n=0}^\infty C_n(x; \beta|q)t^n, \ x = \cos\theta. \tag{3.10.19}$$

$$\frac{1}{(e^{-i\theta}t; q)_\infty} {}_2\phi_1\left(\begin{array}{c} \beta, \beta e^{-2i\theta} \\ \beta^2 \end{array} \Big| q; e^{i\theta}t\right) = \sum_{n=0}^\infty \frac{C_n(x; \beta|q)}{(\beta^2; q)_n} t^n, \ x = \cos\theta. \tag{3.10.20}$$

$$(e^{-i\theta}t; q)_\infty \cdot {}_2\phi_1\left(\begin{array}{c} \beta, \beta e^{2i\theta} \\ \beta^2 \end{array} \Big| q; e^{-i\theta}t\right) = \sum_{n=0}^\infty \frac{(-1)^n \beta^n q^{\binom{n}{2}}}{(\beta^2; q)_n} C_n(x; \beta|q)t^n, \ x = \cos\theta. \tag{3.10.21}$$

$${}_2\phi_1\left(\begin{array}{c} \beta^{\frac{1}{2}}e^{i\theta}, \beta^{\frac{1}{2}}q^{\frac{1}{2}}e^{i\theta} \\ \beta q^{\frac{1}{2}} \end{array} \Big| q; e^{-i\theta}t\right) {}_2\phi_1\left(\begin{array}{c} -\beta^{\frac{1}{2}}e^{-i\theta}, -\beta^{\frac{1}{2}}q^{\frac{1}{2}}e^{-i\theta} \\ \beta q^{\frac{1}{2}} \end{array} \Big| q; e^{i\theta}t\right)$$
$$= \sum_{n=0}^\infty \frac{(-\beta, -\beta q^{\frac{1}{2}}; q)_n}{(\beta^2, \beta q^{\frac{1}{2}}; q)_n} C_n(x; \beta|q)t^n, \ x = \cos\theta. \tag{3.10.22}$$

$${}_2\phi_1\left(\begin{array}{c} \beta^{\frac{1}{2}}e^{i\theta}, -\beta^{\frac{1}{2}}e^{i\theta} \\ -\beta \end{array} \Big| q; e^{-i\theta}t\right) {}_2\phi_1\left(\begin{array}{c} \beta^{\frac{1}{2}}q^{\frac{1}{2}}e^{-i\theta}, -\beta^{\frac{1}{2}}q^{\frac{1}{2}}e^{-i\theta} \\ -\beta q \end{array} \Big| q; e^{i\theta}t\right)$$
$$= \sum_{n=0}^\infty \frac{(\beta q^{\frac{1}{2}}, -\beta q^{\frac{1}{2}}; q)_n}{(\beta^2, -\beta q; q)_n} C_n(x; \beta|q)t^n, \ x = \cos\theta. \tag{3.10.23}$$

$${}_2\phi_1\left(\begin{array}{c} \beta^{\frac{1}{2}}e^{i\theta}, -\beta^{\frac{1}{2}}q^{\frac{1}{2}}e^{i\theta} \\ -\beta q^{\frac{1}{2}} \end{array} \Big| q; e^{-i\theta}t\right) {}_2\phi_1\left(\begin{array}{c} \beta^{\frac{1}{2}}q^{\frac{1}{2}}e^{-i\theta}, -\beta^{\frac{1}{2}}e^{-i\theta} \\ -\beta q^{\frac{1}{2}} \end{array} \Big| q; e^{i\theta}t\right)$$
$$= \sum_{n=0}^\infty \frac{(-\beta, \beta q^{\frac{1}{2}}; q)_n}{(\beta^2, -\beta q^{\frac{1}{2}}; q)_n} C_n(x; \beta|q)t^n, \ x = \cos\theta. \tag{3.10.24}$$

**Remarks.** The continuous $q$-ultraspherical (or Rogers) polynomials can also be written as :

$$C_n(x; \beta|q) = \sum_{k=0}^n \frac{(\beta; q)_k (\beta; q)_{n-k}}{(q; q)_k (q; q)_{n-k}} e^{i(n-2k)\theta}, \ x = \cos\theta.$$

They can be obtained from the continuous $q$-Jacobi polynomials defined by (3.10.1) in the following way. Set $\beta = \alpha$ in the definition (3.10.1) and change $q^{\alpha+\frac{1}{2}}$ by $\beta$ and we find the continuous $q$-ultraspherical (or Rogers) polynomials with a different normalization. We have

$$P_n^{(\alpha,\alpha)}(x|q) \stackrel{q^{\alpha+\frac{1}{2}} \to \beta}{\longrightarrow} \frac{(\beta q^{\frac{1}{2}}; q)_n}{(\beta^2; q)_n} \beta^{\frac{1}{2}n} C_n(x; \beta|q).$$

If we set $\beta = q^{\alpha+\frac{1}{2}}$ in the definition (3.10.15) of the $q$-ultraspherical (or Rogers) polynomials we find the continuous $q$-Jacobi polynomials given by (3.10.1) with $\beta = \alpha$. In fact we have

$$C_n\left(x; q^{\alpha+\frac{1}{2}}\Big| q\right) = \frac{(q^{2\alpha+1}; q)_n}{(q^{\alpha+1}; q)_n q^{(\frac{1}{2}\alpha+\frac{1}{4})n}} P_n^{(\alpha,\alpha)}(x|q).$$

If we change $q$ to $q^{-1}$ we find

$$C_n(x; \beta|q^{-1}) = (\beta q)^n C_n(x; \beta^{-1}|q).$$

The special case $\beta = q$ of the continuous $q$-ultraspherical (or Rogers) polynomials equals the Chebyshev polynomials of the second kind defined by (1.8.21). In fact we have

$$C_n(x;q|q) = \frac{\sin(n+1)\theta}{\sin\theta} = U_n(x), \ x = \cos\theta.$$

The limit case $\beta \to 1$ leads to the Chebyshev polynomials of the first kind given by (1.8.20) in the following way :

$$\lim_{\beta \to 1} \frac{1-q^n}{2(1-\beta)} C_n(x;\beta|q) = \cos n\theta = T_n(x), \ x = \cos\theta, \ n = 1,2,3,\ldots.$$

The continuous $q$-Jacobi polynomials given by (3.10.2) and the continuous $q$-ultraspherical (or Rogers) polynomials given by (3.10.15) are connected by the quadratic transformations :

$$C_{2n}(x;q^\lambda|q) = \frac{(q^\lambda,-q;q)_n}{(q^{\frac{1}{2}},-q^{\frac{1}{2}};q)_n} q^{-\frac{1}{2}n} P_n^{(\lambda-\frac{1}{2},-\frac{1}{2})}(2x^2-1;q)$$

and

$$C_{2n+1}(x;q^\lambda|q) = \frac{(q^\lambda,-1;q)_{n+1}}{(q^{\frac{1}{2}},-q^{\frac{1}{2}};q)_{n+1}} q^{-\frac{1}{2}n} x P_n^{(\lambda-\frac{1}{2},\frac{1}{2})}(2x^2-1;q).$$

Finally we remark that the continuous $q$-ultraspherical (or Rogers) polynomials are related to the continuous $q$-Legendre polynomials defined by (3.10.25) in the following way :

$$C_n(x;q|q^2) = q^{-\frac{1}{2}n} P_n(x;q).$$

**References.** [10], [11], [12], [25], [31], [32], [37], [38], [39], [41], [45], [60], [62], [63], [107], [108], [109], [110], [112], [113], [114], [127], [137], [159], [179], [180], [181], [202], [204], [207], [208], [209].

### 3.10.2 Continuous $q$-Legendre

**Definition.** The continuous $q$-Legendre polynomials are continuous $q$-Jacobi polynomials with $\alpha = \beta = 0$. If we set $\alpha = \beta = 0$ in the definition (3.10.2) of the continuous $q$-Jacobi polynomials we obtain

$$P_n(x;q) = {}_4\phi_3\left(\begin{array}{c} q^{-n}, q^{n+1}, q^{\frac{1}{2}}e^{i\theta}, q^{\frac{1}{2}}e^{-i\theta} \\ q, -q, -q \end{array} \bigg| q;q\right), \ x = \cos\theta. \tag{3.10.25}$$

If we set $\alpha = \beta = 0$ in the definition (3.10.1) we find

$$P_n(x|q) = {}_4\phi_3\left(\begin{array}{c} q^{-n}, q^{n+1}, q^{\frac{1}{4}}e^{i\theta}, q^{\frac{1}{4}}e^{-i\theta} \\ q, -q^{\frac{1}{2}}, -q \end{array} \bigg| q;q\right), \ x = \cos\theta,$$

but these are not really different in view of the quadratic transformation

$$P_n(x|q^2) = P_n(x;q).$$

**Orthogonality.**

$$\frac{1}{2\pi} \int_{-1}^{1} \frac{w(x)}{\sqrt{1-x^2}} P_m(x;q) P_n(x;q) dx = \frac{(q;q)_{2n}(q^{2n+2};q)_\infty q^n}{(-q;q)_\infty^4 (q;q)_\infty^3} \delta_{mn}, \tag{3.10.26}$$

where

$$w(x) = \left|\frac{(e^{2i\theta};q)_\infty}{(q^{\frac{1}{2}}e^{i\theta}, q^{\frac{1}{2}}e^{i\theta}, -q^{\frac{1}{2}}e^{i\theta}, -q^{\frac{1}{2}}e^{i\theta};q)_\infty}\right|^2 = \left|\frac{(e^{i\theta},-e^{i\theta};q)_\infty}{(q^{\frac{1}{2}}e^{i\theta},-q^{\frac{1}{2}}e^{i\theta};q)_\infty}\right|^2 = \frac{h(x,1)h(x,-1)}{h(x,q^{\frac{1}{2}})h(x,-q^{\frac{1}{2}})},$$

with
$$h(x,\alpha) := \prod_{k=0}^{\infty} \left[1 - 2\alpha x q^k + \alpha^2 q^{2k}\right] = \left(\alpha e^{i\theta}, \alpha e^{-i\theta}; q\right)_{\infty}, \ x = \cos\theta.$$

**Recurrence relation.**
$$2(1-q^{2n+1})xP_n(x;q) = q^{-\frac{1}{2}}(1-q^{2n+2})P_{n+1}(x;q) + q^{\frac{1}{2}}(1-q^{2n})P_{n-1}(x;q). \tag{3.10.27}$$

**$q$-Difference equation.**
$$(1-q)^2 D_q \left[\tilde{w}(x;q;q) D_q y(x)\right] + \lambda_n \tilde{w}(x;q^{\frac{1}{2}};q)y(x) = 0, \ y(x) = P_n(x;q), \tag{3.10.28}$$

where
$$\lambda_n = 4q^{-n+1}(1-q^n)(1-q^{n+1})$$

and
$$\tilde{w}(x;a;q) := \frac{w(x;a;q)}{\sqrt{1-x^2}},$$

where
$$w(x;a;q) = \left|\frac{(e^{2i\theta};q)_\infty}{(ae^{i\theta},ae^{i\theta},-ae^{i\theta},-ae^{i\theta};q)_\infty}\right|^2 = \frac{h(x,1)h(x,-1)h(x,q^{\frac{1}{2}})h(x,-q^{\frac{1}{2}})}{h(x,a)h(x,a)h(x,-a)h(x,-a)},$$

with
$$h(x,\alpha) := \prod_{k=0}^{\infty} \left[1 - 2\alpha x q^k + \alpha^2 q^{2k}\right] = \left(\alpha e^{i\theta}, \alpha e^{-i\theta}; q\right)_{\infty}, \ x = \cos\theta$$

and
$$D_q f(x) := \frac{\delta_q f(x)}{\delta_q x} \ \text{with} \ \delta_q f(e^{i\theta}) = f(q^{\frac{1}{2}}e^{i\theta}) - f(q^{-\frac{1}{2}}e^{i\theta}), \ x = \cos\theta.$$

**Generating functions.**
$$_2\phi_1\left(\begin{array}{c} q^{\frac{1}{2}}e^{i\theta}, -q^{\frac{1}{2}}e^{i\theta} \\ -q \end{array} \Bigg| q; e^{-i\theta}t\right) {}_2\phi_1\left(\begin{array}{c} q^{\frac{1}{2}}e^{-i\theta}, -q^{\frac{1}{2}}e^{-i\theta} \\ -q \end{array} \Bigg| q; e^{i\theta}t\right)$$
$$= \sum_{n=0}^{\infty} \frac{P_n(x;q)}{q^{\frac{1}{2}n}} t^n, \ x = \cos\theta. \tag{3.10.29}$$

$$_2\phi_1\left(\begin{array}{c} q^{\frac{1}{2}}e^{i\theta}, q^{\frac{1}{2}}e^{i\theta} \\ q \end{array} \Bigg| q; e^{-i\theta}t\right) {}_2\phi_1\left(\begin{array}{c} -q^{\frac{1}{2}}e^{-i\theta}, -q^{\frac{1}{2}}e^{-i\theta} \\ q \end{array} \Bigg| q; e^{i\theta}t\right)$$
$$= \sum_{n=0}^{\infty} \frac{(-q,-q;q)_n}{(q,q;q)_n} \frac{P_n(x;q)}{q^{\frac{1}{2}n}} t^n, \ x = \cos\theta. \tag{3.10.30}$$

$$_2\phi_1\left(\begin{array}{c} q^{\frac{1}{2}}e^{i\theta}, q^{\frac{3}{2}}e^{i\theta} \\ q^2 \end{array} \Bigg| q^2; e^{-i\theta}t\right) {}_2\phi_1\left(\begin{array}{c} -q^{\frac{1}{2}}e^{-i\theta}, -q^{\frac{3}{2}}e^{-i\theta} \\ q^2 \end{array} \Bigg| q^2; e^{i\theta}t\right)$$
$$= \sum_{n=0}^{\infty} \frac{(-q^{n+1};q)_n}{(q,q,-q;q)_n} \frac{P_n(x;q)}{q^{\frac{1}{2}n}} t^n, \ x = \cos\theta. \tag{3.10.31}$$

$$_2\phi_1\left(\begin{array}{c} q^{\frac{1}{2}}e^{i\theta}, -q^{\frac{1}{2}}e^{i\theta} \\ -q \end{array} \Bigg| q^2; e^{-i\theta}t\right) {}_2\phi_1\left(\begin{array}{c} q^{\frac{3}{2}}e^{-i\theta}, -q^{\frac{3}{2}}e^{-i\theta} \\ -q^3 \end{array} \Bigg| q^2; e^{i\theta}t\right)$$
$$= \sum_{n=0}^{\infty} \frac{(-q^2;q^2)_n}{(-q^3;q^2)_n} \frac{P_n(x;q)}{q^{\frac{1}{2}n}} t^n, \ x = \cos\theta. \tag{3.10.32}$$

$$_2\phi_1\left(\begin{array}{c}q^{\frac{1}{2}}e^{i\theta},-q^{\frac{3}{2}}e^{i\theta}\\-q^2\end{array}\bigg|\,q^2;e^{-i\theta}t\right){}_2\phi_1\left(\begin{array}{c}q^{\frac{3}{2}}e^{-i\theta},-q^{\frac{1}{2}}e^{-i\theta}\\-q^2\end{array}\bigg|\,q^2;e^{i\theta}t\right)$$

$$=\sum_{n=0}^{\infty}\frac{(-q;q^2)_n}{(-q^2;q^2)_n}\frac{P_n(x;q)}{q^{\frac{1}{2}n}}t^n,\ x=\cos\theta. \tag{3.10.33}$$

**Remarks.** The continuous $q$-Legendre polynomials can also be written as :

$$P_n(x;q)=q^{\frac{1}{2}n}\sum_{k=0}^{n}\frac{(q;q^2)_k(q;q^2)_{n-k}}{(q^2;q^2)_k(q^2;q^2)_{n-k}}e^{i(n-2k)\theta},\ x=\cos\theta.$$

They are related to the continuous $q$-ultraspherical (or Rogers) polynomials given by (3.10.15) in the following way :

$$P_n(x;q)=q^{\frac{1}{2}n}C_n(x;q|q^2).$$

**References.** [157], [160], [163].

## 3.11 Big $q$-Laguerre

**Definition.**

$$\begin{aligned}P_n(x;a,b;q)&=\ {}_3\phi_2\left(\begin{array}{c}q^{-n},0,x\\aq,bq\end{array}\bigg|\,q;q\right) \tag{3.11.1}\\ &=\ \frac{1}{(b^{-1}q^{-n};q)_n}{}_2\phi_1\left(\begin{array}{c}q^{-n},aqx^{-1}\\aq\end{array}\bigg|\,q;\frac{x}{b}\right).\end{aligned}$$

**Orthogonality.**

$$\int_{bq}^{aq}\frac{(a^{-1}x,b^{-1}x;q)_\infty}{(x;q)_\infty}P_m(x;a,b;q)P_n(x;a,b;q)d_qx$$

$$=\ aq(1-q)\frac{(q,a^{-1}b,ab^{-1}q;q)_\infty}{(aq,bq;q)_\infty}\frac{(q;q)_n}{(aq,bq;q)_n}(-abq^2)^n q^{\binom{n}{2}}\delta_{mn}. \tag{3.11.2}$$

**Recurrence relation.**

$$(x-1)P_n(x;a,b;q)=A_nP_{n+1}(x;a,b;q)-(A_n+C_n)P_n(x;a,b;q)+C_nP_{n-1}(x;a,b;q), \tag{3.11.3}$$

where

$$\begin{cases}A_n=(1-aq^{n+1})(1-bq^{n+1})\\ C_n=-abq^{n+1}(1-q^n).\end{cases}$$

**$q$-Difference equation.**

$$q^{-n}(1-q^n)x^2y(x)=B(x)y(qx)-[B(x)+D(x)]y(x)+D(x)y(q^{-1}x), \tag{3.11.4}$$

where

$$y(x)=P_n(x;a,b;q)$$

and

$$\begin{cases}B(x)=abq(1-x)\\ D(x)=(x-aq)(x-bq).\end{cases}$$

**Generating functions.**

$$(bqt;q)_\infty \cdot {}_2\phi_1\left(\begin{array}{c} aqx^{-1},0 \\ aq \end{array} \bigg| q;xt\right) = \sum_{n=0}^{\infty} \frac{(bq;q)_n}{(q;q)_n} P_n(x;a,b;q)t^n. \qquad (3.11.5)$$

$$(t;q)_\infty \cdot {}_3\phi_2\left(\begin{array}{c} 0,0,x \\ aq,bq \end{array} \bigg| q;t\right) = \sum_{n=0}^{\infty} \frac{(-1)^n q^{\binom{n}{2}}}{(q;q)_n} P_n(x;a,b;q)t^n. \qquad (3.11.6)$$

**Remark.** The big $q$-Laguerre polynomials defined by (3.11.1) and the affine $q$-Krawtchouk polynomials given by (3.16.1) are related in the following way :

$$K_n^{Aff}(q^{-x};p,N;q) = P_n(q^{-x};p,q^{-N-1};q).$$

**References.** [10], [23].

## 3.12 Little $q$-Jacobi

**Definition.**
$$p_n(x;a,b|q) = {}_2\phi_1\left(\begin{array}{c} q^{-n},abq^{n+1} \\ aq \end{array} \bigg| q;qx\right). \qquad (3.12.1)$$

**Orthogonality.**

$$\sum_{k=0}^{\infty} \frac{(bq;q)_k}{(q;q)_k}(aq)^k p_m(q^k;a,b|q)p_n(q^k;a,b|q)$$
$$= \frac{(abq^2;q)_\infty}{(aq;q)_\infty} \frac{(1-abq)(aq)^n}{(1-abq^{2n+1})} \frac{(q,bq;q)_n}{(aq,abq;q)_n}\delta_{mn},\ 0 < aq < 1,\ b < q^{-1}. \qquad (3.12.2)$$

**Recurrence relation.**

$$-xp_n(x;a,b|q) = A_n p_{n+1}(x;a,b|q) - (A_n+C_n)p_n(x;a,b|q) + C_n p_{n-1}(x;a,b|q), \qquad (3.12.3)$$

where
$$\begin{cases} A_n = q^n \dfrac{(1-aq^{n+1})(1-abq^{n+1})}{(1-abq^{2n+1})(1-abq^{2n+2})} \\ C_n = aq^n \dfrac{(1-q^n)(1-bq^n)}{(1-abq^{2n})(1-abq^{2n+1})}. \end{cases}$$

$q$-**Difference equation.**

$$q^{-n}(1-q^n)(1-abq^{n+1})xy(x)$$
$$= B(x)y(qx) - [B(x)+D(x)]y(x) + D(x)y(q^{-1}x),\ y(x) = p_n(x;a,b|q), \qquad (3.12.4)$$

where
$$\begin{cases} B(x) = a(bqx-1) \\ D(x) = x-1. \end{cases}$$

**Generating function.**

$${}_2\phi_1\left(\begin{array}{c} 0,0 \\ aq \end{array} \bigg| q;xt\right) {}_1\phi_1\left(\begin{array}{c} bqx \\ bq \end{array} \bigg| q;t\right) = \sum_{n=0}^{\infty} \frac{(-1)^n q^{\binom{n}{2}}}{(bq,q;q)_n} p_n(x;a,b|q)t^n. \qquad (3.12.5)$$

**Remarks.** The little $q$-Jacobi polynomials defined by (3.12.1) and the big $q$-Jacobi polynomials given by (3.5.1) are related in the following way :

$$p_n(x;a,b|q) = \frac{(bq;q)_n}{(aq;q)_n}(-1)^n b^{-n} q^{-n-\binom{n}{2}} P_n(bqx;b,a,0;q).$$

The little $q$-Jacobi polynomials and the $q$-Meixner polynomials defined by (3.13.1) are related in the following way :

$$M_n(q^{-x};b,c;q) = p_n(-c^{-1}q^n;b,b^{-1}q^{-n-x-1}|q).$$

**References.** [8], [10], [18], [19], [24], [25], [31], [111], [114], [121], [127], [134], [140], [145], [158], [160], [161], [163], [172], [176], [180], [197], [212], [214].

# Special case

## 3.12.1 Little $q$-Legendre

**Definition.** The little $q$-Legendre polynomials are little $q$-Jacobi polynomials with $a = b = 1$ :

$$p_n(x|q) = {}_2\phi_1\left(\begin{array}{c}q^{-n},q^{n+1}\\q\end{array}\bigg|\,q;qx\right). \tag{3.12.6}$$

**Orthogonality.**

$$\sum_{k=0}^{\infty} q^k p_m(q^k|q) p_n(q^k|q) = \frac{q^n}{(1-q^{2n+1})}\delta_{mn}. \tag{3.12.7}$$

**Recurrence relation.**

$$-xp_n(x|q) = A_n p_{n+1}(x|q) - (A_n + C_n) p_n(x|q) + C_n p_{n-1}(x|q), \tag{3.12.8}$$

where

$$\begin{cases} A_n = q^n \dfrac{(1-q^{n+1})}{(1+q^{n+1})(1-q^{2n+1})} \\[2ex] C_n = q^n \dfrac{(1-q^n)}{(1+q^n)(1-q^{2n+1})}. \end{cases}$$

**$q$-Difference equation.**

$$q^{-n}(1-q^n)(1-q^{n+1})xy(x) = B(x)y(qx) - [B(x)+D(x)]y(x) + D(x)y(q^{-1}x), \tag{3.12.9}$$

where

$$y(x) = p_n(x|q)$$

and

$$\begin{cases} B(x) = qx - 1 \\ D(x) = x - 1. \end{cases}$$

**Generating function.**

$${}_2\phi_1\left(\begin{array}{c}0,0\\q\end{array}\bigg|\,q;xt\right){}_1\phi_1\left(\begin{array}{c}qx\\q\end{array}\bigg|\,q;t\right) = \sum_{n=0}^{\infty}\frac{(-1)^n q^{\binom{n}{2}}}{(q,q;q)_n}p_n(x|q)t^n. \tag{3.12.10}$$

**References.** [160], [161], [201], [221].

## 3.13   $q$-Meixner

**Definition.**
$$M_n(q^{-x};b,c;q) = {}_2\phi_1\left(\begin{array}{c} q^{-n},q^{-x} \\ bq \end{array} \Big| \, q;-\frac{q^{n+1}}{c}\right). \tag{3.13.1}$$

**Orthogonality.**
$$\sum_{x=0}^{\infty} \frac{(bq;q)_x}{(q,-bcq;q)_x} c^x q^{\binom{x}{2}} M_m(q^{-x};b,c;q) M_n(q^{-x};b,c;q)$$
$$= \frac{(-c;q)_\infty}{(-bcq;q)_\infty} \frac{(q,-c^{-1}q;q)_n}{(bq;q)_n} q^{-n}\delta_{mn}, \; 0 < bq < 1, \; c > 0. \tag{3.13.2}$$

**Recurrence relation.**
$$q^{2n+1}(1-q^{-x})M_n(q^{-x}) = c(1-bq^{n+1})M_{n+1}(q^{-x}) +$$
$$- \left[c(1-bq^{n+1}) + q(1-q^n)(c+q^n)\right] M_n(q^{-x}) + q(1-q^n)(c+q^n)M_{n-1}(q^{-x}), \tag{3.13.3}$$

where
$$M_n(q^{-x}) := M_n(q^{-x};b,c;q).$$

**$q$-Difference equation.**
$$-(1-q^n)y(x) = B(x)y(x+1) - [B(x)+D(x)]y(x) + D(x)y(x-1), \tag{3.13.4}$$

where
$$y(x) = M_n(q^{-x};b,c;q)$$

and
$$\begin{cases} B(x) = cq^x(1-bq^{x+1}) \\ D(x) = (1-q^x)(1+bcq^x). \end{cases}$$

**Generating functions.**
$$\frac{1}{(t;q)_\infty} {}_1\phi_1\left(\begin{array}{c} q^{-x} \\ bq \end{array} \Big| \, q;-\frac{qt}{c}\right) = \sum_{n=0}^{\infty} \frac{M_n(q^{-x};b,c;q)}{(q;q)_n} t^n. \tag{3.13.5}$$

$$\frac{1}{(t;q)_\infty} {}_1\phi_1\left(\begin{array}{c} -b^{-1}c^{-1}q^{-x} \\ -c^{-1}q \end{array} \Big| \, q;bqt\right) = \sum_{n=0}^{\infty} \frac{(bq;q)_n}{(-c^{-1}q,q;q)_n} M_n(q^{-x};b,c;q)t^n. \tag{3.13.6}$$

**Remarks.** The $q$-Meixner polynomials defined by (3.13.1) and the little $q$-Jacobi polynomials given by (3.12.1) are related in the following way :
$$M_n(q^{-x};b,c;q) = p_n(-c^{-1}q^n;b,b^{-1}q^{-n-x-1}|q).$$

The $q$-Meixner polynomials and the quantum $q$-Krawtchouk polynomials defined by (3.14.1) are related in the following way :
$$K_n^{qtm}(q^{-x};p,N;q) = M_n(q^{-x};q^{-N-1},-p^{-1};q).$$

**References.** [10], [22], [23], [114], [121], [180].

## 3.14 Quantum $q$-Krawtchouk

**Definition.**
$$K_n^{qtm}(q^{-x};p,N;q) = {}_2\tilde{\phi}_1\left(\begin{array}{c}q^{-n},q^{-x}\\q^{-N}\end{array}\bigg|\,q;pq^{n+1}\right),\ n=0,1,2,\ldots,N. \tag{3.14.1}$$

**Orthogonality.**
$$\sum_{x=0}^{N}\frac{(pq;q)_{N-x}}{(q;q)_x(q;q)_{N-x}}(-1)^{N-x}q^{\binom{x}{2}}K_m^{qtm}(q^{-x};p,N;q)K_n^{qtm}(q^{-x};p,N;q)$$
$$= \frac{(-1)^n p^N (q;q)_{N-n}(q,pq;q)_n}{(q,q;q)_N}q^{\binom{N+1}{2}-\binom{n+1}{2}+Nn}\delta_{mn}. \tag{3.14.2}$$

**Recurrence relation.**
$$-pq^{2n+1}(1-q^{-x})K_n^{qtm}(q^{-x}) = (1-q^{n-N})K_{n+1}^{qtm}(q^{-x}) +$$
$$-\left[(1-q^{n-N})+q(1-q^n)(1-pq^n)\right]K_n^{qtm}(q^{-x})+q(1-q^n)(1-pq^n)K_{n-1}^{qtm}(q^{-x}), \tag{3.14.3}$$

where
$$K_n^{qtm}(q^{-x}) := K_n^{qtm}(q^{-x};p,N;q).$$

**$q$-Difference equation.**
$$-p(1-q^n)y(x) = B(x)y(x+1) - [B(x)+D(x)]y(x) + D(x)y(x-1), \tag{3.14.4}$$

where
$$y(x) = K_n^{qtm}(q^{-x};p,N;q)$$

and
$$\begin{cases} B(x) = -q^x(1-q^{x-N}) \\ D(x) = (1-q^x)(p-q^{x-N-1}). \end{cases}$$

**Generating function.**
$$\frac{(q^{-x}t;q)_\infty}{(t;q)_\infty}{}_2\phi_1\left(\begin{array}{c}q^{x-N},0\\pq\end{array}\bigg|\,q;q^{-x}t\right) \simeq \sum_{n=0}^{N}\frac{(q^{-N};q)_n}{(pq,q;q)_n}K_n^{qtm}(q^{-x};p,N;q)t^n. \tag{3.14.5}$$

**Remarks.** The quantum $q$-Krawtchouk polynomials defined by (3.14.1) and the $q$-Meixner polynomials given by (3.13.1) are related in the following way :

$$K_n^{qtm}(q^{-x};p,N;q) = M_n(q^{-x};q^{-N-1},-p^{-1};q).$$

The quantum $q$-Krawtchouk polynomials are related to the affine $q$-Krawtchouk polynomials defined by (3.16.1) by the transformation $q \leftrightarrow q^{-1}$ in the following way :

$$K_n^{qtm}(q^x;p,N;q^{-1}) = (p^{-1}q;q)_n\left(-\frac{p}{q}\right)^n q^{-\binom{n}{2}}K_n^{Aff}(q^{x-N};p^{-1},N;q).$$

For $x = 0,1,2,\ldots,N$ the generating function (3.14.5) can also be written as :

$$(q^{-x}t;q)_x \cdot {}_2\phi_1\left(\begin{array}{c}q^{x-N},0\\pq\end{array}\bigg|\,q;q^{-x}t\right) = \sum_{n=0}^{N}\frac{(q^{-N};q)_n}{(pq,q;q)_n}K_n^{qtm}(q^{-x};p,N;q)t^n.$$

**References.** [114], [158], [160].

## 3.15 $q$-Krawtchouk

**Definition.**

$$K_n(q^{-x};p,N;q) = {}_3\tilde{\phi}_2\left(\begin{array}{c}q^{-n},q^{-x},-pq^n\\q^{-N},0\end{array}\bigg|q;q\right) \qquad (3.15.1)$$

$$= \frac{(q^{x-N};q)_n}{(q^{-N};q)_n q^{nx}}{}_2\phi_1\left(\begin{array}{c}q^{-n},q^{-x}\\q^{N-x-n+1}\end{array}\bigg|q;-pq^{n+N+1}\right),\ n=0,1,2,\ldots,N.$$

**Orthogonality.**

$$\sum_{x=0}^{N}\frac{(q^{-N};q)_x}{(q;q)_x}(-p)^{-x}K_m(q^{-x};p,N;q)K_n(q^{-x};p,N;q)$$

$$= \frac{(q,-pq^{N+1};q)_n}{(-p,q^{-N};q)_n}\frac{(1+p)}{(1+pq^{2n})}\times$$

$$\times (-pq;q)_N p^{-N}q^{-\binom{N+1}{2}}\left(-pq^{-N}\right)^n q^{n^2}\delta_{mn}. \qquad (3.15.2)$$

**Recurrence relation.**

$$-\left(1-q^{-x}\right)K_n(q^{-x}) = A_n K_{n+1}(q^{-x}) - (A_n+C_n)K_n(q^{-x}) + C_n K_{n-1}(q^{-x}), \qquad (3.15.3)$$

where

$$K_n(q^{-x}) := K_n(q^{-x};p,N;q)$$

and

$$\begin{cases}A_n = \dfrac{(1-q^{n-N})(1+pq^n)}{(1+pq^{2n})(1+pq^{2n+1})}\\[2ex] C_n = -pq^{2n-N-1}\dfrac{(1+pq^{n+N})(1-q^n)}{(1+pq^{2n-1})(1+pq^{2n})}.\end{cases}$$

$q$-**Difference equation.**

$$q^{-n}(1-q^n)(1+pq^n)y(x) = (1-q^{x-N})y(x+1) +$$
$$-\left[(1-q^{x-N})-p(1-q^x)\right]y(x) - p(1-q^x)y(x-1), \qquad (3.15.4)$$

where

$$y(x) = K_n(q^{-x};p,N;q).$$

**Generating functions.**

$${}_2\phi_0\left(\begin{array}{c}q^{x-N},0\\-\end{array}\bigg|q;-q^{-x}t\right){}_1\phi_1\left(\begin{array}{c}q^{-x}\\0\end{array}\bigg|q;pqt\right) \sim \sum_{n=0}^{N}\frac{(q^{-N};q)_n}{(q;q)_n}q^{-\binom{n}{2}}K_n(q^{-x};p,N;q)t^n. \qquad (3.15.5)$$

$${}_2\tilde{\phi}_1\left(\begin{array}{c}q^{x-N},0\\q^{-N}\end{array}\bigg|q;q^{-x}t\right){}_0\phi_1\left(\begin{array}{c}-\\-pq^{N+1}\end{array}\bigg|q;-pq^{N+1-x}t\right) \simeq \sum_{n=0}^{N}\frac{K_n(q^{-x};p,N;q)}{(-pq^{N+1},q;q)_n}t^n. \qquad (3.15.6)$$

**Remarks.** The $q$-Krawtchouk polynomials defined by (3.15.1) and the dual $q$-Krawtchouk polynomials given by (3.17.1) are related in the following way :

$$K_n(q^{-x};p,N;q) = K_x(\lambda(n);-pq^N,N|q)$$

with

$$\lambda(n) = q^{-n} - pq^n$$

or
$$K_n(\lambda(x); c, N|q) = K_x(q^{-n}; -cq^{-N}, N; q)$$

with
$$\lambda(x) = q^{-x} + cq^{x-N}.$$

The generating function (3.15.5) must be seen as an equality in terms of formal power series. For $x = 0, 1, 2, \ldots, N$ this generating function can also be written as :

$$_2\phi_0\left(\begin{array}{c} q^{x-N}, 0 \\ - \end{array} \bigg| q; -q^{-x}t\right) {}_1\phi_1\left(\begin{array}{c} q^{-x} \\ 0 \end{array} \bigg| q; pqt\right) = \sum_{n=0}^{N} \frac{(q^{-N}; q)_n}{(q; q)_n} q^{-\binom{n}{2}} K_n(q^{-x}; p, N; q) t^n.$$

**References.** [43], [114], [180], [181], [216], [217].

## 3.16 Affine $q$-Krawtchouk

**Definition.**
$$K_n^{Aff}(q^{-x}; p, N; q) = {}_3\tilde{\phi}_2\left(\begin{array}{c} q^{-n}, 0, q^{-x} \\ pq, q^{-N} \end{array} \bigg| q; q\right) \quad (3.16.1)$$
$$= \frac{(-pq)^n q^{\binom{n}{2}}}{(pq; q)_n} {}_2\tilde{\phi}_1\left(\begin{array}{c} q^{-n}, q^{x-N} \\ q^{-N} \end{array} \bigg| q; \frac{q^{-x}}{p}\right), \; n = 0, 1, 2, \ldots, N.$$

**Orthogonality.**
$$\sum_{x=0}^{N} \frac{(pq; q)_x (q; q)_N}{(q; q)_x (q; q)_{N-x}} (pq)^{-x} K_m^{Aff}(q^{-x}; p, N; q) K_n^{Aff}(q^{-x}; p, N; q)$$
$$= (pq)^{n-N} \frac{(q; q)_n (q; q)_{N-n}}{(pq; q)_n (q; q)_N} \delta_{mn}, \; 0 < pq < 1. \quad (3.16.2)$$

**Recurrence relation.**
$$-(1 - q^{-x}) K_n^{Aff}(q^{-x}) = (1 - q^{n-N})(1 - pq^{n+1}) K_{n+1}^{Aff}(q^{-x}) +$$
$$- \left[(1 - q^{n-N})(1 - pq^{n+1}) - pq^{n-N}(1 - q^n)\right] K_n^{Aff}(q^{-x}) - pq^{n-N}(1 - q^n) K_{n-1}^{Aff}(q^{-x}), (3.16.3)$$

where
$$K_n^{Aff}(q^{-x}) := K_n^{Aff}(q^{-x}; p, N; q).$$

**$q$-Difference equation.**
$$q^{-n}(1 - q^n) y(x) = B(x) y(x + 1) - [B(x) + D(x)] y(x) + D(x) y(x - 1), \quad (3.16.4)$$

where
$$y(x) = K_n^{Aff}(q^{-x}; p, N; q)$$

and
$$\begin{cases} B(x) = (1 - q^{x-N})(1 - pq^{x+1}) \\ D(x) = -p(1 - q^x) q^{x-N}. \end{cases}$$

**Generating functions.**
$$\frac{(q^{-N}t; q)_\infty}{(q^{-x}t; q)_\infty} {}_1\phi_1\left(\begin{array}{c} q^{-x} \\ pq \end{array} \bigg| q; pqt\right) \simeq \sum_{n=0}^{N} \frac{(q^{-N}; q)_n}{(q; q)_n} K_n^{Aff}(q^{-x}; p, N; q) t^n. \quad (3.16.5)$$

$$(pqt;q)_\infty \cdot {}_2\tilde{\phi}_1\left(\begin{array}{c}q^{x-N},0\\q^{-N}\end{array}\bigg|q;q^{-x}t\right) \simeq \sum_{n=0}^{N}\frac{(pq;q)_n}{(q;q)_n}K_n^{Aff}(q^{-x};p,N;q)t^n. \qquad (3.16.6)$$

**Remarks.** The affine $q$-Krawtchouk polynomials defined by (3.16.1) and the big $q$-Laguerre polynomials given by (3.11.1) are related in the following way :

$$K_n^{Aff}(q^{-x};p,N;q) = P_n(q^{-x};p,q^{-N-1};q).$$

The affine $q$-Krawtchouk polynomials are related to the quantum $q$-Krawtchouk polynomials defined by (3.14.1) by the transformation $q \leftrightarrow q^{-1}$ in the following way :

$$K_n^{Aff}(q^x;p,N;q^{-1}) = \frac{1}{(p^{-1}q;q)_n}K_n^{qtm}(q^{x-N};p^{-1},N;q).$$

For $x = 0, 1, 2, \ldots, N$ the generating function (3.16.5) can also be written as :

$$(q^{-N}t;q)_{N-x} \cdot {}_1\phi_1\left(\begin{array}{c}q^{-x}\\pq\end{array}\bigg|q;pqt\right) = \sum_{n=0}^{N}\frac{(q^{-N};q)_n}{(q;q)_n}K_n^{Aff}(q^{-x};p,N;q)t^n.$$

**References.** [85], [86], [88], [114], [217].

## 3.17 Dual $q$-Krawtchouk

**Definition.**

$$\begin{array}{rcl}K_n(\lambda(x);c,N|q) & = & {}_3\tilde{\phi}_2\left(\begin{array}{c}q^{-n},q^{-x},cq^{x-N}\\q^{-N},0\end{array}\bigg|q;q\right) \qquad (3.17.1)\\ & = & \dfrac{(q^{x-N};q)_n}{(q^{-N};q)_n q^{nx}}{}_2\phi_1\left(\begin{array}{c}q^{-n},q^{-x}\\q^{N-x-n+1}\end{array}\bigg|q;cq^{x+1}\right), \ n=0,1,2,\ldots,N,\end{array}$$

where

$$\lambda(x) := q^{-x} + cq^{x-N}.$$

**Orthogonality.**

$$\sum_{x=0}^{N}\frac{(cq^{-N},q^{-N};q)_x}{(q,cq;q)_x}\frac{(1-cq^{2x-N})}{(1-cq^{-N})}c^{-x}q^{x(2N-x)}K_m(\lambda(x))K_n(\lambda(x))$$
$$= (c^{-1};q)_N\frac{(q;q)_n}{(q^{-N};q)_n}(cq^{-N})^n\delta_{mn}, \qquad (3.17.2)$$

where

$$K_n(\lambda(x)) := K_n(\lambda(x);c,N|q).$$

**Recurrence relation.**

$$-(1-q^{-x})(1-cq^{x-N})K_n(\lambda(x)) = (1-q^{n-N})K_{n+1}(\lambda(x)) +$$
$$-\left[(1-q^{n-N})+cq^{-N}(1-q^n)\right]K_n(\lambda(x))+cq^{-N}(1-q^n)K_{n-1}(\lambda(x)), \qquad (3.17.3)$$

where

$$K_n(\lambda(x)) := K_n(\lambda(x);c,N|q).$$

**$q$-Difference equation.**

$$q^{-n}(1-q^n)y(x) = B(x)y(x+1) - [B(x)+D(x)]y(x) + D(x)y(x-1), \qquad (3.17.4)$$

where
$$y(x) = K_n(\lambda(x); c, N|q)$$

and
$$\begin{cases} B(x) = \dfrac{(1-q^{x-N})(1-cq^{x-N})}{(1-cq^{2x-N})(1-cq^{2x-N+1})} \\ \\ D(x) = cq^{2x-2N-1}\dfrac{(1-q^x)(1-cq^x)}{(1-cq^{2x-N-1})(1-cq^{2x-N})}. \end{cases}$$

**Generating functions.**

$$\frac{(q^{-N}t, cq^{-N}t; q)_\infty}{(q^{-x}t, cq^{x-N}t; q)_\infty} \simeq \sum_{n=0}^{N} \frac{(q^{-N};q)_n}{(q;q)_n} K_n(\lambda(x); c, N|q) t^n. \tag{3.17.5}$$

$$\frac{1}{(q^{-x}t;q)_\infty} {}_2\tilde{\phi}_1\left(\begin{array}{c} q^{-x}, c^{-1}q^{-x} \\ q^{-N} \end{array} \middle| q; cq^{x-N}t \right) \simeq \sum_{n=0}^{N} \frac{K_n(\lambda(x); c, N|q)}{(q;q)_n} t^n. \tag{3.17.6}$$

**Remark.** The dual $q$-Krawtchouk polynomials defined by (3.17.1) and the $q$-Krawtchouk polynomials given by (3.15.1) are related in the following way :

$$K_n(q^{-x}; p, N; q) = K_x(\lambda(n); -pq^N, N|q)$$

with
$$\lambda(n) = q^{-n} - pq^n$$

or
$$K_n(\lambda(x); c, N|q) = K_x(q^{-n}; -cq^{-N}, N; q)$$

with
$$\lambda(x) = q^{-x} + cq^{x-N}.$$

For $x = 0, 1, 2, \ldots, N$ the generating function (3.17.5) can also be written as :

$$(q^{-N}t; q)_{N-x} \cdot (cq^{-N}t; q)_x = \sum_{n=0}^{N} \frac{(q^{-N};q)_n}{(q;q)_n} K_n(\lambda(x); c, N|q) t^n.$$

**References.** [160], [163].

## 3.18 Continuous big $q$-Hermite

**Definition.**
$$\begin{aligned} H_n(x; a|q) &= a^{-n} {}_3\phi_2\left(\begin{array}{c} q^{-n}, ae^{i\theta}, ae^{-i\theta} \\ 0, 0 \end{array} \middle| q; q\right) \\ &= e^{in\theta} {}_2\phi_0\left(\begin{array}{c} q^{-n}, ae^{i\theta} \\ - \end{array} \middle| q; q^n e^{-2i\theta}\right), \quad x = \cos\theta. \end{aligned} \tag{3.18.1}$$

**Orthogonality.** When $a$ is real and $|a| < 1$, then we have the following orthogonality relation

$$\frac{1}{2\pi} \int_{-1}^{1} \frac{w(x)}{\sqrt{1-x^2}} H_m(x; a|q) H_n(x; a|q) dx = \frac{\delta_{mn}}{(q^{n+1}; q)_\infty}, \tag{3.18.2}$$

where
$$w(x) := w(x; a|q) = \left|\frac{(e^{2i\theta}; q)_\infty}{(ae^{i\theta}; q)_\infty}\right|^2 = \frac{h(x,1)h(x,-1)h(x,q^{\frac{1}{2}})h(x,-q^{\frac{1}{2}})}{h(x,a)},$$

with
$$h(x,\alpha) := \prod_{k=0}^{\infty} \left[1 - 2\alpha x q^k + \alpha^2 q^{2k}\right] = \left(\alpha e^{i\theta}, \alpha e^{-i\theta}; q\right)_{\infty}, \quad x = \cos\theta.$$

If $a > 1$, then we have another orthogonality relation given by :

$$\frac{1}{2\pi} \int_{-1}^{1} \frac{w(x)}{\sqrt{1-x^2}} H_m(x;a|q) H_n(x;a|q) dx +$$

$$+ \sum_{\substack{k \\ 1 < aq^k \le a}} w_k H_m(x_k;a|q) H_n(x_k;a|q) = \frac{\delta_{mn}}{(q^{n+1};q)_{\infty}}, \qquad (3.18.3)$$

where $w(x)$ is as before,
$$x_k = \frac{aq^k + (aq^k)^{-1}}{2}$$
and
$$w_k = \frac{(a^{-2};q)_{\infty}}{(q;q)_{\infty}} \frac{(1-a^2 q^{2k})(a^2;q)_k}{(1-a^2)(q;q)_k} q^{-\frac{3}{2}k^2 - \frac{1}{2}k} \left(-\frac{1}{a^4}\right)^k.$$

**Recurrence relation.**

$$2x H_n(x;a|q) = H_{n+1}(x;a|q) + aq^n H_n(x;a|q) + (1-q^n) H_{n-1}(x;a|q). \qquad (3.18.4)$$

**$q$-Difference equations.**

$$(1-q)^2 D_q \left[\tilde{w}(x;aq^{\frac{1}{2}}|q) D_q y(x)\right] + 4q^{-n+1}(1-q^n)\tilde{w}(x;a|q) y(x) = 0, \quad y(x) = H_n(x;a|q), \quad (3.18.5)$$

where
$$\tilde{w}(x;a|q) := \frac{w(x;a|q)}{\sqrt{1-x^2}}$$

and
$$D_q f(x) := \frac{\delta_q f(x)}{\delta_q x} \quad \text{with} \quad \delta_q f(e^{i\theta}) = f(q^{\frac{1}{2}} e^{i\theta}) - f(q^{-\frac{1}{2}} e^{i\theta}), \quad x = \cos\theta.$$

If we define
$$P_n(z) := a^{-n} {}_3\phi_2\left(\begin{matrix} q^{-n}, az, az^{-1} \\ 0, 0 \end{matrix} \bigg| q; q\right)$$

then the $q$-difference equation can also be written in the form

$$q^{-n}(1-q^n) P_n(z) = A(z) P_n(qz) - \left[A(z) + A(z^{-1})\right] P_n(z) + A(z^{-1}) P_n(q^{-1}z), \qquad (3.18.6)$$

where
$$A(z) = \frac{(1-az)}{(1-z^2)(1-qz^2)}.$$

**Generating function.**

$$\frac{(at;q)_{\infty}}{(e^{i\theta}t, e^{-i\theta}t; q)_{\infty}} = \sum_{n=0}^{\infty} \frac{H_n(x;a|q)}{(q;q)_n} t^n, \quad x = \cos\theta. \qquad (3.18.7)$$

**References.**

## 3.19 Continuous $q$-Laguerre

**Definitions.** We have two kinds of continuous $q$-Laguerre polynomials coming from the continuous $q$-Jacobi polynomials defined by (3.10.1) and (3.10.2) :

$$P_n^{(\alpha)}(x|q) = \frac{(q^{\alpha+1};q)_n}{(q;q)_n} {}_3\phi_2 \left( \begin{array}{c} q^{-n}, q^{\frac{1}{2}\alpha+\frac{1}{4}}e^{i\theta}, q^{\frac{1}{2}\alpha+\frac{1}{4}}e^{-i\theta} \\ q^{\alpha+1}, 0 \end{array} \bigg| q; q \right) \qquad (3.19.1)$$

$$= \frac{(q^{\frac{1}{2}\alpha+\frac{3}{4}}e^{-i\theta};q)_n}{(q;q)_n} q^{(\frac{1}{2}\alpha+\frac{1}{4})n} e^{in\theta} {}_2\phi_1 \left( \begin{array}{c} q^{-n}, q^{\frac{1}{2}\alpha+\frac{1}{4}}e^{i\theta} \\ q^{-\frac{1}{2}\alpha+\frac{1}{4}-n}e^{i\theta} \end{array} \bigg| q; q^{-\frac{1}{2}\alpha+\frac{1}{4}}e^{-i\theta} \right), \; x = \cos\theta$$

and

$$P_n^{(\alpha)}(x;q) = \frac{(q^{\alpha+1};q)_n}{(q;q)_n} {}_3\phi_2 \left( \begin{array}{c} q^{-n}, q^{\frac{1}{2}}e^{i\theta}, q^{\frac{1}{2}}e^{-i\theta} \\ q^{\alpha+1}, -q \end{array} \bigg| q; q \right), \; x = \cos\theta. \qquad (3.19.2)$$

These two $q$-analogues of the Laguerre polynomials are connected by the following quadratic transformation :

$$P_n^{(\alpha)}(x|q^2) = q^{n\alpha} P_n^{(\alpha)}(x;q).$$

**Orthogonality.** For $\alpha \geq -\frac{1}{2}$ the orthogonality relations are respectively

$$\frac{1}{2\pi} \int_{-1}^{1} \frac{w(x|q)}{\sqrt{1-x^2}} P_m^{(\alpha)}(x|q) P_n^{(\alpha)}(x|q) dx = \frac{1}{(q, q^{\alpha+1};q)_\infty} \frac{(q^{\alpha+1};q)_n}{(q;q)_n} q^{(\alpha+\frac{1}{2})n} \delta_{mn}, \qquad (3.19.3)$$

where

$$w(x|q) := w(x;q^\alpha|q) = \left| \frac{(e^{2i\theta};q)_\infty}{(q^{\frac{1}{2}\alpha+\frac{1}{4}}e^{i\theta}, q^{\frac{1}{2}\alpha+\frac{3}{4}}e^{i\theta};q)_\infty} \right|^2 = \left| \frac{(e^{i\theta}, -e^{i\theta};q^{\frac{1}{2}})_\infty}{(q^{\frac{1}{2}\alpha+\frac{1}{4}}e^{i\theta};q^{\frac{1}{2}})_\infty} \right|^2$$

$$= \frac{h(x,1)h(x,-1)h(x,q^{\frac{1}{2}})h(x,-q^{\frac{1}{2}})}{h(x,q^{\frac{1}{2}\alpha+\frac{1}{4}})h(x,q^{\frac{1}{2}\alpha+\frac{3}{4}})},$$

with

$$h(x,\alpha) := \prod_{k=0}^{\infty} \left[1 - 2\alpha x q^k + \alpha^2 q^{2k}\right] = \left(\alpha e^{i\theta}, \alpha e^{-i\theta};q\right)_\infty, \; x = \cos\theta$$

and

$$\frac{1}{2\pi} \int_{-1}^{1} \frac{w(x;q)}{\sqrt{1-x^2}} P_m^{(\alpha)}(x;q) P_n^{(\alpha)}(x;q) dx$$

$$= \frac{1}{(q,-q,q^{\alpha+1},-q^{\alpha+1};q)_\infty} \frac{(q^{\alpha+1},-q^{\alpha+1};q)_n}{(q,-q;q)_n} q^n \delta_{mn}, \qquad (3.19.4)$$

where

$$w(x;q) := w(x;q^\alpha;q) = \left| \frac{(e^{2i\theta};q)_\infty}{(q^{\alpha+\frac{1}{2}}e^{i\theta}, q^{\frac{1}{2}}e^{i\theta}, -q^{\frac{1}{2}}e^{i\theta};q)_\infty} \right|^2 = \left| \frac{(e^{i\theta}, -e^{i\theta};q)_\infty}{(q^{\alpha+\frac{1}{2}}e^{i\theta};q)_\infty} \right|^2 = \frac{h(x,1)h(x,-1)}{h(x,q^{\alpha+\frac{1}{2}})},$$

with

$$h(x,\alpha) := \prod_{k=0}^{\infty} \left[1 - 2\alpha x q^k + \alpha^2 q^{2k}\right] = \left(\alpha e^{i\theta}, \alpha e^{-i\theta};q\right)_\infty, \; x = \cos\theta.$$

**Recurrence relations.**

$$2xP_n^{(\alpha)}(x|q) = q^{-\frac{1}{2}\alpha-\frac{1}{4}}(1-q^{n+1})P_{n+1}^{(\alpha)}(x|q) +$$
$$+ q^{n+\frac{1}{2}\alpha+\frac{1}{4}}(1+q^{\frac{1}{2}})P_n^{(\alpha)}(x|q) + q^{\frac{1}{2}\alpha+\frac{1}{4}}(1-q^{n+\alpha})P_{n-1}^{(\alpha)}(x|q). \quad (3.19.5)$$

$$2xP_n^{(\alpha)}(x;q) = q^{-\frac{1}{2}}(1-q^{2n+2})P_{n+1}^{(\alpha)}(x;q) +$$
$$+ q^{2n+\alpha+\frac{1}{2}}(1+q)P_n^{(\alpha)}(x;q) + q^{\frac{1}{2}}(1-q^{2n+2\alpha})P_{n-1}^{(\alpha)}(x;q). \quad (3.19.6)$$

**$q$-Difference equations.**

$$(1-q)^2 D_q\left[\tilde{w}(x;q^{\alpha+\frac{1}{2}}|q)D_q y(x)\right] + 4q^{-n+1}(1-q^n)\tilde{w}(x;q^\alpha|q)y(x) = 0, \ y(x) = P_n^{(\alpha)}(x|q), \quad (3.19.7)$$

where

$$\tilde{w}(x;q^\alpha|q) := \frac{w(x;q^\alpha|q)}{\sqrt{1-x^2}}$$

and

$$D_q f(x) := \frac{\delta_q f(x)}{\delta_q x} \ \text{ with } \ \delta_q f(e^{i\theta}) = f(q^{\frac{1}{2}}e^{i\theta}) - f(q^{-\frac{1}{2}}e^{i\theta}), \ x = \cos\theta.$$

$$(1-q)^2 D_q\left[\tilde{w}(x;q^{\alpha+\frac{1}{2}};q)D_q y(x)\right] + 4q^{-n+1}(1-q^n)\tilde{w}(x;q^\alpha;q)y(x) = 0, \ y(x) = P_n^{(\alpha)}(x;q), \quad (3.19.8)$$

where

$$\tilde{w}(x;q^\alpha;q) := \frac{w(x;q^\alpha;q)}{\sqrt{1-x^2}}$$

and

$$D_q f(x) := \frac{\delta_q f(x)}{\delta_q x} \ \text{ with } \ \delta_q f(e^{i\theta}) = f(q^{\frac{1}{2}}e^{i\theta}) - f(q^{-\frac{1}{2}}e^{i\theta}), \ x = \cos\theta.$$

**Generating functions.**

$$\frac{(q^{\alpha+\frac{1}{2}}t, q^{\alpha+1}t;q)_\infty}{(q^{\frac{1}{2}\alpha+\frac{1}{4}}e^{i\theta}t, q^{\frac{1}{2}\alpha+\frac{1}{4}}e^{-i\theta}t;q)_\infty} = \sum_{n=0}^\infty P_n^{(\alpha)}(x|q)t^n, \ x = \cos\theta. \quad (3.19.9)$$

$$\frac{1}{(e^{i\theta}t;q)_\infty} {}_2\phi_1\left(\begin{array}{c} q^{\frac{1}{2}\alpha+\frac{1}{4}}e^{i\theta}, q^{\frac{1}{2}\alpha+\frac{3}{4}}e^{i\theta} \\ q^{\alpha+1} \end{array} \bigg| q; e^{-i\theta}t\right) = \sum_{n=0}^\infty \frac{P_n^{(\alpha)}(x|q)t^n}{(q^{\alpha+1};q)_n q^{(\frac{1}{2}\alpha+\frac{1}{4})n}}, \ x = \cos\theta. \quad (3.19.10)$$

$$\frac{(q^{\alpha+1}t;q)_\infty}{(q^{\frac{1}{2}}e^{i\theta}t;q)_\infty} {}_2\phi_1\left(\begin{array}{c} q^{\frac{1}{2}}e^{i\theta}, -q^{\frac{1}{2}}e^{i\theta} \\ -q \end{array} \bigg| q; q^{\frac{1}{2}}e^{-i\theta}t\right) = \sum_{n=0}^\infty P_n^{(\alpha)}(x;q)t^n, \ x = \cos\theta. \quad (3.19.11)$$

$$\frac{(-q^{\frac{1}{2}}t;q)_\infty}{(e^{i\theta}t;q)_\infty} {}_2\phi_1\left(\begin{array}{c} q^{\frac{1}{2}}e^{i\theta}, q^{\alpha+\frac{1}{2}}e^{i\theta} \\ q^{\alpha+1} \end{array} \bigg| q; e^{-i\theta}t\right) = \sum_{n=0}^\infty \frac{(-q;q)_n}{(q^{\alpha+1};q)_n} \frac{P_n^{(\alpha)}(x;q)}{q^{\frac{1}{2}n}} t^n, \ x = \cos\theta. \quad (3.19.12)$$

$$\frac{(q^{\frac{1}{2}}t;q)_\infty}{(e^{-i\theta}t;q)_\infty} {}_2\phi_1\left(\begin{array}{c} -q^{\frac{1}{2}}e^{-i\theta}, q^{\alpha+\frac{1}{2}}e^{-i\theta} \\ -q^{\alpha+1} \end{array} \bigg| q; e^{i\theta}t\right)$$
$$= \sum_{n=0}^\infty \frac{(-q;q)_n}{(-q^{\alpha+1};q)_n} \frac{P_n^{(\alpha)}(x;q)}{q^{\frac{1}{2}n}} t^n, \ x = \cos\theta. \quad (3.19.13)$$

**References.** [45].

## 3.20 Little $q$-Laguerre / Wall

**Definition.**

$$p_n(x;a|q) = {}_2\phi_1\left(\begin{array}{c}q^{-n},0\\aq\end{array}\bigg| q;qx\right) \qquad (3.20.1)$$

$$= \frac{1}{(a^{-1}q^{-n};q)_n}{}_2\phi_0\left(\begin{array}{c}q^{-n},x^{-1}\\-\end{array}\bigg| q;\frac{x}{a}\right).$$

**Orthogonality.**

$$\sum_{k=0}^{\infty}\frac{(aq)^k}{(q;q)_k}p_m(q^k;a|q)p_n(q^k;a|q) = \frac{(aq)^n}{(aq;q)_\infty}\frac{(q;q)_n}{(aq;q)_n}\delta_{mn},\ 0<aq<1. \qquad (3.20.2)$$

**Recurrence relation.**

$$-xp_n(x;a|q) = A_n p_{n+1}(x;a|q) - (A_n+C_n)p_n(x;a|q) + C_n p_{n-1}(x;a|q), \qquad (3.20.3)$$

where

$$\begin{cases} A_n = q^n(1-aq^{n+1}) \\ C_n = aq^n(1-q^n). \end{cases}$$

**$q$-Difference equation.**

$$-q^{-n}(1-q^n)xy(x) = ay(qx) + (x-a-1)y(x) + (1-x)y(q^{-1}x),\ y(x) = p_n(x;a|q). \qquad (3.20.4)$$

**Generating function.**

$$(t;q)_\infty \cdot {}_2\phi_1\left(\begin{array}{c}0,0\\aq\end{array}\bigg| q;xt\right) = \sum_{n=0}^{\infty}\frac{(-1)^n q^{\binom{n}{2}}}{(q;q)_n}p_n(x;a|q)t^n. \qquad (3.20.5)$$

**Remark.** If we set $a = q^\alpha$ and change $q$ to $q^{-1}$ we find the $q$-Laguerre polynomials defined by (3.21.1) in the following way :

$$p_n(x;q^{-\alpha}|q^{-1}) = \frac{(q;q)_n}{(q^{\alpha+1};q)_n}L_n^{(\alpha)}(-x;q).$$

**References.** [10], [21], [75], [77], [78], [114], [160], [161], [180], [221], [224].

## 3.21 $q$-Laguerre

**Definition.**

$$L_n^{(\alpha)}(x;q) = \frac{(q^{\alpha+1};q)_n}{(q;q)_n}{}_1\phi_1\left(\begin{array}{c}q^{-n}\\q^{\alpha+1}\end{array}\bigg| q;-xq^{n+\alpha+1}\right) \qquad (3.21.1)$$

$$= \frac{1}{(q;q)_n}{}_2\phi_1\left(\begin{array}{c}q^{-n},-x\\0\end{array}\bigg| q;q^{n+\alpha+1}\right).$$

**Orthogonality.** The $q$-Laguerre polynomials satisfy two kinds of orthogonality relations, an absolutely continuous one and a discrete one. These orthogonality relations are given by, respectively :

$$\int_0^\infty \frac{x^\alpha}{(-x;q)_\infty}L_m^{(\alpha)}(x;q)L_n^{(\alpha)}(x;q)dx = \frac{(q^{-\alpha};q)_\infty}{(q;q)_\infty}\frac{(q^{\alpha+1};q)_n}{(q;q)_n q^n}\Gamma(-\alpha)\Gamma(\alpha+1)\delta_{mn},\ \alpha>-1 \qquad (3.21.2)$$

and

$$\sum_{k=-\infty}^{\infty} \frac{q^{k\alpha+k}}{(-cq^k;q)_\infty} L_m^{(\alpha)}(cq^k;q) L_n^{(\alpha)}(cq^k;q)$$
$$= \frac{(q,-cq^{\alpha+1},-c^{-1}q^{-\alpha};q)_\infty}{(q^{\alpha+1},-c,-c^{-1}q;q)_\infty} \frac{(q^{\alpha+1};q)_n}{(q;q)_n q^n} \delta_{mn}, \ \alpha > -1, \ c > 0. \tag{3.21.3}$$

**Recurrence relation.**

$$-q^{2n+\alpha+1} x L_n^{(\alpha)}(x;q) = (1-q^{n+1}) L_{n+1}^{(\alpha)}(x;q) +$$
$$- \left[(1-q^{n+1}) + q(1-q^{n+\alpha})\right] L_n^{(\alpha)}(x;q) + q(1-q^{n+\alpha}) L_{n-1}^{(\alpha)}(x;q). \tag{3.21.4}$$

**$q$-Difference equation.**

$$-q^\alpha (1-q^n) x y(x) = q^\alpha (1+x) y(qx) - [1 + q^\alpha (1+x)] y(x) + y(q^{-1}x), \tag{3.21.5}$$

where

$$y(x) = L_n^{(\alpha)}(x;q).$$

**Generating functions.**

$$\frac{(-xt;q)_\infty}{(t;q)_\infty} {}_2\phi_1 \left( \begin{array}{c} 0,0 \\ q^{\alpha+1} \end{array} \Big| q; -xt \right) = \sum_{n=0}^\infty \frac{L_n^{(\alpha)}(x;q)}{(q^{\alpha+1};q)_n} t^n. \tag{3.21.6}$$

$$\frac{1}{(t;q)_\infty} {}_1\phi_1 \left( \begin{array}{c} -x \\ 0 \end{array} \Big| q; q^{\alpha+1}t \right) = \sum_{n=0}^\infty L_n^{(\alpha)}(x;q) t^n. \tag{3.21.7}$$

**Remarks.** The $q$-Laguerre polynomials are sometimes called the generalized Stieltjes-Wigert polynomials.

If we change $q$ to $q^{-1}$ we obtain the little $q$-Laguerre (or Wall) polynomials given by (3.20.1) in the following way :

$$L_n^{(\alpha)}(x;q^{-1}) = \frac{(q^{\alpha+1};q)_n}{(q;q)_n q^{n\alpha}} p_n(-x;q^\alpha|q).$$

The $q$-Laguerre polynomials defined by (3.21.1) and the alternative $q$-Charlier polynomials given by (3.22.1) are related in the following way :

$$\frac{K_n(q^x;a;q)}{(q;q)_n} = L_n^{(x-n)}(aq^n;q).$$

The $q$-Laguerre polynomials defined by (3.21.1) and the $q$-Charlier polynomials given by (3.23.1) are related in the following way :

$$\frac{C_n(-x;-q^{-\alpha};q)}{(q;q)_n} = L_n^{(\alpha)}(x;q).$$

Since the Stieltjes and Hamburger moment problems corresponding to the $q$-Laguerre polynomials are indeterminate there exist many different weight functions.

**References.** [8], [10], [30], [31], [45], [75], [77], [78], [96], [114], [177].

## 3.22 Alternative $q$-Charlier

**Definition.**

$$K_n(x;a;q) = {}_2\phi_1\left(\begin{array}{c}q^{-n},-aq^n\\0\end{array}\bigg|\,q;qx\right) \tag{3.22.1}$$

$$= (xq^{1-n};q)_n \cdot {}_1\phi_1\left(\begin{array}{c}q^{-n}\\xq^{1-n}\end{array}\bigg|\,q;-axq^{n+1}\right)$$

$$= (-axq^n)^n \cdot {}_2\phi_1\left(\begin{array}{c}q^{-n},x^{-1}\\0\end{array}\bigg|\,q;-\frac{q^{1-n}}{a}\right).$$

**Orthogonality.**

$$\sum_{k=0}^{\infty}\frac{a^k}{(q;q)_k}q^{\binom{k+1}{2}}K_m(q^k;a;q)K_n(q^k;a;q) = (q;q)_n(-aq^n;q)_\infty\frac{a^n q^{\binom{n+1}{2}}}{(1+aq^{2n})}\delta_{mn},\ a>0. \tag{3.22.2}$$

**Recurrence relation.**

$$-xK_n(x;a;q) = A_n K_{n+1}(x;a;q) - (A_n+C_n)K_n(x;a;q) + C_n K_{n-1}(x;a;q), \tag{3.22.3}$$

where

$$\begin{cases} A_n = q^n\dfrac{(1+aq^n)}{(1+aq^{2n})(1+aq^{2n+1})} \\[2mm] C_n = aq^{2n-1}\dfrac{(1-q^n)}{(1+aq^{2n-1})(1+aq^{2n})}. \end{cases}$$

**$q$-Difference equation.**

$$-q^{-n}(1-q^n)(1+aq^n)xy(x) = axy(qx) - (ax+1-x)y(x) + (1-x)y(q^{-1}x), \tag{3.22.4}$$

where

$$y(x) = K_n(x;a;q).$$

**Generating functions.**

$${}_2\phi_0\left(\begin{array}{c}x^{-1},0\\-\end{array}\bigg|\,q;xt\right){}_0\phi_1\left(\begin{array}{c}-\\0\end{array}\bigg|\,q;-aqxt\right) \sim \sum_{n=0}^{\infty}\frac{K_n(x;a;q)}{(q;q)_n}t^n. \tag{3.22.5}$$

$$\frac{(t;q)_\infty}{(xt;q)_\infty}{}_1\phi_3\left(\begin{array}{c}xt\\0,0,t\end{array}\bigg|\,q;-aqxt\right) = \sum_{n=0}^{\infty}\frac{(-1)^n q^{\binom{n}{2}}}{(q;q)_n}K_n(x;a;q)t^n. \tag{3.22.6}$$

**Remarks.** The alternative $q$-Charlier polynomials defined by (3.22.1) and the $q$-Laguerre polynomials given by (3.21.1) related in the following way :

$$\frac{K_n(q^x;a;q)}{(q;q)_n} = L_n^{(x-n)}(aq^n;q).$$

The generating function (3.22.5) must be seen as an equality in terms of formal power series. For $x=0,1,2,\ldots,N$ this generating function can also be written as :

$${}_2\phi_0\left(\begin{array}{c}q^{-x},0\\-\end{array}\bigg|\,q;q^x t\right){}_0\phi_1\left(\begin{array}{c}-\\0\end{array}\bigg|\,q;-aq^{x+1}t\right) = \sum_{n=0}^{\infty}\frac{K_n(q^x;a;q)}{(q;q)_n}t^n.$$

**References.**

## 3.23  $q$-Charlier

**Definition.**
$$C_n(q^{-x};a;q) = {}_2\phi_1\left(\begin{array}{c}q^{-n},q^{-x}\\0\end{array}\bigg|\, q;-\frac{q^{n+1}}{a}\right) \qquad (3.23.1)$$
$$= (-a^{-1}q;q)_n \cdot {}_1\phi_1\left(\begin{array}{c}q^{-n}\\-a^{-1}q\end{array}\bigg|\, q;-\frac{q^{n+1-x}}{a}\right).$$

**Orthogonality.**

$$\sum_{x=0}^{\infty}\frac{a^x}{(q;q)_x}q^{\binom{x}{2}}C_m(q^{-x};a;q)C_n(q^{-x};a;q)=q^{-n}(-a;q)_\infty(-a^{-1}q,q;q)_n\delta_{mn},\ a>0. \qquad (3.23.2)$$

**Recurrence relation.**

$$q^{2n+1}(1-q^{-x})C_n(q^{-x}) = aC_{n+1}(q^{-x}) +$$
$$-\left[a+q(1-q^n)(a+q^n)\right]C_n(q^{-x})+q(1-q^n)(a+q^n)C_{n-1}(q^{-x}), \quad (3.23.3)$$

where
$$C_n(q^{-x}) := C_n(q^{-x};a;q).$$

**$q$-Difference equation.**

$$q^n y(x) = aq^x y(x+1) - q^x(a-1)y(x) + (1-q^x)y(x-1),\ y(x)=C_n(q^{-x};a;q). \qquad (3.23.4)$$

**Generating functions.**

$$\frac{1}{(t;q)_\infty}{}_1\phi_1\left(\begin{array}{c}q^{-x}\\0\end{array}\bigg|\, q;-a^{-1}qt\right) = \sum_{n=0}^{\infty}\frac{C_n(q^{-x};a;q)}{(q;q)_n}t^n. \qquad (3.23.5)$$

$$\frac{(q^{-x}t;q)_\infty}{(t;q)_\infty}{}_2\phi_1\left(\begin{array}{c}0,0\\-a^{-1}q\end{array}\bigg|\, q;q^{-x}t\right)=\sum_{n=0}^{\infty}\frac{C_n(q^{-x};a;q)}{(-a^{-1}q,q;q)_n}t^n. \qquad (3.23.6)$$

**Remark.** The $q$-Charlier polynomials defined by (3.23.1) and the $q$-Laguerre polynomials given by (3.21.1) are related in the following way :

$$\frac{C_n(-x;-q^{-\alpha};q)}{(q;q)_n}=L_n^{(\alpha)}(x;q).$$

**References.** [114], [121], [180].

## 3.24  Al-Salam-Carlitz I

**Definition.**
$$U_n^{(a)}(x;q) = (-a)^n q^{\binom{n}{2}}{}_2\phi_1\left(\begin{array}{c}q^{-n},x^{-1}\\0\end{array}\bigg|\, q;\frac{qx}{a}\right). \qquad (3.24.1)$$

**Orthogonality.**

$$\int_a^1 (qx,a^{-1}qx;q)_\infty U_m^{(a)}(x;q)U_n^{(a)}(x;q)d_qx$$
$$= (-a)^n(1-q)(q;q)_n(q,a,a^{-1}q;q)_\infty q^{\binom{n}{2}}\delta_{mn},\ a<0. \qquad (3.24.2)$$

**Recurrence relation.**

$$xU_n^{(a)}(x;q) = U_{n+1}^{(a)}(x;q) + (a+1)q^n U_n^{(a)}(x;q) - aq^{n-1}(1-q^n)U_{n-1}^{(a)}(x;q). \qquad (3.24.3)$$

**$q$-Difference equation.**

$$(1-q^n)x^2 y(x) = aq^{n-1}y(qx) - \left[aq^{n-1} + q^n(1-x)(a-x)\right]y(x) + \\ + q^n(1-x)(a-x)y(q^{-1}x),\ y(x) = U_n^{(a)}(x;q). \qquad (3.24.4)$$

**Generating function.**

$$\frac{(t;q)_\infty (at;q)_\infty}{(xt;q)_\infty} = \sum_{n=0}^{\infty} \frac{U_n^{(a)}(x;q)}{(q;q)_n} t^n. \qquad (3.24.5)$$

**Remark.** The Al-Salam-Carlitz I polynomials are related to the Al-Salam-Carlitz II polynomials defined by (3.25.1) in the following way :

$$U_n^{(a)}(x;q^{-1}) = V_n^{(a)}(x;q).$$

**References.** [10], [13], [15], [75], [77], [84], [114], [125], [135].

## 3.25 Al-Salam-Carlitz II

**Definition.**

$$V_n^{(a)}(x;q) = (-a)^n q^{-\binom{n}{2}} {}_2\phi_0\left(\begin{array}{c} q^{-n}, x \\ - \end{array} \bigg| q; \frac{q^n}{a}\right). \qquad (3.25.1)$$

**Orthogonality.**

$$\sum_{k=0}^{\infty} \frac{q^{k^2} a^k}{(q;q)_k (aq;q)_k} V_m^{(a)}(q^{-k};q) V_n^{(a)}(q^{-k};q) = \frac{(q;q)_n a^n}{(aq;q)_\infty q^{n^2}} \delta_{mn},\ a > 0. \qquad (3.25.2)$$

**Recurrence relation.**

$$xV_n^{(a)}(x;q) = V_{n+1}^{(a)}(x;q) + (a+1)q^{-n} V_n^{(a)}(x;q) + aq^{-2n+1}(1-q^n)V_{n-1}^{(a)}(x;q). \qquad (3.25.3)$$

**$q$-Difference equation.**

$$-(1-q^n)x^2 y(x) = (1-x)(a-x)y(qx) - \left[(1-x)(a-x) + aq\right]y(x) + \\ + aqy(q^{-1}x),\ y(x) = V_n^{(a)}(x;q). \qquad (3.25.4)$$

**Generating functions.**

$$\frac{(xt;q)_\infty}{(t;q)_\infty (at;q)_\infty} = \sum_{n=0}^{\infty} \frac{(-1)^n q^{\binom{n}{2}}}{(q;q)_n} V_n^{(a)}(x;q) t^n. \qquad (3.25.5)$$

$$(at;q)_\infty \cdot {}_1\phi_1\left(\begin{array}{c} x \\ at \end{array} \bigg| q; t\right) = \sum_{n=0}^{\infty} \frac{q^{2\binom{n}{2}}}{(q;q)_n} V_n^{(a)}(x;q) t^n. \qquad (3.25.6)$$

**Remark.** The Al-Salam-Carlitz II polynomials are related to the Al-Salam-Carlitz I polynomials defined by (3.24.1) in the following way :

$$V_n^{(a)}(x;q^{-1}) = U_n^{(a)}(x;q).$$

**References.** [10], [13], [15], [74], [75], [77], [84], [125].

## 3.26 Continuous $q$-Hermite

**Definition.**
$$H_n(x|q) = e^{in\theta} {}_2\phi_0 \left( \begin{matrix} q^{-n}, 0 \\ - \end{matrix} \bigg| q; q^n e^{-2i\theta} \right), \ x = \cos\theta. \tag{3.26.1}$$

**Orthogonality.**
$$\frac{1}{2\pi} \int_{-1}^{1} \frac{w(x)}{\sqrt{1-x^2}} H_m(x|q) H_n(x|q) dx = \frac{\delta_{mn}}{(q^{n+1};q)_\infty}, \tag{3.26.2}$$

where
$$w(x) = \left| \left(e^{2i\theta}; q\right)_\infty \right|^2 = h(x,1)h(x,-1)h(x,q^{\frac{1}{2}})h(x,-q^{\frac{1}{2}}),$$

with
$$h(x,\alpha) := \prod_{k=0}^{\infty} \left[1 - 2\alpha x q^k + \alpha^2 q^{2k}\right] = \left(\alpha e^{i\theta}, \alpha e^{-i\theta}; q\right)_\infty, \ x = \cos\theta.$$

**Recurrence relation.**
$$2xH_n(x|q) = H_{n+1}(x|q) + (1-q^n)H_{n-1}(x|q). \tag{3.26.3}$$

**$q$-Difference equation.**
$$(1-q)^2 D_q \left[\tilde{w}(x) D_q y(x)\right] + 4q^{-n+1}(1-q^n)\tilde{w}(x)y(x) = 0, \ y(x) = H_n(x|q), \tag{3.26.4}$$

where
$$\tilde{w}(x) := \frac{w(x)}{\sqrt{1-x^2}}$$

and
$$D_q f(x) := \frac{\delta_q f(x)}{\delta_q x} \ \text{with} \ \delta_q f(e^{i\theta}) = f(q^{\frac{1}{2}} e^{i\theta}) - f(q^{-\frac{1}{2}} e^{i\theta}), \ x = \cos\theta.$$

**Generating functions.**
$$\frac{1}{\left|(e^{i\theta} t; q)_\infty\right|^2} = \sum_{n=0}^{\infty} \frac{H_n(x|q)}{(q;q)_n} t^n, \ x = \cos\theta. \tag{3.26.5}$$

$$(e^{i\theta} t; q)_\infty \cdot {}_1\phi_1 \left( \begin{matrix} 0 \\ e^{i\theta} t \end{matrix} \bigg| q; e^{-i\theta} t \right) = \sum_{n=0}^{\infty} \frac{(-1)^n q^{\binom{n}{2}}}{(q;q)_n} H_n(x|q) t^n, \ x = \cos\theta. \tag{3.26.6}$$

**Remark.** The continuous $q$-Hermite polynomials can also be written as :
$$H_n(x|q) = \sum_{k=0}^{n} \frac{(q;q)_n}{(q;q)_k (q;q)_{n-k}} e^{i(n-2k)\theta}, \ x = \cos\theta.$$

**References.** [6], [10], [20], [25], [31], [32], [37], [38], [45], [59], [62], [110], [114], [128], [137], [138], [180], [207], [208], [209].

## 3.27 Stieltjes-Wigert

**Definition.**
$$S_n(x;q) = \frac{1}{(q;q)_n} {}_1\phi_1\left(\begin{array}{c} q^{-n} \\ 0 \end{array} \bigg| q; -xq^{n+1}\right). \qquad (3.27.1)$$

**Orthogonality.**
$$\int_0^\infty \frac{S_m(x;q)S_n(x;q)}{(-x;q)_\infty(-qx^{-1};q)_\infty} dx = -\frac{\ln q}{q^n} \frac{(q;q)_\infty}{(q;q)_n} \delta_{mn}. \qquad (3.27.2)$$

**Recurrence relation.**
$$-q^{2n+1}xS_n(x;q) = (1-q^{n+1})S_{n+1}(x;q) - [1+q-q^{n+1}]S_n(x;q) + qS_{n-1}(x;q). \qquad (3.27.3)$$

**$q$-Difference equation.**
$$-x(1-q^n)y(x) = xy(qx) - (x+1)y(x) + y(q^{-1}x), \ y(x) = S_n(x;q). \qquad (3.27.4)$$

**Generating functions.**
$$\frac{1}{(t;q)_\infty} {}_0\phi_1\left(\begin{array}{c} - \\ 0 \end{array} \bigg| q; -qxt\right) = \sum_{n=0}^\infty S_n(x;q)t^n. \qquad (3.27.5)$$

$$(t;q)_\infty \cdot {}_0\phi_2\left(\begin{array}{c} - \\ 0,t \end{array} \bigg| q; -qxt\right) = \sum_{n=0}^\infty (-1)^n q^{\binom{n}{2}} S_n(x;q)t^n. \qquad (3.27.6)$$

**Remark.** Since the Stieltjes and Hamburger moment problems corresponding to the Stieltjes-Wigert polynomials are indeterminate there exist many different weight functions. For instance, they are also orthogonal with respect to the weight function
$$w(x) = \frac{\gamma}{\sqrt{\pi}} \exp\left(-\gamma^2 \ln^2 x\right), \ x > 0, \text{ with } \gamma^2 = -\frac{1}{2\ln q}.$$

**References.** [30], [31], [76], [77], [84], [180], [219], [220], [226].

## 3.28 Discrete $q$-Hermite I

**Definition.** The discrete $q$-Hermite I polynomials are Al-Salam-Carlitz I polynomials with $a = -1$ :
$$\begin{aligned} h_n(x;q) = U_n^{(-1)}(x;q) &= q^{\binom{n}{2}} {}_2\phi_1\left(\begin{array}{c} q^{-n}, x^{-1} \\ 0 \end{array} \bigg| q; -qx\right) \qquad (3.28.1) \\ &= x^n {}_2\phi_0\left(\begin{array}{c} q^{-n}, q^{-n+1} \\ - \end{array} \bigg| q^2; \frac{q^{2n-1}}{x^2}\right). \end{aligned}$$

**Orthogonality.**
$$\int_{-1}^1 (qx, -qx; q)_\infty h_m(x;q)h_n(x;q) d_q x = (1-q)(q;q)_n (q,-1,-q;q)_\infty q^{\binom{n}{2}} \delta_{mn}. \qquad (3.28.2)$$

**Recurrence relation.**
$$xh_n(x;q) = h_{n+1}(x;q) + q^{n-1}(1-q^n)h_{n-1}(x;q). \qquad (3.28.3)$$

**$q$-Difference equation.**

$$-q^{-n+1}x^2 y(x) = y(qx) - (1+q)y(x) + q(1-x^2)y(q^{-1}x), \ y(x) = h_n(x;q). \tag{3.28.4}$$

**Generating function.**

$$\frac{(t;q)_\infty (-t;q)_\infty}{(xt;q)_\infty} = \sum_{n=0}^{\infty} \frac{h_n(x;q)}{(q;q)_n} t^n. \tag{3.28.5}$$

**Remark.** The discrete $q$-Hermite I polynomials are related to the discrete $q$-Hermite II polynomials defined by (3.29.1) in the following way :

$$h_n(ix;q^{-1}) = i^n \tilde{h}_n(x;q).$$

**References.** [10], [13], [62], [114], [121].

## 3.29 Discrete $q$-Hermite II

**Definition.** The discrete $q$-Hermite II polynomials are Al-Salam-Carlitz II polynomials with $a = -1$ :

$$\tilde{h}_n(x;q) = i^{-n} V_n^{(-1)}(ix;q) = i^{-n} q^{-\binom{n}{2}} {}_2\phi_0 \left( \begin{array}{c} q^{-n}, ix \\ - \end{array} \bigg| q; -q^n \right) \tag{3.29.1}$$

$$= x^n {}_2\phi_1 \left( \begin{array}{c} q^{-n}, q^{-n+1} \\ 0 \end{array} \bigg| q^2; -\frac{q^2}{x^2} \right).$$

**Orthogonality.**

$$\sum_{k=-\infty}^{\infty} \left[ \tilde{h}_m(cq^k;q)\tilde{h}_n(cq^k;q) + \tilde{h}_m(-cq^k;q)\tilde{h}_n(-cq^k;q) \right] w(cq^k) q^k$$

$$= 2 \frac{(q^2, -c^2 q, -c^{-2} q; q^2)_\infty}{(q, -c^2, -c^{-2} q^2; q^2)_\infty} \frac{(q;q)_n}{q^{n^2}} \delta_{mn}, \ c > 0, \tag{3.29.2}$$

where

$$w(x) = \frac{1}{(ix;q)_\infty (-ix;q)_\infty} = \frac{1}{(-x^2;q^2)_\infty}.$$

**Recurrence relation.**

$$x\tilde{h}_n(x;q) = \tilde{h}_{n+1}(x;q) + q^{-2n+1}(1-q^n)\tilde{h}_{n-1}(x;q). \tag{3.29.3}$$

**$q$-Difference equation.**

$$-(1-q^n)x^2 \tilde{h}_n(x;q) = (1+x^2)\tilde{h}_n(qx;q) - (1+x^2+q)\tilde{h}_n(x;q) + q\tilde{h}_n(q^{-1}x;q). \tag{3.29.4}$$

**Generating functions.**

$$\frac{(-xt;q)_\infty}{(-t^2;q^2)_\infty} = \sum_{n=0}^{\infty} \frac{q^{\binom{n}{2}}}{(q;q)_n} \tilde{h}_n(x;q) t^n. \tag{3.29.5}$$

$$(-it;q)_\infty \cdot {}_1\phi_1 \left( \begin{array}{c} ix \\ -it \end{array} \bigg| q; it \right) = \sum_{n=0}^{\infty} \frac{(-1)^n q^{2\binom{n}{2}}}{(q;q)_n} \tilde{h}_n(x;q) t^n. \tag{3.29.6}$$

**Remark.** The discrete $q$-Hermite II polynomials are related to the discrete $q$-Hermite I polynomials defined by (3.28.1) in the following way :

$$\tilde{h}_n(x;q^{-1}) = i^{-n} h_n(ix;q).$$

**References.**

# Chapter 4

# Limit relations between basic hypergeometric orthogonal polynomials

## 4.1 Askey-Wilson → Continuous dual $q$-Hahn

The continuous dual $q$-Hahn polynomials defined by (3.3.1) simply follow from the Askey-Wilson polynomials given by (3.1.1) by setting $d = 0$ in (3.1.1) :

$$p_n(x; a, b, c, 0|q) = p_n(x; a, b, c|q).$$

## 4.2 Askey-Wilson → Continuous $q$-Hahn

The continuous $q$-Hahn polynomials defined by (3.4.1) can be obtained from the Askey-Wilson polynomials given by (3.1.1) by the substitutions $\theta \to \theta + \phi$, $a \to ae^{i\phi}$, $b \to be^{i\phi}$, $c \to ce^{-i\phi}$ and $d \to de^{-i\phi}$ :

$$p_n(\cos(\theta + \phi); ae^{i\phi}, be^{i\phi}, ce^{-i\phi}, de^{-i\phi}|q) = p_n(\cos(\theta + \phi); a, b, c, d; q).$$

## 4.3 Askey-Wilson → Big $q$-Jacobi

The big $q$-Jacobi polynomials defined by (3.5.1) can be obtained from the Askey-Wilson polynomials by setting $x \to \frac{1}{2}a^{-1}x$, $b = a^{-1}\alpha q$, $c = a^{-1}\gamma q$ and $d = a\beta\gamma^{-1}$ in

$$\tilde{p}_n(x; a, b, c, d|q) = \frac{a^n p_n(x; a, b, c, d|q)}{(ab, ac, ad; q)_n}$$

defined by (3.1.1) and then taking the limit $a \to 0$ :

$$\lim_{a \to 0} \tilde{p}_n\left(\frac{x}{2a}; a, \frac{\alpha q}{a}, \frac{\gamma q}{a}, \frac{a\beta}{\gamma}\bigg| q\right) = P_n(x; \alpha, \beta, \gamma; q).$$

## 4.4 Askey-Wilson → Continuous $q$-Jacobi

If we take $a = q^{\frac{1}{2}\alpha + \frac{1}{4}}$, $b = q^{\frac{1}{2}\alpha + \frac{3}{4}}$, $c = -q^{\frac{1}{2}\beta + \frac{1}{4}}$ and $d = -q^{\frac{1}{2}\beta + \frac{3}{4}}$ in the definition (3.1.1) of the Askey-Wilson polynomials and change the normalization we find the continuous $q$-Jacobi polynomials given by (3.10.1) :

$$\frac{q^{(\frac{1}{2}\alpha+\frac{1}{4})n}p_n\left(x;q^{\frac{1}{2}\alpha+\frac{1}{4}},q^{\frac{1}{2}\alpha+\frac{3}{4}},-q^{\frac{1}{2}\beta+\frac{1}{4}},-q^{\frac{1}{2}\beta+\frac{3}{4}}\Big|q\right)}{(q,-q^{\frac{1}{2}(\alpha+\beta+1)},-q^{\frac{1}{2}(\alpha+\beta+2)};q)_n}=P_n^{(\alpha,\beta)}(x|q).$$

In [196] M. Rahman takes $a=q^{\frac{1}{2}}$, $b=q^{\alpha+\frac{1}{2}}$, $c=-q^{\beta+\frac{1}{2}}$ and $d=-q^{\frac{1}{2}}$ to obtain after a change of normalization the continuous $q$-Jacobi polynomials defined by (3.10.2) :

$$\frac{q^{\frac{1}{2}n}p_n\left(x;q^{\frac{1}{2}},q^{\alpha+\frac{1}{2}},-q^{\beta+\frac{1}{2}},-q^{\frac{1}{2}}\Big|q\right)}{(q,-q,-q;q)_n}=P_n^{(\alpha,\beta)}(x;q).$$

As was pointed out in section 0.6 these two $q$-analogues of the Jacobi polynomials are not really different, since they are connected by the quadratic transformation

$$P_n^{(\alpha,\beta)}(x|q^2)=\frac{(-q;q)_n}{(-q^{\alpha+\beta+1};q)_n}q^{n\alpha}P_n^{(\alpha,\beta)}(x;q).$$

## 4.5  Askey-Wilson → Continuous $q$-ultraspherical / Rogers

If we set $a=\beta^{\frac{1}{2}}$, $b=\beta^{\frac{1}{2}}q^{\frac{1}{2}}$, $c=-\beta^{\frac{1}{2}}$ and $d=-\beta^{\frac{1}{2}}q^{\frac{1}{2}}$ in the definition (3.1.1) of the Askey-Wilson polynomials and change the normalization we obtain the continuous $q$-ultraspherical (or Rogers) polynomials defined by (3.10.15). In fact we have :

$$\frac{(\beta^2;q)_n p_n\left(x;\beta^{\frac{1}{2}},\beta^{\frac{1}{2}}q^{\frac{1}{2}},-\beta^{\frac{1}{2}},-\beta^{\frac{1}{2}}q^{\frac{1}{2}}\Big|q\right)}{(\beta q^{\frac{1}{2}},-\beta,-\beta q^{\frac{1}{2}},q;q)_n}=C_n(x;\beta|q).$$

## 4.6  $q$-Racah → Big $q$-Jacobi

The big $q$-Jacobi polynomials defined by (3.5.1) can be obtained from the $q$-Racah polynomials by setting $\delta=0$ in the definition (3.2.1) :

$$R_n(\mu(x);a,b,c,0|q)=P_n(q^{-x};a,b,c;q).$$

## 4.7  $q$-Racah → $q$-Hahn

The $q$-Hahn polynomials follow from the $q$-Racah polynomials by the substitution $\delta=0$ and $\gamma q=q^{-N}$ in the definition (3.2.1) of the $q$-Racah polynomials :

$$R_n(\mu(x);\alpha,\beta,q^{-N-1},0|q)=Q_n(q^{-x};\alpha,\beta,N|q).$$

Another way to obtain the $q$-Hahn polynomials from the $q$-Racah polynomials is by setting $\gamma=0$ and $\delta=\beta^{-1}q^{-N-1}$ in the definition (3.2.1) :

$$R_n(\mu(x);\alpha,\beta,0,\beta^{-1}q^{-N-1}|q)=Q_n(q^{-x};\alpha,\beta,N|q).$$

And if we take $\alpha q=q^{-N}$, $\beta\to\beta\gamma q^{N+1}$ and $\delta=0$ in the definition (3.2.1) of the $q$-Racah polynomials we find the $q$-Hahn polynomials given by (3.6.1) in the following way :

$$R_n(\mu(x);q^{-N-1},\beta\gamma q^{N+1},\gamma,0|q)=Q_n(q^{-x};\gamma,\beta,N|q).$$

Note that $\mu(x)=q^{-x}$ in each case.

## 4.8 $q$-Racah $\to$ Dual $q$-Hahn

To obtain the dual $q$-Hahn polynomials from the $q$-Racah polynomials we have to take $\beta = 0$ and $\alpha q = q^{-N}$ in (3.2.1) :

$$R_n(\mu(x); q^{-N-1}, 0, \gamma, \delta | q) = R_n(\mu(x); \gamma, \delta, N | q),$$

with

$$\mu(x) = q^{-x} + \gamma \delta q^{x+1}.$$

We may also take $\alpha = 0$ and $\beta = \delta^{-1} q^{-N-1}$ in (3.2.1) to obtain the dual $q$-Hahn polynomials from the $q$-Racah polynomials :

$$R_n(\mu(x); 0, \delta^{-1} q^{-N-1}, \gamma, \delta | q) = R_n(\mu(x); \gamma, \delta, N | q),$$

with

$$\mu(x) = q^{-x} + \gamma \delta q^{x+1}.$$

And if we take $\gamma q = q^{-N}$, $\delta \to \alpha \delta q^{N+1}$ and $\beta = 0$ in the definition (3.2.1) of the $q$-Racah polynomials we find the dual $q$-Hahn polynomials given by (3.7.1) in the following way :

$$R_n(\mu(x); \alpha, 0, q^{-N-1}, \alpha \delta q^{N+1} | q) = R_n(\mu(x); \alpha, \delta, N | q),$$

with

$$\mu(x) = q^{-x} + \alpha \delta q^{x+1}.$$

## 4.9 $q$-Racah $\to$ $q$-Krawtchouk

The $q$-Krawtchouk polynomials defined by (3.15.1) can be obtained from the $q$-Racah polynomials by setting $\alpha q = q^{-N}$, $\beta = -pq^N$ and $\gamma = \delta = 0$ in the definition (3.2.1) of the $q$-Racah polynomials :

$$R_n(q^{-x}; q^{-N-1}, -pq^N, 0, 0 | q) = K_n(q^{-x}; p, N; q).$$

Note that $\mu(x) = q^{-x}$ in this case.

## 4.10 $q$-Racah $\to$ Dual $q$-Krawtchouk

The dual $q$-Krawtchouk polynomials defined by (3.17.1) easily follow from the $q$-Racah polynomials given by (3.2.1) by using the substitutions $\alpha = \beta = 0$, $\gamma q = q^{-N}$ and $\delta = c$ :

$$R_n(\mu(x); 0, 0, q^{-N-1}, c | q) = K_n(\lambda(x); c, N | q).$$

Note that

$$\mu(x) = \lambda(x) = q^{-x} + c q^{x-N}.$$

## 4.11 Continuous dual $q$-Hahn $\to$ Al-Salam-Chihara

The Al-Salam-Chihara polynomials defined by (3.8.1) simply follow from the continuous dual $q$-Hahn polynomials by taking $c = 0$ in the definition (3.3.1) of the continuous dual $q$-Hahn polynomials :

$$p_n(x; a, b, 0 | q) = Q_n(x; a, b | q).$$

## 4.12 Continuous $q$-Hahn $\to$ $q$-Meixner-Pollaczek

The $q$-Meixner-Pollaczek polynomials defined by (3.9.1) simply follow from the continuous $q$-Hahn polynomials if we set $d = a$ and $b = c = 0$ in the definition (3.4.1) of the continuous $q$-Hahn polynomials :

$$\frac{p_n(\cos(\theta+\phi); a, 0, 0, a; q)}{(q;q)_n} = P_n(\cos(\theta+\phi); a|q).$$

## 4.13 Big $q$-Jacobi $\to$ Big $q$-Laguerre

If we set $b = 0$ in the definition (3.5.1) of the big $q$-Jacobi polynomials we obtain the big $q$-Laguerre polynomials given by (3.11.1) :

$$P_n(x; a, 0, c; q) = P_n(x; a, c; q).$$

## 4.14 Big $q$-Jacobi $\to$ Little $q$-Jacobi

The little $q$-Jacobi polynomials defined by (3.12.1) can be obtained from the big $q$-Jacobi polynomials by the substitution $x \to cqx$ in the definition (3.5.1) and then by the limit $c \to \infty$ :

$$\lim_{c \to \infty} P_n(cqx; a, b, c; q) = p_n(x; a, b|q).$$

## 4.15 Big $q$-Jacobi $\to$ $q$-Meixner

If we take the limit $a \to \infty$ in the definition (3.5.1) of the big $q$-Jacobi polynomials we simply obtain the $q$-Meixner polynomials defined by (3.13.1) :

$$\lim_{a \to \infty} P_n(q^{-x}; a, b, c; q) = M_n(q^{-x}; c, -b^{-1}; q).$$

## 4.16 $q$-Hahn $\to$ Little $q$-Jacobi

If we set $x \to N - x$ in the definition (3.6.1) of the $q$-Hahn polynomials and take the limit $N \to \infty$ we find the little $q$-Jacobi polynomials :

$$\lim_{N \to \infty} Q_n(q^{x-N}; \alpha, \beta, N|q) = p_n(q^x; \alpha, \beta|q),$$

where $p_n(q^x; \alpha, \beta|q)$ is defined by (3.12.1).

## 4.17 $q$-Hahn $\to$ $q$-Meixner

The $q$-Meixner polynomials defined by (3.13.1) can be obtained from the $q$-Hahn polynomials by setting $\alpha = b$ and $\beta = -b^{-1}c^{-1}q^{-N-1}$ in the definition (3.6.1) of the $q$-Hahn polynomials and letting $N \to \infty$ :

$$\lim_{N \to \infty} Q_n\left(q^{-x}; b, -b^{-1}c^{-1}q^{-N-1}, N|q\right) = M_n(q^{-x}; b, c; q).$$

## 4.18 $q$-Hahn $\to$ Quantum $q$-Krawtchouk

The quantum $q$-Krawtchouk polynomials defined by (3.14.1) simply follow from the $q$-Hahn polynomials by setting $\beta = p$ in the definition (3.6.1) of the $q$-Hahn polynomials and taking the limit $\alpha \to \infty$ :

$$\lim_{\alpha \to \infty} Q_n(q^{-x}; \alpha, p, N|q) = K_n^{qtm}(q^{-x}; p, N; q).$$

## 4.19 $q$-Hahn $\to$ $q$-Krawtchouk

If we set $\beta = -\alpha^{-1}q^{-1}p$ in the definition (3.6.1) of the $q$-Hahn polynomials and then let $\alpha \to 0$ we obtain the $q$-Krawtchouk polynomials defined by (3.15.1) :

$$\lim_{\alpha \to 0} Q_n \left( q^{-x}; \alpha, -\frac{p}{\alpha q}, N \bigg| q \right) = K_n(q^{-x}; p, N; q).$$

## 4.20 $q$-Hahn $\to$ Affine $q$-Krawtchouk

The affine $q$-Krawtchouk polynomials defined by (3.16.1) can be obtained from the $q$-Hahn polynomials by the substitution $\alpha = p$ and $\beta = 0$ in (3.6.1) :

$$Q_n(q^{-x}; p, 0, N|q) = K_n^{Aff}(q^{-x}; p, N; q).$$

## 4.21 Dual $q$-Hahn $\to$ Affine $q$-Krawtchouk

The affine $q$-Krawtchouk polynomials defined by (3.16.1) can be obtained from the dual $q$-Hahn polynomials by the substitution $\gamma = p$ and $\delta = 0$ in (3.7.1) :

$$R_n(\mu(x); p, 0, N|q) = K_n^{Aff}(q^{-x}; p, N; q).$$

Note that $\mu(x) = q^{-x}$ in this case.

## 4.22 Dual $q$-Hahn $\to$ Dual $q$-Krawtchouk

The dual $q$-Krawtchouk polynomials defined by (3.17.1) can be obtained from the dual $q$-Hahn polynomials by setting $\delta = c\gamma^{-1}q^{-N-1}$ in (3.7.1) and then letting $\gamma \to 0$ :

$$\lim_{\gamma \to 0} R_n \left( \mu(x); \gamma, \frac{c}{\gamma}q^{-N-1} \bigg| q \right) = K_n(\lambda(x); c, N|q).$$

## 4.23 Al-Salam-Chihara $\to$ Continuous big $q$-Hermite

If we take the limit $b \to 0$ in the definition (3.8.1) of the Al-Salam-Chihara polynomials we simply obtain the continuous big $q$-Hermite polynomials given by (3.18.1) :

$$\lim_{b \to 0} Q_n(x; a, b|q) = H_n(x; a|q).$$

## 4.24 Al-Salam-Chihara $\to$ Continuous $q$-Laguerre

The continuous $q$-Laguerre polynomials defined by (3.19.1) can be obtained from the Al-Salam-Chihara polynomials given by (3.8.1) by taking $a = q^{\frac{1}{2}\alpha + \frac{1}{4}}$ and $b = q^{\frac{1}{2}\alpha + \frac{3}{4}}$ :

$$Q_n \left( x; q^{\frac{1}{2}\alpha + \frac{1}{4}}, q^{\frac{1}{2}\alpha + \frac{3}{4}} \bigg| q \right) = \frac{(q;q)_n}{q^{(\frac{1}{2}\alpha + \frac{1}{4})n}} P_n^{(\alpha)}(x|q).$$

## 4.25 $q$-Meixner-Pollaczek $\to$ Continuous $q$-ultraspherical / Rogers

If we take $\theta = 0$ and $a = \beta$ in the definition (3.9.1) of the $q$-Meixner-Pollaczek polynomials we obtain the continuous $q$-ultraspherical (or Rogers) polynomials given by (3.10.15) :

$$P_n(\cos\phi; \beta|q) = C_n(\cos\phi; \beta|q).$$

## 4.26 Continuous $q$-Jacobi $\to$ Continuous $q$-Laguerre

The continuous $q$-Laguerre polynomials given by (3.19.1) and (3.19.2) follow simply from the continuous $q$-Jacobi polynomials defined by (3.10.1) and (3.10.2) respectively by taking the limit $\beta \to \infty$ :

$$\lim_{\beta \to \infty} P_n^{(\alpha,\beta)}(x|q) = P_n^{(\alpha)}(x|q)$$

and

$$\lim_{\beta \to \infty} P_n^{(\alpha,\beta)}(x;q) = \frac{P_n^{(\alpha)}(x;q)}{(-q;q)_n}.$$

## 4.27 Continuous $q$-ultraspherical / Rogers $\to$ Continuous $q$-Hermite

The continuous $q$-Hermite polynomials defined by (3.26.1) can be obtained from the continuous $q$-ultraspherical (or Rogers) polynomials given by (3.10.15) by taking the limit $\beta \to 0$. In fact we have

$$\lim_{\beta \to 0} C_n(x;\beta|q) = \frac{H_n(x|q)}{(q;q)_n}.$$

## 4.28 Big $q$-Laguerre $\to$ Little $q$-Laguerre / Wall

The little $q$-Laguerre (or Wall) polynomials defined by (3.20.1) can be obtained from the big $q$-Laguerre polynomials by taking $x \to bqx$ in (3.11.1) and then letting $b \to \infty$ :

$$\lim_{b \to \infty} P_n(bqx;a,b;q) = p_n(x;a|q).$$

## 4.29 Big $q$-Laguerre $\to$ Al-Salam-Carlitz I

If we set $x \to aqx$ and $b \to ab$ in the definition (3.11.1) of the big $q$-Laguerre polynomials and take the limit $a \to 0$ we obtain the Al-Salam-Carlitz I polynomials given by (3.24.1) :

$$\lim_{a \to 0} \frac{P_n(aqx;a,ab;q)}{a^n} = U_n^{(b)}(x;q).$$

## 4.30 Little $q$-Jacobi $\to$ Little $q$-Laguerre / Wall

The little $q$-Laguerre (or Wall) polynomials defined by (3.20.1) are little $q$-Jacobi polynomials with $b = 0$. So if we set $b = 0$ in the definition (3.12.1) of the little $q$-Jacobi polynomials we obtain the little $q$-Laguerre (or Wall) polynomials :

$$p_n(x;a,0|q) = p_n(x;a|q).$$

## 4.31 Little $q$-Jacobi $\to$ $q$-Laguerre

If we substitute $a = q^\alpha$ and $x \to -b^{-1}q^{-1}x$ in the definition (3.12.1) of the little $q$-Jacobi polynomials and then let $b$ tend to infinity we find the $q$-Laguerre polynomials given by (3.21.1) :

$$\lim_{b \to \infty} p_n\left(-\frac{x}{bq};q^\alpha,b \middle| q\right) = \frac{(q;q)_n}{(q^{\alpha+1};q)_n} L_n^{(\alpha)}(x;q).$$

## 4.32 Little $q$-Jacobi $\to$ Alternative $q$-Charlier

If we set $b \to -a^{-1}q^{-1}b$ in the definition (3.12.1) of the little $q$-Jacobi polynomials and then take the limit $a \to 0$ we obtain the alternative $q$-Charlier polynomials given by (3.22.1) :

$$\lim_{a \to 0} p_n\left(x; a, -\frac{b}{aq} \Big| q\right) = K_n(x; b; q).$$

## 4.33 $q$-Meixner $\to$ $q$-Laguerre

The $q$-Laguerre polynomials defined by (3.21.1) can be obtained from the $q$-Meixner polynomials given by (3.13.1) by setting $b = q^\alpha$ and $q^{-x} \to cq^\alpha x$ in the definition (3.13.1) of the $q$-Meixner polynomials and then taking the limit $c \to \infty$ :

$$\lim_{c \to \infty} M_n(cq^\alpha x; q^\alpha, c; q) = \frac{(q;q)_n}{(q^{\alpha+1};q)_n} L_n^{(\alpha)}(x;q).$$

## 4.34 $q$-Meixner $\to$ $q$-Charlier

The $q$-Meixner polynomials and the $q$-Charlier polynomials defined by (3.13.1) and (3.23.1) respectively are simply related by the limit $b \to 0$ in the definition (3.13.1) of the $q$-Meixner polynomials. In fact we have

$$M_n(x; 0, a; q) = C_n(x; a; q).$$

## 4.35 $q$-Meixner $\to$ Al-Salam-Carlitz II

The Al-Salam-Carlitz II polynomials defined by (3.25.1) can be obtained from the $q$-Meixner polynomials defined by (3.13.1) by setting $b = -c^{-1}a$ in the definition (3.13.1) of the $q$-Meixner polynomials and then taking the limit $c \downarrow 0$ :

$$\lim_{c \downarrow 0} M_n\left(x; -\frac{a}{c}, c; q\right) = \left(-\frac{1}{a}\right)^n q^{\binom{n}{2}} V_n^{(a)}(x; q).$$

## 4.36 Quantum $q$-Krawtchouk $\to$ Al-Salam-Carlitz II

If we set $p = a^{-1}q^{-N-1}$ in the definition (3.14.1) of the quantum $q$-Krawtchouk polynomials and let $N \to \infty$ we obtain the Al-Salam-Carlitz II polynomials given by (3.25.1). In fact we have

$$\lim_{N \to \infty} K_n^{qtm}(x; a^{-1}q^{-N-1}, N; q) = \left(-\frac{1}{a}\right)^n q^{\binom{n}{2}} V_n^{(a)}(x; q).$$

## 4.37 $q$-Krawtchouk $\to$ Alternative $q$-Charlier

If we set $x \to N - x$ in the definition (3.15.1) of the $q$-Krawtchouk polynomials and then take the limit $N \to \infty$ we obtain the alternative $q$-Charlier polynomials defined by (3.22.1) :

$$\lim_{N \to \infty} K_n\left(q^{x-N}; p, N; q\right) = K_n(q^x; p; q).$$

## 4.38 $q$-Krawtchouk $\to$ $q$-Charlier

The $q$-Charlier polynomials given by (3.23.1) can be obtained from the $q$-Krawtchouk polynomials defined by (3.15.1) by setting $p = a^{-1}q^{-N}$ in the definition (3.15.1) of the $q$-Krawtchouk polynomials and then taking the limit $N \to \infty$ :

$$\lim_{N \to \infty} K_n\left(q^{-x}; a^{-1}q^{-N}, N; q\right) = C_n(q^{-x}; a; q).$$

## 4.39 Affine $q$-Krawtchouk $\to$ Little $q$-Laguerre / Wall

If we set $x \to N - x$ in the definition (3.16.1) of the affine $q$-Krawtchouk polynomials and take the limit $N \to \infty$ we simply obtain the little $q$-Laguerre (or Wall) polynomials defined by (3.20.1) :

$$\lim_{N \to \infty} K_n^{Aff}(q^{x-N}; p, N; q) = p_n(q^x; p; q).$$

## 4.40 Dual $q$-Krawtchouk $\to$ Al-Salam-Carlitz I

If we set $c = a^{-1}$ in the definition (3.17.1) of the dual $q$-Krawtchouk polynomials and take the limit $N \to \infty$ we simply obtain the Al-Salam-Carlitz I polynomials given by (3.24.1) :

$$\lim_{N \to \infty} K_n\left(\lambda(x); \frac{1}{a}, N \middle| q\right) = \left(-\frac{1}{a}\right)^n q^{-\binom{n}{2}} U_n^{(a)}(q^x; q).$$

Note that $\lambda(x) = q^{-x} + a^{-1} q^{x-N}$.

## 4.41 Continuous big $q$-Hermite $\to$ Continuous $q$-Hermite

The continuous $q$-Hermite polynomials defined by (3.26.1) can easily be obtained from the continuous big $q$-Hermite polynomials given by (3.18.1) by taking $a = 0$ :

$$H_n(x; 0|q) = H_n(x|q).$$

## 4.42 Continuous $q$-Laguerre $\to$ Continuous $q$-Hermite

The continuous $q$-Hermite polynomials given by (3.26.1) can be obtained from the continuous $q$-Laguerre polynomials defined by (3.19.1) by taking the limit $\alpha \to \infty$ in the following way :

$$\lim_{\alpha \to \infty} \frac{P_n^{(\alpha)}(x|q)}{q^{(\frac{1}{2}\alpha + \frac{1}{4})n}} = \frac{H_n(x|q)}{(q; q)_n}.$$

## 4.43 $q$-Laguerre $\to$ Stieltjes-Wigert

If we set $x \to xq^{-\alpha}$ in the definition (3.21.1) of the $q$-Laguerre polynomials and take the limit $\alpha \to \infty$ we simply obtain the Stieltjes-Wigert polynomials given by (3.27.1) :

$$\lim_{\alpha \to \infty} L_n^{(\alpha)}\left(xq^{-\alpha}; q\right) = S_n(x; q).$$

## 4.44 Alternative $q$-Charlier $\to$ Stieltjes-Wigert

The Stieltjes-Wigert polynomials defined by (3.27.1) can be obtained from the alternative $q$-Charlier polynomials by setting $x \to a^{-1}x$ in the definition (3.22.1) of the alternative $q$-Charlier polynomials and then taking the limit $a \to \infty$. In fact we have

$$\lim_{a \to \infty} K_n\left(\frac{x}{a}; a; q\right) = (q; q)_n S_n(x; q).$$

## 4.45 $q$-Charlier $\to$ Stieltjes-Wigert

If we set $q^{-x} \to ax$ in the definition (3.23.1) of the $q$-Charlier polynomials and take the limit $a \to \infty$ we obtain the Stieltjes-Wigert polynomials given by (3.27.1) in the following way :

$$\lim_{a \to \infty} C_n(ax; a; q) = (q; q)_n S_n(x; q).$$

## 4.46 Al-Salam-Carlitz I → Discrete $q$-Hermite I

The discrete $q$-Hermite I polynomials defined by (3.28.1) can easily be obtained from the Al-Salam-Carlitz I polynomials given by (3.24.1) by the substitution $a = -1$ :

$$U_n^{(-1)}(x;q) = h_n(x;q).$$

## 4.47 Al-Salam-Carlitz II → Discrete $q$-Hermite II

The discrete $q$-Hermite II polynomials defined by (3.29.1) follow from the Al-Salam-Carlitz II polynomials given by (3.25.1) by the substitution $a = -1$ in the following way :

$$i^{-n}V_n^{(-1)}(ix;q) = \tilde{h}_n(x;q).$$

# Chapter 5

# From basic to classical hypergeometric orthogonal polynomials

## 5.1 Askey-Wilson → Wilson

To find the Wilson polynomials defined by (1.1.1) from the Askey-Wilson polynomials we set $a \to q^a$, $b \to q^b$, $c \to q^c$, $d \to q^d$ and $e^{i\theta} = q^{ix}$ (or $\theta = \ln q^x$) in the definition (3.1.1) and take the limit $q \uparrow 1$ :

$$\lim_{q\uparrow 1} \frac{p_n(\frac{1}{2}\left(q^{ix}+q^{-ix}\right);q^a,q^b,q^c,q^d|q)}{(1-q)^{3n}} = W_n(x^2;a,b,c,d).$$

## 5.2 $q$-Racah → Racah

If we set $\alpha \to q^\alpha$, $\beta \to q^\beta$, $\gamma \to q^\gamma$, $\delta \to q^\delta$ in the definition (3.2.1) of the $q$-Racah polynomials and let $q \uparrow 1$ we easily obtain the Racah polynomials defined by (1.2.1) :

$$\lim_{q\uparrow 1} R_n(\mu(x);q^\alpha,q^\beta,q^\gamma,q^\delta|q) = R_n(\lambda(x);\alpha,\beta,\gamma,\delta),$$

where

$$\begin{cases} \mu(x) = q^{-x} + q^{x+\gamma+\delta+1} \\ \lambda(x) = x(x+\gamma+\delta+1). \end{cases}$$

## 5.3 Continuous dual $q$-Hahn → Continuous dual Hahn

To find the continuous dual Hahn polynomials defined by (1.3.1) from the continuous dual $q$-Hahn polynomials we set $a \to q^a$, $b \to q^b$, $c \to q^c$ and $e^{i\theta} = q^{ix}$ (or $\theta = \ln q^x$) in the definition (3.3.1) and take the limit $q \uparrow 1$ :

$$\lim_{q\uparrow 1} \frac{p_n(\frac{1}{2}\left(q^{ix}+q^{-ix}\right);q^a,q^b,q^c|q)}{(1-q)^{2n}} = S_n(x^2;a,b,c).$$

## 5.4 Continuous $q$-Hahn → Continuous Hahn

If we set $a \to q^a$, $b \to q^b$, $c \to q^c$, $d \to q^d$ and $e^{-i\theta} = q^{ix}$ (or $\theta = \ln q^{-x}$) in the definition (3.4.1) of the continuous $q$-Hahn polynomials and take the limit $q \uparrow 1$ we find the continuous Hahn polynomials given by (1.4.1) in the following way :

$$\lim_{q\uparrow 1}\frac{p_n(\cos(\ln q^{-x}+\phi);q^a,q^b,q^c,q^d;q)}{(1-q)^n(q;q)_n}=(-2\sin\phi)^n p_n(x;a,b,c,d).$$

## 5.5 Big $q$-Jacobi $\to$ Jacobi

If we set $c=0$, $a=q^\alpha$ and $b=q^\beta$ in the definition (3.5.1) of the big $q$-Jacobi polynomials and let $q\uparrow 1$ we find the Jacobi polynomials given by (1.8.1) :

$$\lim_{q\uparrow 1}P_n(x;q^\alpha,q^\beta,0;q)=\frac{P_n^{(\alpha,\beta)}(2x-1)}{P_n^{(\alpha,\beta)}(1)}.$$

If we take $c=-q^\gamma$ for arbitrary real $\gamma$ instead of $c=0$ we find

$$\lim_{q\uparrow 1}P_n(x;q^\alpha,q^\beta,-q^\gamma;q)=\frac{P_n^{(\alpha,\beta)}(x)}{P_n^{(\alpha,\beta)}(1)}.$$

### 5.5.1 Big $q$-Legendre $\to$ Legendre / Spherical

If we set $c=0$ in the definition (3.5.7) of the big $q$-Legendre polynomials and let $q\uparrow 1$ we simply obtain the Legendre (or spherical) polynomials defined by (1.8.40) :

$$\lim_{q\uparrow 1}P_n(x;0;q)=P_n(2x-1).$$

If we take $c=-q^\gamma$ for arbitrary real $\gamma$ instead of $c=0$ we find

$$\lim_{q\uparrow 1}P_n(x;-q^\gamma;q)=P_n(x).$$

## 5.6 $q$-Hahn $\to$ Hahn

The Hahn polynomials defined by (1.5.1) simply follow from the $q$-Hahn polynomials given by (3.6.1), after setting $\alpha\to q^\alpha$ and $\beta\to q^\beta$, in the following way :

$$\lim_{q\uparrow 1}Q_n(q^{-x};q^\alpha,q^\beta,N|q)=Q_n(x;\alpha,\beta,N).$$

## 5.7 Dual $q$-Hahn $\to$ Dual Hahn

The dual Hahn polynomials given by (1.6.1) follow from the dual $q$-Hahn polynomials by simply taking the limit $q\uparrow 1$ in the definition (3.7.1) of the dual $q$-Hahn polynomials after applying the substitution $\gamma\to q^\gamma$ and $\delta\to q^\delta$ :

$$\lim_{q\uparrow 1}R_n(\mu(x);q^\gamma,q^\delta,N|q)=R_n(\lambda(x);\gamma,\delta,N),$$

where

$$\begin{cases}\mu(x)=q^{-x}+q^{x+\gamma+\delta+1}\\ \lambda(x)=x(x+\gamma+\delta+1).\end{cases}$$

## 5.8 Al-Salam-Chihara $\to$ Meixner-Pollaczek

If we set $a=q^\lambda e^{-i\phi}$, $b=q^\lambda e^{i\phi}$ and $e^{i\theta}=q^{ix}e^{i\phi}$ in the definition (3.8.1) of the Al-Salam-Chihara polynomials and take the limit $q\uparrow 1$ we obtain the Meixner-Pollaczek polynomials given by (1.7.1) in the following way :

$$\lim_{q\uparrow 1}\frac{Q_n\left(\cos(\ln q^x+\phi);q^\lambda e^{i\phi},q^\lambda e^{-i\phi}|q\right)}{(q;q)_n}=P_n^{(\lambda)}(x;\phi).$$

## 5.9 $q$-Meixner-Pollaczek $\rightarrow$ Meixner-Pollaczek

To find the Meixner-Pollaczek polynomials defined by (1.7.1) from the $q$-Meixner-Pollaczek polynomials we substitute $a = q^\lambda$ and $e^{i\theta} = q^{-ix}$ (or $\theta = \ln q^{-x}$) in the definition (3.9.1) of the $q$-Meixner-Pollaczek polynomials and take the limit $q \uparrow 1$ to find :

$$\lim_{q \uparrow 1} P_n(\cos(\ln q^{-x} + \phi); q^\lambda | q) = P_n^{(\lambda)}(x; -\phi).$$

## 5.10 Continuous $q$-Jacobi $\rightarrow$ Jacobi

If we take the limit $q \uparrow 1$ in the definitions (3.10.1) and (3.10.2) of the continuous $q$-Jacobi polynomials we simply find the Jacobi polynomials defined by (1.8.1) :

$$\lim_{q \uparrow 1} P_n^{(\alpha,\beta)}(x|q) = P_n^{(\alpha,\beta)}(x)$$

and

$$\lim_{q \uparrow 1} P_n^{(\alpha,\beta)}(x;q) = P_n^{(\alpha,\beta)}(x).$$

### 5.10.1 Continuous $q$-ultraspherical / Rogers $\rightarrow$ Gegenbauer / Ultraspherical

If we set $\beta = q^\lambda$ in the definition (3.10.15) of the continuous $q$-ultraspherical (or Rogers) polynomials and let $q$ tend to one we obtain the Gegenbauer (or ultraspherical) polynomials given by (1.8.10) :

$$\lim_{q \uparrow 1} C_n(x; q^\lambda | q) = C_n^{(\lambda)}(x).$$

### 5.10.2 Continuous $q$-Legendre $\rightarrow$ Legendre / Spherical

The Legendre (or spherical) polynomials defined by (1.8.40) easily follow from the continuous $q$-Legendre polynomials given by (3.10.25) by taking the limit $q \uparrow 1$ :

$$\lim_{q \uparrow 1} P_n(x; q) = P_n(x).$$

Of course, we also have

$$\lim_{q \uparrow 1} P_n(x|q) = P_n(x).$$

## 5.11 Big $q$-Laguerre $\rightarrow$ Laguerre

The Laguerre polynomials defined by (1.11.1) can be obtained from the big $q$-Laguerre polynomials by the substitution $a = q^\alpha$ and $b = (1-q)^{-1} q^\beta$ in the definition (3.11.1) of the big $q$-Laguerre polynomials and the limit $q \uparrow 1$ :

$$\lim_{q \uparrow 1} P_n(x; q^\alpha, (1-q)^{-1} q^\beta; q) = \frac{L_n^{(\alpha)}(x-1)}{L_n^{(\alpha)}(0)}.$$

## 5.12 Little $q$-Jacobi $\rightarrow$ Jacobi

The Jacobi polynomials defined by (1.8.1) simply follow from the little $q$-Jacobi polynomials defined by (3.12.1) in the following way :

$$\lim_{q \uparrow 1} p_n(x; q^\alpha, q^\beta | q) = \frac{P_n^{(\alpha,\beta)}(1-2x)}{P_n^{(\alpha,\beta)}(1)}.$$

### 5.12.1 Little $q$-Legendre $\to$ Legendre / Spherical

If we take the limit $q \uparrow 1$ in the definition (3.12.6) of the little $q$-Legendre polynomials we simply find the Legendre (or spherical) polynomials given by (1.8.40) :

$$\lim_{q \uparrow 1} p_n(x|q) = P_n(1 - 2x).$$

### 5.12.2 Little $q$-Jacobi $\to$ Laguerre

If we take $a = q^\alpha$, $b = -q^\beta$ for arbitrary real $\beta$ and $x \to \frac{1}{2}(1-q)x$ in the definition (3.12.1) of the little $q$-Jacobi polynomials and then take the limit $q \uparrow 1$ we obtain the Laguerre polynomials given by (1.11.1) :

$$\lim_{q \uparrow 1} p_n \left( \frac{1}{2}(1-q)x; q^\alpha, -q^\beta \middle| q \right) = \frac{L_n^{(\alpha)}(x)}{L_n^{(\alpha)}(0)}.$$

## 5.13 $q$-Meixner $\to$ Meixner

To find the Meixner polynomials defined by (1.9.1) from the $q$-Meixner polynomials given by (3.13.1) we set $b = q^{\beta-1}$ and $c \to (1-c)^{-1}c$ and let $q \uparrow 1$ :

$$\lim_{q \uparrow 1} M_n \left( q^{-x}; q^{\beta-1}, \frac{c}{1-c}; q \right) = M_n(x; \beta, c).$$

## 5.14 Quantum $q$-Krawtchouk $\to$ Krawtchouk

The Krawtchouk polynomials given by (1.10.1) easily follow from the quantum $q$-Krawtchouk polynomials defined by (3.14.1) in the following way :

$$\lim_{q \uparrow 1} K_n^{qtm}(q^{-x}; p, N; q) = K_n(x; p^{-1}, N).$$

## 5.15 $q$-Krawtchouk $\to$ Krawtchouk

If we take the limit $q \uparrow 1$ in the definition (3.15.1) of the $q$-Krawtchouk polynomials we simply find the Krawtchouk polynomials given by (1.10.1) in the following way :

$$\lim_{q \uparrow 1} K_n(q^{-x}; p, N; q) = K_n \left( x; \frac{1}{p+1}, N \right).$$

## 5.16 Affine $q$-Krawtchouk $\to$ Krawtchouk

If we let $q \uparrow 1$ in the definition (3.16.1) of the affine $q$-Krawtchouk polynomials we obtain :

$$\lim_{q \uparrow 1} K_n^{Aff}(q^{-x}; p, N|q) = K_n(x; 1-p, N),$$

where $K_n(x; 1-p, N)$ is the Krawtchouk polynomial defined by (1.10.1).

## 5.17 Dual $q$-Krawtchouk $\to$ Krawtchouk

If we set $c = 1 - p^{-1}$ in the definition (3.17.1) of the dual $q$-Krawtchouk polynomials and take the limit $q \uparrow 1$ we simply find the Krawtchouk polynomials given by (1.10.1) :

$$\lim_{q \uparrow 1} K_n \left( \lambda(x); 1 - \frac{1}{p}, N | q \right) = K_n(x; p, N).$$

## 5.18 Continuous big $q$-Hermite → Hermite

If we set $a = 0$ and $x \to x\sqrt{\frac{1}{2}(1-q)}$ in the definition (3.18.1) of the continuous big $q$-Hermite polynomials and let $q$ tend to one, we obtain the Hermite polynomials given by (1.13.1) in the following way :

$$\lim_{q \uparrow 1} \frac{H_n\left(x\left(\frac{1-q}{2}\right)^{\frac{1}{2}}; 0 \Big| q\right)}{\left(\frac{1-q}{2}\right)^{\frac{n}{2}}} = H_n(x).$$

If we take $a \to a\sqrt{2(1-q)}$ and $x \to x\sqrt{\frac{1}{2}(1-q)}$ in the definition (3.18.1) of the continuous big $q$-Hermite polynomials and take the limit $q \uparrow 1$ we find the Hermite polynomials defined by (1.13.1) with shifted argument :

$$\lim_{q \uparrow 1} \frac{H_n\left(x\left(\frac{1-q}{2}\right)^{\frac{1}{2}}; a\sqrt{2(1-q)} \Big| q\right)}{\left(\frac{1-q}{2}\right)^{\frac{n}{2}}} = H_n(x-a).$$

## 5.19 Continuous $q$-Laguerre → Laguerre

If we set $x \to q^x$ in the definitions (3.19.1) and (3.19.2) of the continuous $q$-Laguerre polynomials and take the limit $q \uparrow 1$ we find the Laguerre polynomials defined by (1.11.1). In fact we have :

$$\lim_{q \uparrow 1} P_n^{(\alpha)}(q^x|q) = L_n^{(\alpha)}(2x)$$

and

$$\lim_{q \uparrow 1} P_n^{(\alpha)}(q^x; q) = L_n^{(\alpha)}(x).$$

## 5.20 Little $q$-Laguerre / Wall → Laguerre

If we set $a = q^\alpha$ and $x \to (1-q)x$ in the definition (3.20.1) of the little $q$-Laguerre (or Wall) polynomials and let $q$ tend to one, we obtain the Laguerre polynomials given by (1.11.1) :

$$\lim_{q \uparrow 1} p_n((1-q)x; q^\alpha|q) = \frac{L_n^{(\alpha)}(x)}{L_n^{(\alpha)}(0)}.$$

### 5.20.1 Little $q$-Laguerre / Wall → Charlier

If we set $a \to (q-1)a$ and $x \to q^x$ in the definition (3.20.1) of the little $q$-Laguerre (or Wall) polynomials and take the limit $q \uparrow 1$ we obtain the Charlier polynomials given by (1.12.1) in the following way :

$$\lim_{q \uparrow 1} \frac{p_n(q^x; (q-1)a|q)}{(1-q)^n} = \frac{C_n(x; a)}{a^n}.$$

## 5.21 $q$-Laguerre → Laguerre

If we set $x \to (1-q)x$ in the definition (3.21.1) of the $q$-Laguerre polynomials and take the limit $q \uparrow 1$ we obtain the Laguerre polynomials given by (1.11.1) :

$$\lim_{q \uparrow 1} L_n^{(\alpha)}((1-q)x; q) = L_n^{(\alpha)}(x).$$

### 5.21.1 $q$-Laguerre $\to$ Charlier

If we set $x \to -q^{-x}$ and $q^\alpha = a^{-1}(q-1)^{-1}$ (or $\alpha = -(\ln q)^{-1}\ln(q-1)a$) in the definition (3.21.1) of the $q$-Laguerre polynomials, multiply by $(q;q)_n$, and take the limit $q \uparrow 1$ we obtain the Charlier polynomials given by (1.12.1) :

$$\lim_{q \uparrow 1}(q;q)_n L_n^{(\alpha)}(-q^{-x};q) = C_n(x;a),\ q^\alpha = \frac{1}{a(q-1)}\ \text{ or }\ \alpha = -\frac{\ln(q-1)a}{\ln q}.$$

## 5.22 Alternative $q$-Charlier $\to$ Charlier

If we set $x \to q^x$ and $a \to a(1-q)$ in the definition (3.22.1) of the alternative $q$-Charlier polynomials and take the limit $q \uparrow 1$ we find the Charlier polynomials given by (1.12.1) :

$$\lim_{q \uparrow 1}\frac{K_n(q^x;a(1-q);q)}{(q-1)^n} = a^n C_n(x;a).$$

## 5.23 $q$-Charlier $\to$ Charlier

If we set $a \to a(1-q)$ in the definition (3.23.1) of the $q$-Charlier polynomials and take the limit $q \uparrow 1$ we obtain the Charlier polynomials defined by (1.12.1) :

$$\lim_{q \uparrow 1} C_n(q^{-x};a(1-q);q) = C_n(x;a).$$

## 5.24 Al-Salam-Carlitz I $\to$ Charlier

If we set $a \to a(q-1)$ and $x \to q^x$ in the definition (3.24.1) of the Al-Salam-Carlitz I polynomials and take the limit $q \uparrow 1$ after dividing by $a^n(1-q)^n$ we obtain the Charlier polynomials defined by (1.12.1) :

$$\lim_{q \uparrow 1}\frac{U_n^{(a(q-1))}(q^x;q)}{(1-q)^n} = a^n C_n(x;a).$$

### 5.24.1 Al-Salam-Carlitz I $\to$ Hermite

If we set $x \to x\sqrt{1-q^2}$ and $a \to a\sqrt{1-q^2}-1$ in the definition (3.24.1) of the Al-Salam-Carlitz I polynomials, divide by $(1-q^2)^{\frac{n}{2}}$, and let $q$ tend to one we obtain the Hermite polynomials given by (1.13.1) with shifted argument. In fact we have

$$\lim_{q \uparrow 1}\frac{U_n^{(a\sqrt{1-q^2}-1)}(x\sqrt{1-q^2};q)}{(1-q^2)^{\frac{n}{2}}} = \frac{H_n(x-a)}{2^n}.$$

## 5.25 Al-Salam-Carlitz II $\to$ Charlier

If we set $a \to a(1-q)$ and $x \to q^{-x}$ in the definition (3.25.1) of the Al-Salam-Carlitz II polynomials and taking the limit $q \uparrow 1$ we find

$$\lim_{q \uparrow 1}\frac{V_n^{(a(1-q))}(q^{-x};q)}{(q-1)^n} = a^n C_n(x;a).$$

### 5.25.1 Al-Salam-Carlitz II → Hermite

If we set $x \to x\sqrt{1-q^2}$ and $a \to a\sqrt{1-q^2}+1$ in the definition (3.25.1) of the Al-Salam-Carlitz II polynomials, divide by $(1-q^2)^{\frac{n}{2}}$, and let $q$ tend to one we obtain the Hermite polynomials given by (1.13.1) with shifted argument. In fact we have

$$\lim_{q \uparrow 1} \frac{V_n^{(a\sqrt{1-q^2}+1)}(x\sqrt{1-q^2};q)}{(1-q^2)^{\frac{n}{2}}} = \frac{H_n(x-2)}{2^n}.$$

## 5.26 Continuous $q$-Hermite → Hermite

The Hermite polynomials defined by (1.13.1) can be obtained from the continuous $q$-Hermite polynomials given by (3.26.1) by setting $x \to x\sqrt{\frac{1}{2}(1-q)}$. In fact we have

$$\lim_{q \uparrow 1} \frac{H_n\left(x\left(\frac{1-q}{2}\right)^{\frac{1}{2}}\Big| q\right)}{\left(\frac{1-q}{2}\right)^{\frac{n}{2}}} = H_n(x).$$

## 5.27 Stieltjes-Wigert → Hermite

The Hermite polynomials defined by (1.13.1) can be obtained from the Stieltjes-Wigert polynomials given by (3.27.1) by setting $x \to q^{-1}x\sqrt{2(1-q)}+1$ and taking the limit $q \uparrow 1$ in the following way :

$$\lim_{q \uparrow 1} \frac{(q;q)_n S_n(q^{-1}x\sqrt{2(1-q)}+1;q)}{\left(\frac{1-q}{2}\right)^{\frac{n}{2}}} = (-1)^n H_n(x).$$

## 5.28 Discrete $q$-Hermite I → Hermite

The Hermite polynomials defined by (1.13.1) can be found from the discrete $q$-Hermite I polynomials given by (3.28.1) in the following way :

$$\lim_{q \uparrow 1} \frac{h_n\left(x\sqrt{1-q^2};q\right)}{(1-q^2)^{\frac{n}{2}}} = \frac{H_n(x)}{2^n}.$$

## 5.29 Discrete $q$-Hermite II → Hermite

The Hermite polynomials defined by (1.13.1) can also be found from the discrete $q$-Hermite II polynomials given by (3.29.1) in a similar way :

$$\lim_{q \uparrow 1} \frac{\tilde{h}_n\left(x\sqrt{1-q^2};q\right)}{(1-q^2)^{\frac{n}{2}}} = \frac{H_n(x)}{2^n}.$$

# Index

# The Askey-scheme of hypergeometric orthogonal polynomials and its $q$-analogue

Roelof Koekoek    René F. Swarttouw

February 20, 1996

**List of errata in report no. 94-05**

- Page 7, line 1. Replace "all sets" by "all known sets".

- Page 7, line -1. Replace "positive definite" by "positive".

- Page 12, line -2. This should read

$$\lim_{q \uparrow 1} {}_r\phi_s \left( \begin{array}{c} q^{a_1}, \ldots, q^{a_r} \\ q^{b_1}, \ldots, q^{b_s} \end{array} \bigg| q; (q-1)^{1+s-r} z \right) = {}_rF_s \left( \begin{array}{c} a_1, \ldots, a_r \\ b_1, \ldots, b_s \end{array} \bigg| z \right).$$

- Page 28, formula (1.4.2). The right-hand side should read

$$\frac{\Gamma(n+a+c)\Gamma(n+a+d)\Gamma(n+b+c)\Gamma(n+b+d)}{(2n+a+b+c+d-1)\Gamma(n+a+b+c+d-1)n!} \delta_{mn}.$$

- Page 51, formula (3.1.5). This should read

$$(1-q)^2 D_q \left[ \tilde{w}(x; aq^{\frac{1}{2}}, bq^{\frac{1}{2}}, cq^{\frac{1}{2}}, dq^{\frac{1}{2}}|q) D_q y(x) \right] +$$
$$+ \lambda_n \tilde{w}(x; a, b, c, d|q) y(x) = 0, \ y(x) = p_n(x; a, b, c, d|q),$$

where

$$\tilde{w}(x; a, b, c, d|q) := \frac{w(x; a, b, c, d|q)}{\sqrt{1-x^2}}.$$

- Page 55, formula (3.3.5). This should read

$$(1-q)^2 D_q \left[ \tilde{w}(x; aq^{\frac{1}{2}}, bq^{\frac{1}{2}}, cq^{\frac{1}{2}}|q) D_q y(x) \right] +$$
$$+ 4q^{-n+1}(1-q^n) \tilde{w}(x; a, b, c|q) y(x) = 0, \ y(x) = p_n(x; a, b, c|q),$$

where

$$\tilde{w}(x; a, b, c|q) := \frac{w(x; a, b, c|q)}{\sqrt{1-x^2}}.$$

- Page 56, formula (3.4.4). This should read

$$(1-q)^2 D_q \left[ \tilde{w}(x; aq^{\frac{1}{2}}, bq^{\frac{1}{2}}, cq^{\frac{1}{2}}, dq^{\frac{1}{2}}; q) D_q y(x) \right] +$$
$$+ \lambda_n \tilde{w}(x; a, b, c, d; q) y(x) = 0, \ y(x) = p_n(x; a, b, c, d; q),$$

where

$$\tilde{w}(x; a, b, c, d; q) := \frac{w(x; a, b, c, d; q)}{\sqrt{1-x^2}}.$$

- Page 57, formula (3.5.2). Replace $(-ac)^{-n}$ by $(-acq^2)^n$.

- Page 58, line 9. The weight function should read

$$\frac{(c^{-1}qx, -d^{-1}qx; q)_\infty}{(ac^{-1}qx, -bd^{-1}qx; q)_\infty} d_q x.$$

- Page 58, formula (3.5.8). Replace $(-c)^{-n}$ by $(-cq^2)^n$.

- Page 63, formula (3.8.5). This should read

$$(1-q)^2 D_q \left[ \tilde{w}(x; aq^{\frac{1}{2}}, bq^{\frac{1}{2}} | q) D_q y(x) \right] +$$
$$+ 4q^{-n+1}(1-q^n) \tilde{w}(x; a, b|q) y(x) = 0, \ y(x) = Q_n(x; a, b|q),$$

where

$$\tilde{w}(x; a, b|q) := \frac{w(x; a, b|q)}{\sqrt{1-x^2}}.$$

- Page 64, formula (3.9.4). This should read

$$(1-q)^2 D_q \left[ \tilde{w}(x; aq^{\frac{1}{2}} | q) D_q y(x) \right] + 4q^{-n+1}(1-q^n) \tilde{w}(x; a|q) y(x) = 0, \ y(x) = P_n(x; a|q),$$

where

$$\tilde{w}(x; a|q) := \frac{w(x; a|q)}{\sqrt{1-x^2}}.$$

- Page 66, formula (3.10.7). This should read

$$(1-q)^2 D_q \left[ \tilde{w}(x; q^{\alpha+\frac{1}{2}}, q^{\beta+\frac{1}{2}} | q) D_q y(x) \right] + \lambda_n \tilde{w}(x; q^\alpha, q^\beta | q) y(x) = 0, \ y(x) = P_n^{(\alpha,\beta)}(x|q),$$

where

$$\tilde{w}(x; q^\alpha, q^\beta | q) := \frac{w(x; q^\alpha, q^\beta | q)}{\sqrt{1-x^2}}.$$

- Page 66, fomula (3.10.8). This should read

$$(1-q)^2 D_q \left[ \tilde{w}(x; q^{\alpha+\frac{1}{2}}, q^{\beta+\frac{1}{2}}; q) D_q y(x) \right] + \lambda_n \tilde{w}(x; q^\alpha, q^\beta; q) y(x) = 0, \ y(x) = P_n^{(\alpha,\beta)}(x; q),$$

where

$$\tilde{w}(x; q^\alpha, q^\beta; q) := \frac{w(x; q^\alpha, q^\beta; q)}{\sqrt{1-x^2}}.$$

- Page 68, formula (3.10.18). This should read

$$(1-q)^2 D_q \left[ \tilde{w}(x; \beta q^{\frac{1}{2}} | q) D_q y(x) \right] + \lambda_n \tilde{w}(x; \beta|q) y(x) = 0, \ y(x) = C_n(x; \beta|q),$$

where

$$\tilde{w}(x; \beta|q) := \frac{w(x; \beta|q)}{\sqrt{1-x^2}}.$$

- Page 70, formula (3.10.27). This can be written as

$$2(1-q^{2n+1}) x P_n(x; q) = q^{-\frac{1}{2}}(1-q^{2n+2}) P_{n+1}(x; q) + q^{\frac{1}{2}}(1-q^{2n}) P_{n-1}(x; q).$$

- Page 70, formula (3.10.28). This should read

$$(1-q)^2 D_q\left[\tilde{w}(x;q;q)D_q y(x)\right] + \lambda_n \tilde{w}(x;q^{\frac{1}{2}};q)y(x) = 0, \; y(x) = P_n(x;q),$$

  where

$$\tilde{w}(x;a;q) := \frac{w(x;a;q)}{\sqrt{1-x^2}}.$$

- Page 71, formula (3.11.2). Replace $(-ab)^{-n}$ by $(-abq^2)^n$.

- Page 80, formula (3.18.5). This should read

$$(1-q)^2 D_q\left[\tilde{w}(x;aq^{\frac{1}{2}}|q)D_q y(x)\right] + 4q^{-n+1}(1-q^n)\tilde{w}(x;a|q)y(x) = 0, \; y(x) = H_n(x;a|q),$$

  where

$$\tilde{w}(x;a|q) := \frac{w(x;a|q)}{\sqrt{1-x^2}}.$$

- Page 82, formula (3.19.5). This can also be written as

$$2xP_n^{(\alpha)}(x|q) = q^{-\frac{1}{2}\alpha-\frac{1}{4}}(1-q^{n+1})P_{n+1}^{(\alpha)}(x|q) +$$
$$+ q^{n+\frac{1}{2}\alpha+\frac{1}{4}}(1+q^{\frac{1}{2}})P_n^{(\alpha)}(x|q) + q^{\frac{1}{2}\alpha+\frac{1}{4}}(1-q^{n+\alpha})P_{n-1}^{(\alpha)}(x|q).$$

- Page 82, formula (3.19.6). This can also be written as

$$2xP_n^{(\alpha)}(x;q) = q^{-\frac{1}{2}}(1-q^{2n+2})P_{n+1}^{(\alpha)}(x;q) +$$
$$+ q^{2n+\alpha+\frac{1}{2}}(1+q)P_n^{(\alpha)}(x;q) + q^{\frac{1}{2}}(1-q^{2n+2\alpha})P_{n-1}^{(\alpha)}(x;q).$$

- Page 82, formula (3.19.7). This should read

$$(1-q)^2 D_q\left[\tilde{w}(x;q^{\alpha+\frac{1}{2}}|q)D_q y(x)\right] + 4q^{-n+1}(1-q^n)\tilde{w}(x;q^\alpha|q)y(x) = 0, \; y(x) = P_n^{(\alpha)}(x|q),$$

  where

$$\tilde{w}(x;q^\alpha|q) := \frac{w(x;q^\alpha|q)}{\sqrt{1-x^2}}.$$

- Page 82, formula (3.19.8). This should read

$$(1-q)^2 D_q\left[\tilde{w}(x;q^{\alpha+\frac{1}{2}};q)D_q y(x)\right] + 4q^{-n+1}(1-q^n)\tilde{w}(x;q^\alpha;q)y(x) = 0, \; y(x) = P_n^{(\alpha)}(x;q),$$

  where

$$\tilde{w}(x;q^\alpha;q) := \frac{w(x;q^\alpha;q)}{\sqrt{1-x^2}}.$$

- Page 87, formula (3.24.4). This should read

$$(1-q^n)x^2 y(x) = aq^{n-1}y(qx) - \left[aq^{n-1} + q^n(1-x)(a-x)\right]y(x) +$$
$$+ q^n(1-x)(a-x)y(q^{-1}x), \; y(x) = U_n^{(a)}(x;q).$$

- Page 87, formula (3.25.4). This should read

$$-(1-q^n)x^2 y(x) = (1-x)(a-x)y(qx) - [(1-x)(a-x) + aq]y(x) +$$
$$+ aqy(q^{-1}x), \; y(x) = V_n^{(a)}(x;q).$$

- Page 88, formula (3.26.2). Replace $d\theta$ by $dx$.

- Page 88, formula (3.26.4). This should read
$$(1-q)^2 D_q \left[\tilde{w}(x) D_q y(x)\right] + 4q^{-n+1}(1-q^n)\tilde{w}(x)y(x) = 0, \ y(x) = H_n(x|q),$$
where
$$\tilde{w}(x) := \frac{w(x)}{\sqrt{1-x^2}}.$$

- Page 90, formula (3.28.4). This should read
$$-q^{-n+1}x^2 y(x) = y(qx) - (1+q)y(x) + q(1-x^2)y(q^{-1}x), \ y(x) = h_n(x;q).$$

- Page 90, formula (3.29.1). This should read
$$\begin{aligned}\tilde{h}_n(x;q) = i^{-n} V_n^{(-1)}(ix;q) &= i^{-n} q^{-\binom{n}{2}} {}_2\phi_0\left(\begin{matrix}q^{-n}, ix \\ -\end{matrix} \bigg| q; -q^n\right) \\ &= x^n {}_2\phi_1\left(\begin{matrix}q^{-n}, q^{-n+1} \\ 0\end{matrix} \bigg| q^2; -\frac{q^2}{x^2}\right).\end{aligned}$$

- Page 90, formula (3.29.2). This can be written as
$$\sum_{k=-\infty}^{\infty} \left[\tilde{h}_m(cq^k;q)\tilde{h}_n(cq^k;q) + \tilde{h}_m(-cq^k;q)\tilde{h}_n(-cq^k;q)\right] w(cq^k)q^k$$
$$= 2\frac{(q^2, -c^2q, -c^{-2}q; q^2)_\infty}{(q, -c^2, -c^{-2}q^2; q^2)_\infty} \frac{(q;q)_n}{q^{n^2}}\delta_{mn}, \ c > 0,$$
where
$$w(x) = \frac{1}{(ix;q)_\infty(-ix;q)_\infty} = \frac{1}{(-x^2;q^2)_\infty}.$$

- Page 90, formula (3.29.4). This should read
$$-(1-q^n)x^2 y(x) = (1+x^2)y(qx) - (1+x^2+q)y(x) + qy(q^{-1}x), \ y(x) = \tilde{h}_n(x;q).$$

- Page 109. Reference [40] : "bf A 25" should read **A 25**.

### Acknowledgement

We thank G. Gasper, J. Koekoek, H.T. Koelink, and T.H. Koornwinder for pointing us to some of these errata.4